\documentclass[review,onefignum,onetabnum]{siamart171218}

\usepackage{booktabs}
\usepackage{amsmath}
\usepackage{chngcntr}
\usepackage{amssymb}
\usepackage{float}
\usepackage{cleveref}
\usepackage{algorithmic}
\usepackage{theorem}
\usepackage{overpic}
\usepackage{multirow}

\newtheorem{remark}{Remark}[section]

\title{Scalable Bayesian Physics-Informed Kolmogorov-Arnold Networks}
\author{%
  Zhiwei Gao\thanks{Engineering, Brown University, Providence, RI, USA. 
  Email: \texttt{zhiwei\_gao@brown.edu}.} 
  \and 
  George Em Karniadakis\thanks{Applied Mathematics, Brown University, Providence, RI, USA. 
  Corresponding author. Email: \texttt{george\_karniadakis@brown.edu}.}
}

\nolinenumbers
\begin{document}

\maketitle

\begin{abstract}
Uncertainty quantification (UQ) plays a pivotal role in scientific machine learning, especially when surrogate models are used to approximate complex systems. Although multilayer perceptions (MLPs) are commonly employed as surrogates, they often suffer from overfitting due to their large number of parameters. Kolmogorov-Arnold networks (KANs) offer an alternative solution with fewer parameters. However, gradient-based inference methods, such as Hamiltonian Monte Carlo (HMC), may result in computational inefficiency when applied to KANs, especially for large-scale datasets, due to the high cost of back-propagation.  
To address these challenges, we propose a novel approach, combining the dropout Tikhonov ensemble Kalman inversion (DTEKI) with Chebyshev KANs. This gradient-free method effectively mitigates overfitting and enhances numerical stability. Additionally, we incorporate the active subspace method to reduce the parameter-space dimensionality, allowing us to improve the accuracy of predictions and obtain more reliable uncertainty estimates.  
Extensive experiments demonstrate the efficacy of our approach in various test cases, including scenarios with large datasets and high noise levels. Our results show that the new method achieves comparable or better accuracy, much higher efficiency as well as stability compared to HMC, in addition to scalability. Moreover, by leveraging the low-dimensional parameter subspace, our method preserves prediction accuracy while substantially reducing further the computational cost. 
\end{abstract}

\begin{keywords}
  Uncertainty quantification, Ensemble Kalman inversion, Active Subspace, PINNs, KANs
\end{keywords}

\begin{AMS}
  68Q25, 68R10, 68U05
\end{AMS}

\section{Introduction}
Partial differential equations (PDEs) -- arising in many applications, such as electromagnetism,  heat transfer, and fluid mechanics -- are effective representations describing well-established conservation laws.   The behavior of the solution of PDEs is typically described by their parameters, which cannot always be directly determined. Instead, these parameters can be inferred by partial observations with noise, leading to the emergence of inverse problems, which are generally ill-posed and thus hard to solve. Traditional numerical solvers for inverse problems generally include regularization methods or statistical inference methods,  e.g., Bayesian inference. In practice, by imposing a prior distribution to the parameters, the solution to the Bayesian inverse problems can be obtained by exploring the posterior distribution condition on the observations, which not only gives an accurate estimate of the ground truth but can provide reliable uncertainty quantification that is quite crucial in some real applications.

However, the posterior exploration of Bayesian inverse problems typically runs into the so-called curse-of-dimensionality (CoD) when employing MCMC-based sampling methods. That is; because the sampling process requires repeated evaluations \cite{willard2022integrating} of the expensive forward model governed by the PDEs, leading to prohibitive computational costs. To alleviate these issues, common approaches include (i) reduced-order methods \cite{cui2015data,cui2016scalable,lieberman2010parameter,marzouk2009dimensionality}, which capture the low dimensionality of space, (ii) surrogate modeling \cite{li2014adaptive, yan_adaptive, Yan_2017, MARZOUK2007560, conrad2016accelerating}, and (iii) direct posterior approximation methods \cite{schillings2020convergence,bui2013computational} such as Laplace approximation. While constructing surrogates is the most effective method to be explored,  traditional surrogates such as the Gaussian process \cite{bai2024gaussian}, and polynomial chaos expansion \cite{marzouk2009dimensionality} are difficult to handle high-dimensional problems, and they are hard to encode the physical information, especially for non-linear PDEs. 

The recently emerging physics-informed neural networks (PINNs) \cite{RAISSI2019686, wang2022and, gao2023failure,cuomo2022scientific} provide a new perspective to construct surrogates based on the expressive power of neural networks \cite{raghu2017expressive, li2017convergence}. PINNs reformulate the original PDE-solving problem to an optimization problem by minimizing a loss function, consisting of the PDE loss and data loss. While PINNs have shown success in solving PDEs, they face significant limitations \cite{krishnapriyan2021characterizing}. One major issue is that in the presence of noisy data PINNs are prone to overfitting \cite{zhang2019quantifying, psaros2023uncertainty, doumeche2023convergence}, resulting in inaccurate predictions and a lack of uncertainty estimates. To address this, non-Bayesian methods \cite{abdar2021review, psaros2023uncertainty}, such as dropout and ensemble learning are employed to obtain uncertainty estimates. However, these methods are sensitive to hyper-parameters, such as the prior, and typically underestimate uncertainty. Later on, Bayesian physics-informed neural networks \cite{YANG2021109913} (B-PINNs) were introduced, incorporating a Bayesian framework to handle such cases. The key idea behind B-PINNs is to treat the parameters of the neural network as random variables. Therefore, by imposing a prior distribution on these parameters and incorporating observations, the original problem is reformulated as a high-dimensional Bayesian inverse problem, where the likelihood encodes the physical laws and captures the noise in the data. The posterior estimation of the parameters can then be obtained by sampling the posterior using the Hamiltonian Monte Carlo (HMC) \cite{cobb2021scaling} sampler, allowing for simultaneous inference of the solution and the unknown parameters. Additionally, the samples provide UQ estimates, hence, enhancing the reliability of the results.    

Despite these advantages, B-PINNs face several challenges. First, the multilayer perceptron (MLP) models commonly used in B-PINNs contain a large number of parameters, making them prone to overfitting. Another issue is that while HMC is considered the gold standard for accurate prediction and reliable UQ in small-noise, low-dimensional problems, its performance degrades in the presence of large datasets and substantial noise \cite{Papamarkou2022425}, making the inference process inefficient and possibly erroneous.  Moreover, obtaining the gradient information required for the sampling process is computationally expensive, exacerbating the existing challenges. To address these issues, some studies \cite{PENSONEAULT2024113006} have proposed to use the gradient-free ensemble Kalman inversion \cite{chada2020tikhonov,Kovachki_2019, liu2023dropout} (EKI) method for inference in B-PINNs. While EKI performs well in low-dimensional and small-noise scenarios, it also struggles to scale to more complex, higher-dimensional problems.

To alleviate these issues, we first propose to employ Chebyshev Kolmogorov-Arnold Networks (cKANs) \cite{ss2024chebyshev} as the network architecture in B-PINNs, which typically require fewer parameters with comparable accuracy as MLPs. Unfortunately, the back-propagation process for KANs is typically much slower compared to MLPs. Therefore, we propose to adopt EKI to replace HMC to accelerate training. However, it is well known that EKI may go through the so-called ensemble collapse \cite{iglesias2013ensemble, schillings2018convergence}, which is not practical for problems that require reliable UQ. To this end, we propose alleviating this issue by using the stochastic version with Tikhonov regularization and dropout constraint. Moreover, we employ the active subspace method \cite{constantine2014active} to provide a more accurate UQ estimation by reducing the overfitting phenomenon. The key idea lies in identifying the directions that contribute the most, as determined by the derivatives of the output with respect to the parameters. 

To demonstrate the efficiency of the proposed method, we perform  several experiments, including both forward and inverse problems. The results demonstrate that our method not only provides performance comparable to HMC when the noise scale is small but also scales effectively to large, high-dimensional problems with significantly reduced computational cost. Our main contributions can be summarized as follows:
\begin{itemize}
    \item We propose to replace the MLPs with cKANs to reduce the number of parameters and also accelerate convergence. The EKI method is used to achieve gradient-free inference with comparable performance with respect to HMC.
    \item To further reduce overfitting, we employ the Tikhonov regularization and also the dropout constraint to achieve better performance. Moreover, we adopt the active subspace method  to find the most sensitive directions that can accelerate training and reduce overfitting as well.
    \item To test the effectiveness of the proposed method, several experiments have been implemented. The results show that our method can achieve comparable performance, and can also scale well for large-scale and high-dimensional problems. 
\end{itemize}

The remainder of this paper is organized as follows. In Section 2, we introduce the original B-PINNs as well as cKANs. In Section 3, we introduce the ensemble inversion method and also the active subspace method. To confirm the efficiency of our method, several benchmarks are tested in Section 4. We conclude the paper with a summary in Section 5.

\section{Background}
In this section, we first give a brief review of B-PINNs. Then, we will introduce the basic concepts of cKANs.

\subsection{Bayesian physics-informed neural networks \cite{YANG2021109913} (B-PINNs)}
We consider the following PDE:
\begin{equation}
    \label{PDE}
    \begin{split}
    \mathcal{N}_{x}(u(\mathbf{x});\lambda)&=f(\mathbf{x})\quad\mathbf{x}\in\Omega,\\ 
    B_{x}(u(\mathbf{x});\lambda)&=b(\mathbf{x})\quad\mathbf{x}\in\partial\Omega,
    \end{split}
\end{equation}
where $\mathcal{N}_x$ and $B_x$ denote the differential and boundary operators, respectively. The spatial domain $\Omega\subseteq\mathbb{R}^d$ has a boundary $\partial\Omega$, and $\lambda\in\mathbb{R}^{N_\lambda}$ represents a vector of unknown physical parameters. Suppose the forcing term $f(\mathbf{x})$ and boundary term $b(\mathbf{x})$ are given. When $\lambda$ is known, the forward problem is to infer the true solution $u$ with noisy observations $\mathcal{D}_{f} = \{\mathbf{x}_{f}^{i}, f(\mathbf{x}_{f}^{i})\}_{i=1}^{N_{f}}$ with $\mathbf{x}_{f}^{i}\in \Omega$ and $\mathcal{D}_{b} = \{\mathbf{x}_{b}^{i}, b(\mathbf{x}_{b}^{i})\}_{i=1}^{N_{b}}$ with $\mathbf{x}_{b}^{i}\in \partial \Omega$. Moreover, when $\lambda$ is unknown, then with additional observations $\mathcal{D}_{u} = \{\mathbf{x}_{u}^{i}, u_{i}\}_{i=1}^{N_{u}}$, the goal of the inverse problems is to infer $u$ and $\lambda$ simultaneously via the Bayesian formulation. 

B-PINNs typically employ multilayer perception networks (MLPs) $\tilde{u}(\mathbf{x},\mathbf{\theta})$ to approximate the solution $u(\mathbf{x})$, where $\mathbf{\theta}\in \mathbb{R}^{N_{\theta}}$ is a vector of parameters in the network. By imposing a prior distribution $p(\mathbf{\theta})$ on the network parameters, the original optimization problem is transformed into a Bayesian inverse problem. Together with the parameter $\lambda$, the goal of the inverse problem is to infer the solution $u(\mathbf{x})$ and parameters $\mathbf{\xi} = \{\mathbf{\theta}, \lambda\}$ by exploring the posterior conditioned on the residual measurements $\mathcal{D}_{f}$, boundary measurements $\mathcal{D}_{b}$ and solution measurements $\mathcal{D}_{u}$. Specifically, by Bayes' formula and assuming the independence of the datasets, we have 
\begin{equation}
    \label{posterior}
p(\mathbf{\xi}|\mathcal{D}_{u}, \mathcal{D}_{b}, \mathcal{D}_{f})\propto p(\mathbf{\xi})p(\mathcal{D}_{u}|\mathbf{\xi})p(\mathcal{D}_{f}|\mathbf{\xi})p(\mathcal{D}_{b}|\mathbf{\xi}).
\end{equation}
Suppose the parameters are independent and $\mathbf{\theta}$ follows the zero mean Gaussian distribution, the prior can be decomposed as 
\begin{equation*}
    \label{prior}
p(\mathbf{\xi}) = p(\theta)p(\lambda) =p(\lambda)\prod_{i=1}^{N_\theta}p(\theta^i),\quad p(\theta^i)\sim\mathcal{N}\left(0,\sigma_\theta^i\right),
\end{equation*}
where $\sigma_{\theta}^{i}$ is standard deviation of the corresponding network parameter $\theta^{i}$. Similarly, by assuming the noises added to the data points follow the Gaussian distribution, we can specify the likelihoods as 
\begin{equation*}
    \label{likelihood}
    \begin{aligned}&p(\mathcal{D}_{u}|\xi)=\prod_{i=1}^{N_{u}}p(u^{i}|\xi),\quad p(\mathcal{D}_{f}|\xi)=\prod_{i=1}^{N_{f}}p(f^{i}|\xi),\quad p(\mathcal{D}_{b}|\xi)=\prod_{i=1}^{N_{b}}p(b^{i}|\xi),\\&p(u^{i}|\xi)=\frac{1}{\sqrt{2\pi\sigma_{\eta_{u}}^{2}}}\exp\left(-\frac{\left(u^{i}-\tilde{u}(\mathbf{x}_{u}^{i};\mathbf{\theta})\right)^{2}}{2\sigma_{\eta_{u}}^{2}}\right),\\&p(f^{i}|\xi)=\frac{1}{\sqrt{2\pi\sigma_{\eta_{f}}^{2}}}\exp\left(-\frac{\left(f^{i}-\mathcal{N}_{x}(\tilde{u}(\mathbf{x}_{f}^{i};\mathbf{\theta});\mathbf{\lambda})\right)^{2}}{2\sigma_{\eta_{f}}^{2}}\right),\\&p(b^{i}|\xi)=\frac{1}{\sqrt{2\pi\sigma_{\eta_{b}}^{2}}}\exp\left(-\frac{\left(b^{i}-B_{x}(\tilde{u}(\mathbf{x}_{b}^{i};\mathbf{\theta});\mathbf{\lambda})\right)^{2}}{2\sigma_{\eta_{b}}^{2}}\right),
    \end{aligned}
\end{equation*}
where we assume that the standard deviations $\sigma_{\eta_{u}}, \sigma_{\eta_{f}}, \sigma_{\eta_{b}}$ are the same for different points in the same category. Moreover, the prior for the physical parameters $\lambda$ is problem-dependent. After specifying priors and likelihoods, the posterior can be explored using sampling methods, among which the HMC is the gold standard for addressing the cases with small data with a small noise scale. However, problems corrupted by high noise levels require a lot of data points to increase the accuracy of the prediction. In such cases, the likelihood will dominate, effectively turning the problem into an optimization task to find the MAP point, rendering the prior almost irrelevant. Furthermore, thousands of network parameters will increase the overfitting, leading to very shallow uncertainty estimation. To alleviate these issues, we introduce KANs \cite{liu2024kan, xu2024fourierkan, aghaei2024fkan, ss2024chebyshev}, especially Chebyshev KANs (cKANs) \cite{ss2024chebyshev}, which may be better candidates with much fewer parameters and may possibly converge faster during the training. Thus, in the next subsection, we will briefly introduce the cKANs. 

\subsection{Chebyshev Kolmogorov-Arnold Networks (cKANs)}
Suppose $f:[0,1]^n\to\mathbb{R}, \mathbf{x}\in[0,1]^{n}$ is continuous, the Kolmogorov-Arnold theorem \cite{schmidt2021kolmogorov} guarantees the existence of continuous univariate functions $\Phi_{q}$ and $\phi_{q,p}$ such that:

$$f(x_1,x_2,\ldots,x_n)=\sum_{q=0}^{2n}\Phi_q\left(\sum_{p=1}^n\phi_{q,p}(x_p)\right),$$
where $x=(x_1,x_2,\ldots,x_n)$. Here, $\phi_{q,p}:[0,1]\to\mathbb{R}$
and $\Phi_q:\mathbb{R}\to\mathbb{R}$ are continuous functions. Typical choices of the basis functions include the B-splines and the orthogonal polynomials. To improve the interpretability and numerical stability, the Chebyshev polynomials defined in $[-1,1]$ are better candidates \cite{ss2024chebyshev}. Specifically, the Chebyshev polynomials $T_{n}(x)$ are defined  recursively as
\begin{equation*}
    \label{chebyshev}
    \begin{array}{c}
        T_{0}(x) = 1, T_{1}(x) = x,\\[0.2cm] 
        T_{n}(x) = 2xT_{n-1}(x) - T_{n-2}(x), \quad n\geq 2.
    \end{array}
\end{equation*}
\begin{figure}[htbp]
    \centering
    \includegraphics[width=0.5\linewidth]{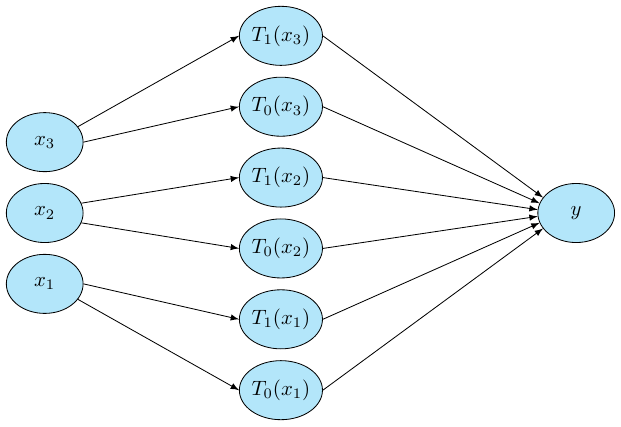}
    \caption{The framework of a single layer Chebyshev KAN with input dimension 3, output dimension 1, and degree 1.}
    \label{kan}
\end{figure}
Combining the definition of Chebyshev polynomials and KANs, the framework, depicted in Fig.~\ref{kan} is straightforward, where the back-propagation is the same as MLPs. By constructing a loss function $\mathcal{L}(\theta)$, the remaining problem is to find
\begin{equation*}
    \theta^{*} = \arg \min_{\theta\in \Theta} \mathcal{L}(\theta)
\end{equation*}
with appropriate optimizers \cite{kingma2014adam, amari1993backpropagation}.
The summary of the training process for a single-layer network is summarized in Algorithm \ref{ChebyKans}.
\begin{algorithm}[htbp]
    \caption{Chebyshev Kolmogorov-Arnold Networks (cKANs)}
    \label{ChebyKans}
    \begin{algorithmic}[1]
        \REQUIRE Input dimension $n_{in}$, output dimension $n_{out}$, degree $d$, input $\mathbf{x}$.
        \STATE Project $\mathbf{x}$ into $[-1,1]$ using the $\mathtt{tanh}$ function:
        \begin{equation*}
            \tilde{\mathbf{x}} = \tanh{\mathbf{x}}, \tilde{\mathbf{x}}\in [-1,1]^{n_{in}}.
        \end{equation*}
        \STATE Forward pass using Chebyshev polynomials with maximum degree $d$:
        \begin{equation*}
            \tilde{f}(\mathbf{x})=\sum_{j=1}^{n_{\mathrm{in}}}\sum_{k=0}^{d}\theta_{j,k}T_{k}(\tilde{x}_{j}),
        \end{equation*}
        where $\mathbf{\theta}\in\mathbb{R}^{n_{in}\times\ n_{out}\times (d+1)}$ are the parameters of the network.

        \STATE Find optimal $\mathbf{\theta}^{*}$ by minimizing the following loss function in parameter space $\Theta$:
        \begin{equation*}
            \mathbf{\theta}^{*}  = \arg\min_{\mathcal{\theta} \in \Theta} \mathcal{L}(\mathbf{\theta}).
        \end{equation*}
    \end{algorithmic}
\end{algorithm}
However, according to the definition, the back-propagation for cKANs is not as efficient as MLPs since the forward pass is based on a three-term recurrence formula. This phenomenon is especially obvious in training with large amounts of data. Therefore, to address problems with big data, the training process can be quite slow, even with mini-batch strategies. To address this issue, we adopt the gradient-free inference method, namely, the ensemble Kalman inversion (EKI) to train the model in parallel, which will be introduced in the following section.

\section{Subspace-enhanced dropout Tikhonov EKI for B-PINNs}
In this section, we introduce the gradient-free dropout Tikhonov EKI \cite{liu2023dropout} (DTEKI) for training the B-PINNs with cKANs, which can greatly save computational cost while achieving performance comparable to HMC.
\subsection{Dropout Tikhonov EKI}
Ensemble Kalman methods are popular methods to solve inverse problems;  they are gradient-free and highly efficient even for high-dimensional problems. Consider the model as follows:
\begin{equation}
    y = \mathcal{G}(\mathbf{\xi}) + \eta_{1}, 
\end{equation}
where $\mathcal{G}:\mathbb{R}^{N_{\xi}}\to \mathbb{R}^{N_{y}}$ is the forward operator and $\eta_{1}\sim \mathcal{N}(0, \Gamma)$ is white noise with covariance matrix $\Gamma\in \mathbb{R}^{N_{y}\times N_{y}}$, which leads to noisy observations $y$. Given a prior distribution $p(\xi)$ and assuming that it also follows a Gaussian distribution, e.g.,  $\xi\sim\mathcal{N}(0, \mathcal{C}_{0})$, the Tikhonov EKI considers the following dynamical system
\begin{equation}
    \label{dynamical_system}
    \begin{split}
    y &= \mathcal{G}(\xi) + \eta_{1},\quad \eta_{1}\sim \mathcal{N}(0, \Gamma),\\[0.1cm] 
    0 &= \xi + \eta_{2},\quad \eta_{2}\sim \mathcal{N}(0, \alpha^{-1}\mathcal{C}_{0}),
    \end{split}
\end{equation}
which leads to the following Tikhonov functional:
\begin{equation}
    \label{potential}
    I(\xi) = \frac{1}{2}\|y-\mathcal{G}(\xi)\|_{\Gamma}^{2} + \frac{\alpha}{2}\|\xi\|_{\mathcal{C}_{0}}^{2},
\end{equation}
where $\alpha$ is an additional regularization parameter, the norm $\|\cdot\|_{A} = \|A^{-\frac{1}{2}}\cdot\|$ and $\|\cdot\|$ denotes the vector $l_{2}$ norm. By denoting $\mathcal{H} = \left[\begin{array}{c}
     \mathcal{G}(\xi)  \\
      \xi
\end{array}\right]$ and $z = \left[\begin{array}{c}
     y  \\
      0
\end{array}\right]$, 
Eq. \eqref{potential} can be rewritten as 
\begin{equation}
    \label{3.4}
    I(\xi) = \frac{1}{2}\|z - \mathcal{H}(\xi)\|^{2}_{\Gamma_{\mathcal{H}}},
\end{equation}
where $\Gamma_{\mathcal{H}} = \left[\begin{array}{cc}
    \Gamma & 0 \\
     0&\alpha^{-1}\mathcal{C}_{0} 
\end{array}\right]$.
The methodology of stochastic Tikhonov EKI to solve this problem is to construct an ensemble $\{\xi_{n}^{(j)}\}_{j=1}^{J}$ from the prior and then apply the following iteration to update the ensemble:
\begin{equation}
    \label{3.5}
    \xi_{n+1}^{(j)}=\xi_n^{(j)}+C_n^{\xi z}(C_n^{zz}+\Gamma_{\mathcal{H}})^{-1}(z-z_{n}^{(j)})),
\end{equation}
with i.i.d predictions
\begin{equation}
\label{3.6}
    z_{n}^{(j)} = \mathcal{H}(\xi_{n}^{(j)}) + \eta_{n}^{(j)}, \quad \eta_{n}^{(j)} \sim \mathcal{N}(0,\Gamma_{\mathcal{H}}),
\end{equation}
where $C_{n}^{\xi z}, C_{n}^{zz}$ are the empirical covariance matrices of the ensemble. In detail, denote
\begin{equation}
\label{3.7}    \begin{aligned}&\overline{\xi}_{n}=\frac{1}{J}\sum_{j=1}^{J}\xi_{n}^{(j)},\quad\tau_{n}^{(j)}=\xi_{n}^{(j)}-\overline{\xi}_{n},\quad\overline{\mathcal{H}(\xi_{n})}=\frac{1}{J}\sum_{j=1}^{J}{\mathcal{H}(\xi_{n}^{(j)})}.\\&\begin{aligned}&C_n^{zz}=\frac{1}{J-1}\sum_{j=1}^J\left(\mathcal{H}(\xi_n^{(j)})-\overline{\mathcal{H}(\xi_n)}\right)\otimes\left(\mathcal{H}(\xi_n^{(j)})-\overline{\mathcal{H}(\xi_n)}\right).\\
    &C_n^{\xi z}=\frac{1}{J-1}\sum_{j=1}^J\tau_n^{(j)}\otimes\left(\mathcal{H}(\xi_n^{(j)})-\overline{\mathcal{H}(\xi_n)}\right).\end{aligned}\end{aligned}
\end{equation}
For low-dimensional cases, the approximation can be accurate, leading to desired convergence rate. However, for high-dimensional cases, the estimated covariance matrix could be rank-deficient when $J < N_{\xi}$, greatly influencing the stability \cite{liu2023dropout}. To deal with this issue, the dropout technique, often used in deep learning to prevent overfitting, is introduced here to give more accurate estimates. The only modification is to apply dropout on the deviations $\tau_{n}^{(j)}$ to obtain the dropout version $\tilde{\tau}_{n}^{(j)}$, then 
\begin{equation}
\label{3.8}
\tilde{\xi}_{n}^{(j)}=\overline{\xi}_{n}+\tilde{\tau}_{n}^{(j)},\quad\tilde{\tau}_{n}^{(j)}=\beta_{n}\circ\tau_{n}^{(j)},\quad\beta_{n}(s)\overset{\mathrm{i.i.d.}}{\operatorname*{\sim}}\mathrm{Bernoulli}(\rho),
\end{equation}
where $1 - \rho\in (0,1)$ is the dropout rate and $\beta_{n}$ is the mask vector at iteration $n$. Hence, the dropout ensemble covariance matrix reads
\begin{equation}
    \label{3.9}
    \begin{aligned}&\overline{\mathcal{H}(\tilde{\xi}_{n})}=\frac{1}{J}\sum_{j=1}^{J}{\mathcal{H}(\tilde{\xi}_{n}^{(j)})},\quad\tilde{C}_{n}^{\xi z}=\frac{1}{J-1}\sum_{j=1}^{J}\tilde{\tau}_{n}^{(j)}\otimes\left(\mathcal{H}(\tilde{\xi}_{n}^{(j)})-\overline{\mathcal{H}(\tilde{\xi}_{n})}\right).\\&\tilde{C}_{n}^{zz}=\frac{1}{J-1}\sum_{j=1}^{J}\left(\mathcal{H}(\tilde{\xi}_{n}^{(j)})-\overline{\mathcal{H}(\tilde{\xi}_{n})}\right)\otimes\left(\mathcal{H}(\tilde{\xi}_{n}^{(j)})-\overline{\mathcal{H}(\tilde{\xi}_{n})}\right).\end{aligned}
\end{equation}
For linear problems, it has been proven that the mean-field limit obtained when $J\to \infty$ will converge to the minimizer of $I(\xi)$ \cite{ding2021ensemble}. However, one limitation of EKI-type methods is the \textit{ensemble collapse}, which means $\lim_{n\to\infty}\tau_{n}^{j} = 0, \forall j.$ Therefore, the ensemble covariance will converge to zero matrices, leading to an underestimation of the true posterior covariance matrix. Moreover, due to the subspace property of EKI, $\xi_{n}^{(j)}\in A = span\{\xi_{0}^{(j)}\}_{j=1}^{J}$, which leads to worse estimates when the initial ensemble is chosen inappropriately. To alleviate these issues, several techniques are used in this paper. Firstly, the mini-batch training during the process can be applied to reduce overfitting, which could greatly accelerate training speed as well. Secondly, at each iteration, additional noises are added to the parameters to force the ensemble to get rid of the original subspace. Specifically, at iteration $n$:  
\begin{equation}
\label{3.10}
\begin{split}
    \hat{\xi}_{n}^{(j)} &= \xi_{n}^{(j)} + \epsilon_{n}^{(j)},\quad \epsilon_{n}^{(j)}\sim \mathcal{N}(0, Q). 
\end{split}
\end{equation}
The ensemble covariance matrix in Eq. \eqref{3.9} can also be calculated accordingly by using these perturbed versions, denoted as $\hat{C}_n^{zz}, \hat{C}_n^{\xi z}$.
To this end, the iteration becomes
\begin{equation}
\label{3.11}
    \xi_{n+1}^{(j)}=\hat{\xi}_n^{(j)}+\hat{C}_n^{\xi z}(\hat{C}_n^{zz}+\Gamma_{\mathcal{H}})^{-1}(z-\mathcal{H}(\hat{\xi}_n^{(j)}) - \eta_{n}^{(j)}).
\end{equation}
In the context of B-PINNs, the parameters $\xi = \{\mathbf{\theta}, \lambda\}$, and the forward operator $\mathcal{G}$ is specified as 
\begin{equation}
\label{forward}
    \mathcal{G}(\xi) = \left[\begin{array}{c}
          \mathcal{N}_{\mathbf{x}}(u(\mathbf{x};\mathbf{\theta});\lambda) \\
          B_{\mathbf{x}}(u(\mathbf{x};\mathbf{\theta});\lambda)\\ 
          u(\mathbf{x};\mathbf{\theta})
    \end{array}\right]
\end{equation}
including the partial differential operator, boundary operator, and solution operator by fixing $x$. Similarly, the observations can be specified as \begin{equation} 
\label{observation}
y = \left[\begin{array}{ccc}
         \mathbf{f}  &
          \mathbf{b} &
          \mathbf{u}\end{array}\right]^{T},\end{equation} where $\mathbf{f} = \{f(\mathbf{x}_{f}^{i})\}_{i=1}^{N_{f}},\mathbf{b} = \{b(\mathbf{x}_{b}^{i})\}_{i=1}^{N_{b}},\mathbf{u} = \{u(\mathbf{x}_{u}^{i})\}_{i=1}^{N_{u}}$. Therefore, $\mathcal{H}$ and ${z}$ can be defined based on the definition. The above procedures can be summarized as Algorithm \ref{DTEKI}.
          \begin{algorithm}
              \centering 
              \caption{Bayesian physics-informed cKANs}
              \label{DTEKI}
              \begin{algorithmic}[1]
                  \REQUIRE Observations $y$, covariance $Q$, ensemble size $J$, dropout rate $\rho$, regularization parameter $\alpha$, number of iterations $N$.
                  \STATE Generate initial ensemble from $p(\xi)$,
                  \begin{equation*}
                      \{\xi_{0}^{(j)}\}_{j=1}^{J}\sim p(\xi).
                  \end{equation*}
                  \STATE $n\leftarrow 0$.
                  \WHILE{$n \leq N$}
                  \STATE Perturb the ensemble Eq.  \eqref{3.10} to get $\{\hat{\xi}_{n}^{(j)}\}_{j=1}^{J}$.
                  \STATE Calculate ensemble covariance matrices $\hat{C}_n^{zz}, \hat{C}_n^{\xi z}$ according to  Eqs. \eqref{3.7}-\eqref{3.9}.
                  \STATE Update ensemble using Eq.  \eqref{3.11} to obtain $\{\xi_{n+1}^{(j)}\}_{j=1}^{J}$.
                  \STATE $n\leftarrow n+1$.
                  \ENDWHILE
                  \RETURN Final ensemble $\{\xi_{N}^{(j)}\}_{j=1}^{J}$.
              \end{algorithmic}
          \end{algorithm}
          
The selection of hyperparameters plays a critical role in the performance of the proposed method. In this study, the dropout rate is consistently set to 0.8 to enhance the stability of the covariance matrix. For the regularization parameter 
$\alpha$, its primary purpose is to mitigate overfitting, which tends to occur when the dataset is small or the network architecture is overly complex. In such scenarios, a higher value of
$\alpha$ is recommended to achieve better results. For the choice of perturbation matrix $Q$, it is reasonable to choose different sensitive levels for the physical parameters $\lambda$ and network parameters $\theta$ respectively. Generally, a larger perturbation level is chosen for $\lambda$ to increase the accuracy of the posterior estimation. 

However, for neural network-based surrogate models, the number of network parameters, often on the order of $\mathcal{O}(10^{3})$, can still be substantial, resulting in significant computational costs. Moreover, an increased number of parameters heightens the risk of overfitting, making the selection of network hyperparameters challenging in practice. To address these challenges, we introduce an approach known as the {\em active subspace method} \cite{constantine2014active} in the next section, which identifies the most sensitive directions in the parameter space, thereby mitigating these issues by reducing the network parameter dimension.
\subsection{Active subspace for cKANs}
Consider a function with $N_{\mathbf{\xi}} = N_{\mathbf{\lambda}} + N_{\mathbf{\theta}}$ continuous inputs,
\begin{equation}
    \label{continuous_functions}
    \mathcal{G}(\mathcal{\xi}) = \mathcal{G}(\mathbf{\lambda},\mathbf{\theta}) = \mathcal{F}_{\mathbf{\lambda}}(\mathbf{\theta}),\quad \mathbf{\xi}\in \mathbb{R}^{N_{\mathbf{\xi}}},\mathbf{\lambda}\in \mathbb{R}^{N_{\mathbf{\lambda}}}, \mathbf{\theta}\in \mathbb{R}^{N_{\mathbf{\theta}}},
\end{equation}
and let $p(\mathbf{\xi}) = p(\mathbf{\lambda})p(\mathbf{\theta})$ be the prior distribution associated with $\mathbf{\xi}$, where $p(\mathbf{\lambda}), p(\mathbf{\theta})$ are the priors for $\mathbf{\lambda}$ and $\mathbf{\theta}$ respectively. To identify the most sensitive directions in the network parameter space $\Theta = \mathbb{R}^{N_{\theta}}$ that capture the majority of the gradient variability of the function $\mathcal{F}_{\lambda}$, we consider the uncentered covariance matrix of the gradients, defined as 
\begin{equation}
    \label{covariance_gradient}
    \begin{aligned}\mathcal{C}&= \mathbb{E}_\theta \left[ \left( \nabla_{\theta}\mathcal{F}_{\lambda}(\mathbf{\theta} \right) \left(\nabla_{\theta}\mathcal{F}_{\lambda}(\mathbf{\theta} \right)^\top \right]\\&= \int_{\Theta} \left( \nabla_{\theta}\mathcal{F}_{\lambda}(\mathbf{\theta} \right) \left( \nabla_{\theta}\mathcal{F}_{\lambda}(\mathbf{\theta} \right)^\top p(\theta) \, d\theta.\end{aligned}
\end{equation}
Since $\mathcal{C}$ is symmetric semi-definite, then it admits an eigendecomposition as 
\begin{equation}
    \label{eigen_decom}
    \mathcal{C}=\mathbf{W\Lambda W}^T=[\mathbf{W}_1\; \mathbf{W}_2]\begin{bmatrix}\mathbf{\Lambda}_1&0\\0&\mathbf{\Lambda}_2\end{bmatrix}[\mathbf{W}_1\; \mathbf{W}_2]^T,
\end{equation}
where $\mathbf{W}, \mathbf{\Lambda}$ contain the eigenfunctions and eigenvalues respectively. Suppose the eigenvalues in $\mathbf{\Lambda}$ are sorted in descending order. In this regard, $\mathbf{W}_{1}$ represents the directions that the function varies most since $\mathbf{W}_{1}$ contains the largest eigenvalues. Formally, $\mathbf{W}_{1}$ forms a subspace with directions that $\mathcal{F}_{\lambda}(\mathbf{\theta})$ is most sensitive to. Suppose we truncate $m$ terms in the decomposition, that is $\mathbf{W}_{1}\in \mathbb{R}^{N_{\mathbf{\theta}}\times m}$, the function $\mathcal{G}$ can be approximated as 
\begin{equation*}
    \mathcal{G}(\mathbf{\xi})\approx \mathcal{G}(\mathbf{\lambda},\mathbf{W_{1}}\mathbf{\omega}) = \mathcal{F}_{\mathbf{\lambda}}(\mathbf{W_{1}}\mathbf{\omega}),\quad \mathbf{\omega}\in \mathbb{R}^{m}, \mathbf{\omega}\sim p(\mathbf{\omega}). 
\end{equation*}
Note that $p(\mathbf{\omega})$ can be Gaussian if $p(\mathbf{\theta})$ is Gaussian due to the linear transformation. In practice, the matrix $\mathcal{C}$ can be estimated using the Monte Carlo method. Specifically, we first sample $M$ samples $\{\mathbf{\theta}_{i}\}_{i=1}^{M}$ from $p(\mathbf{\theta})$, then, the matrix is estimated by 
\begin{equation}
    \label{estimate_matrix}
    \mathcal{C}\approx\tilde{\mathcal{C}}=\frac{1}{M}\sum_{i=1}^{M}(\nabla_{\mathbf{\theta}}\mathcal{F}_{\mathbf{\lambda}}^{i})(\nabla_{\mathbf{\theta}}\mathcal{F}_{\mathbf{\lambda}}^{i})^{T},
\end{equation}
where $\nabla_{\mathbf{\theta}}\mathcal{F}_{\mathbf{\lambda}}^{i} = \nabla_{\mathbf{\theta}}\mathcal{F}_{\mathbf{\lambda}}(\mathbf{\theta}_{i}), i = 1,\ldots, M$. Then, eigenvalue decomposition can be implemented to find the subspace. However, this procedure can also be expensive when $N_{\mathbf{\theta}}$ is large. Therefore, since $\tilde{\mathcal{C}} = \mathbf{GG}^{T}$ and 
$$\mathbf{G} = \frac{1}{\sqrt{M}}[\nabla_{\mathbf{\theta}}\mathcal{F}_{\lambda}^{1}\;  \cdots\; \nabla_{\mathbf{\theta}}\mathcal{F}_{\lambda}^{M}],$$ 
we can instead calculate the singular value decomposition (SVD) of $\mathbf{G}$, that is
\begin{equation}
    \textbf{G} = \tilde{\mathbf{W}}\sqrt{\mathbf{\tilde{\Lambda}}}\mathbf{V}^{T} = [\tilde{\mathbf{W}}_1\; \tilde{\mathbf{W}}_2]\begin{bmatrix}\sqrt{\tilde{\mathbf{\Lambda}}_1}&0\\0&\sqrt{\tilde{\mathbf{\Lambda}}_2}\end{bmatrix}\mathbf{V}^T.
\end{equation}
Thus, 
\begin{equation}
    \label{forward_approximation}
    \mathcal{G}(\mathbf{\xi})\approx \mathcal{G}(\mathbf{\lambda},\mathbf{\tilde{W}_{1}}\mathbf{\omega}) = \mathcal{F}_{\lambda}(\tilde{W}_{1}\mathbf{\omega}). 
\end{equation}
The whole parameters now are $\zeta = \{\lambda, \omega\}$, with prior $p(\zeta) = p(\lambda)p(\omega)$. Substituting Eq. \ref{forward_approximation} into the iteration process of DTEKI, we can get the subspace version of DTEKI (SDTEKI), which could be more efficient and stable. Formally, the subspace-enhanced DTEKI (SDTEKI) can be summarized as Algorithm \ref{algorithm3.2}.
\begin{algorithm}[htbp]
    \centering 
    \caption{Subspace-enhanced DTEKI for Bayesian physics-informed cKANs}
    \label{algorithm3.2}
    \begin{algorithmic}[1]
        \REQUIRE Observations $y$, covariance $Q$, ensemble size $J$, projection matrix $\tilde{\mathbf{W}}_{1}$, dropout rate $\rho$, regularization parameter $\alpha$, number of iterations $N$.
        \STATE Generate initial ensemble from $p(\zeta)$,
                \begin{equation*}
            \{\lambda_{0}^{(j)}, \omega_{0}^{(j)}\} = \{\mathbf{\zeta}_{0}^{(j)}\}_{j=1}^{J}\sim p(\zeta).
                  \end{equation*}.
        \STATE $n\leftarrow 0$.
        \WHILE{$n\leq N$} 
        \STATE Perturb the ensemble and calculate ensemble deviations for $j = 1,\ldots, J$
        \begin{equation*}
        \begin{split}
        \hat{\tau}_{n}^{(j)} = \hat{\zeta}_{n}^{(j)} - \overline{\hat{\zeta}}_{n},\; \overline{\hat{\zeta}}_{n} = \frac{1}{J}\sum_{j=1}^{J}\hat{\mathbf{\zeta}}_{n}^{(j)},\;
        \hat{\zeta}_{n}^{(j)} = \zeta_{n}^{(j)} + \eta_{n}^{(j)}.
        \end{split}
        \end{equation*}
        \STATE Get the dropout version of the ensemble for $j=1,\ldots, J$
        \begin{equation*}
            \tilde{\zeta}_{n}^{(j)}=\overline{\hat{\zeta}}_{n}+\tilde{\tau}_{n}^{(j)},\quad\tilde{\tau}_{n}^{(j)}=\beta_{n}\circ\hat{\tau}_{n}^{(j)},\quad\beta_{n}(s)\overset{\mathrm{i.i.d.}}{\operatorname*{\sim}}\mathrm{Bernoulli}(\rho).
        \end{equation*}
        \STATE Recover $\tilde{\xi}_{n}^{(j)} \approx \{\tilde{\lambda}_{n}^{(j)},\tilde{\mathbf{W}}_{1}\tilde{\omega}_{n}^{(j)}\}$ and calculate the ensemble covariance matrices according to Eqs.~\eqref{3.7}-\eqref{3.9}.
        \STATE Update ensemble using Eq.~\eqref{3.11} to obtain $\{\zeta_{n+1}^{(j)}\}_{j=1}^{J}$.
        \STATE $n\leftarrow n + 1$.
        \ENDWHILE
        \RETURN Final ensemble $\{\zeta_{N}^{(j)}\}_{j=1}^{J}$.
    \end{algorithmic}
\end{algorithm}
The remaining challenge lies in determining the appropriate subspace dimension using SVD. In our study, the singular values exhibit rapid decay, with the top one-third of singular values accounting for over 
$99.9\%$ of the total weight. This behavior will be illustrated in the first example. Consequently, it is reasonable to select the top one-third of the basis vectors in 
$\textbf{W}$ to construct the reduced subspace $\textbf{W}_{1}$.
\begin{remark}
    Note that if $\mathbf{\lambda}$ is an infinite-dimensional stochastic process, we can also parameterize it as a neural network and then construct the active subspace using the method stated above. Then, the network parameters can be inferred by using the DTEKI method.
\end{remark}

\subsection{Complexity analysis}
The computational cost of DTEKI is relatively low compared to HMC, as it primarily focuses on Eqs. \eqref{3.8} and \eqref{3.10}. At each iteration, the computational complexity for DTEKI is given by:
\begin{equation}
\mathcal{O}((N_{\mathbf{\xi}}+N_{y})^{3}+JN_{y}^{2}+ JN_{\mathbf{\xi}}^{2} + JN_{y}N_{\mathbf{\xi}}).
\end{equation} 
This can be divided into the following components:
\begin{itemize}
    \item $\mathcal{O}(JN_{y}^{2}+ JN_{\mathbf{\xi}}^{2} + JN_{y}N_{\mathbf{\xi}})$ for $\hat{C}_n^{zz}, \hat{C}_n^{\xi z}$ in Eq. \eqref{3.8}.
    \item $\mathcal{O}((N_{\mathbf{\xi}}+N_{y})^{3}+JN_{y}^{2}+ JN_{\mathbf{\xi}}^{2} + JN_{y}N_{\mathbf{\xi}})$ for Eq. \eqref{3.10}.
\end{itemize}
When $N_{y}$ is small, the process is highly efficient due to the limited number of parameters, especially when using cKANs. However, for problems with large datasets, applying mini-batching techniques is recommended to accelerate training.

In contrast, when $N_{\mathbf{\xi}}$
  is large, the memory required to store the ensemble matrix can become significant. This issue can be mitigated using the active subspace method. The computational complexity for subspace-enhanced DTEKI at each iteration is:
  \begin{equation}
  \mathcal{O}(JmN_{\mathbb{\theta}}+ (m+N_{y})^{3}+JN_{y}^{2}+ Jm^{2} + JN_{y}m), 
  \end{equation}
where $m$ is the reduced subspace dimension. When $m\ll N_{\mathbf{\theta}}$, this approach significantly reduces the computational cost. Moreover, the memory cost for storing the ensemble matrix is also reduced, which could provide better scalability to large-scale problems. Additionally, the cost of identifying the subspace is negligible, as it is an offline task that can be reused across multiple problems. Thus, SDTEKI can perform better than DTEKI with comparable accuracy.
\begin{remark}
    The above analysis assumes that the forward evaluation of $\mathcal{G}$ is extremely fast and can be considered negligible. If this assumption does not hold, parallelization techniques can be employed to accelerate the computation. Furthermore, compared to HMC, this process is generally more efficient, as it requires solving an ODE with many steps at each iteration. However, the need to compute gradients during each iteration may still slow down the sampling process.
\end{remark}

\section{Results}
In this section, we will demonstrate the performance of our proposed DTEKI and subspace DTEKI (SDTEKI) with cKANs for B-PINNs. In addition, we will compare the results with those obtained using HMC with cKANs. To highlight the advantages of cKANs over MLPs, we will also present comparisons of results using DTEKI (MLP-DTEKI) with different architectures.

To fully demonstrate the advantages of our method, we will introduce different noise levels to the dataset. In addition, we will evaluate the performance on large-scale datasets corrupted by significant noise. Furthermore, we will test a high-dimensional inverse problem to showcase the scalability of our method.

For all examples, the network architectures have two hidden layers with 10 neurons, where the degree of Chebyshev polynomials is at most 7. The activation function is the $\mathtt{sigmoid}$ function. Furthermore, we will assume that the noise scales for the measurements are known. The priors are all standard Gaussian distributions with zero mean. The dropout rate $\rho=0.8$ for all examples. The regularization parameter $\alpha=0.1$ for the small noise scale and $\alpha = 0.01$ for the large noise scale without specification. The perturbation matrix is $$Q = \left(\begin{array}{cc}
0.01^{2}I_{N_{\lambda}} & 0\\ 0 & 0.002^{2}I_{N_{\theta}}\end{array}\right).$$ The ensemble size is $J = 500$ for all examples. For all problems, we will run 1000 iterations using different methods. For SDTEKI, we generate $M=1000$ samples to construct the subspace and choose the subspace dimension to  $m = \frac{1}{3}N_{\mathbf{\theta}}$. 

 For HMC B-PINNs, the leapfrog step \cite{cobb2021scaling} $L = 50$, and the initial step size is 0.1 adaptively in order to reach an acceptance rate of $60\%$ during burn-in steps. We draw a total of 1000 samples following 3000 burn-in steps.

 To compare the performance between different methods, we will use the relative $L_{2}$ error defined as 
 \begin{equation}
     e_u=\sqrt{\frac{\sum_{i=1}^{N}\left|u(\mathbf{x}^i)-\bar{u}(\mathbf{x}^i)\right|^2}{\sum_{i=1}^{N}\left|u(\mathbf{x}^i)\right|^2}},\quad e_\lambda=\frac{\left|\lambda-\bar{\lambda}\right|}{\left|\lambda\right|},
 \end{equation}
 where $\{\mathbf{x}_{i}\}_{i=1}^{N}$ are the test points in the physical  domain, and $u,\lambda$ are the reference solutions; $\overline{\lambda}$ and $\overline{u}$ are the ensemble mean of the predictions. Furthermore, we used the standard deviation of the ensemble to demonstrate the uncertainty estimate. All computations are performed on a standalone NVIDIA GPU RTX 6000. 

 \subsection{Inverse transport equation}
 In this subsection,  the transport equation is tested to verify the correctness of the proposed method. Consider the following equation:

 \begin{equation}
 \left\{
 \begin{split}
     &\frac{\partial u}{\partial t} + a\frac{\partial u}{\partial x} = 0, x\in [0,1], t\in [0,1],\\
     &u(x,0) = x.
\end{split} 
\right.
 \end{equation}
 Then, the true equation can be written as $u(x,t) = x - at$. By giving some observations on the boundary and interior, we aim to infer the posterior distribution of $a$. Specifically, we set $a_{true} = 1$, then we take $N_{b} = 30$ measurements from the boundary and initial conditions. Moreover, we take $N_{u} = 60$ interior measurements together with $N_{f} = 500$ residual points. The noises added to these measurements follow the Gaussian distribution $\mathcal{N}(0, 0.1^{2})$, that is $\sigma_{u} = \sigma_{b} = \sigma_{f} = 0.1$. The batch size for the residual points is 20.

\begin{figure}[htbp]
     \centering
    \includegraphics[width=0.4\linewidth]{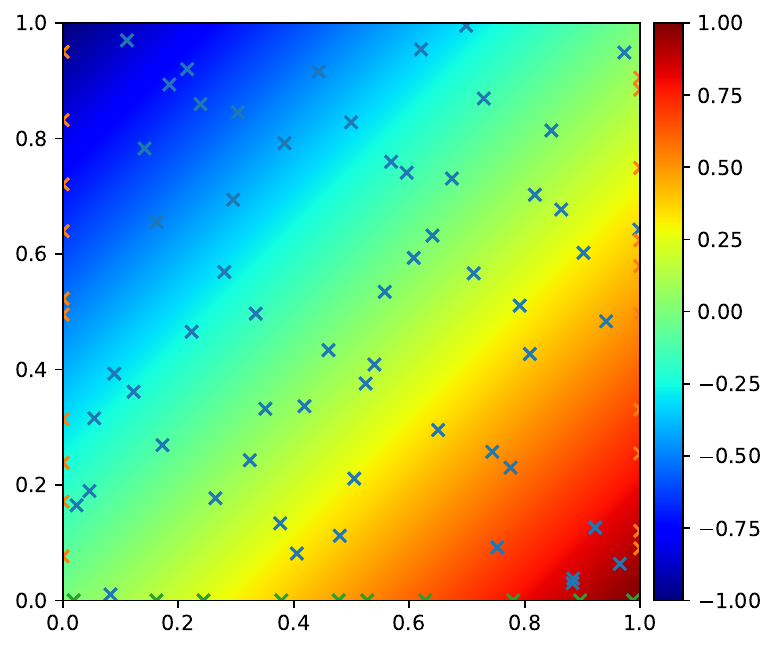}
     \caption{Transport equation: true solution field with the interior, boundary, and initial measurements. The noise scale is 0.1 for all data.}
     \label{transport_true}
 \end{figure}

 We assume that the prior is $p(a) = \mathcal{N}(0,1)$, and since the solution $u$ is linear on $a$, the posterior of $a$ can be written explicitly as 
 \begin{equation}
     p(a|\mathcal{D}_{u}, \mathcal{D}_{b})\propto \exp\left(-\frac{\sum_{i=1}^{N_{b}}(u_{b}^{i} - x_{b}^{i} + at_{b}^{i})^{2}}{2\sigma_{b}^{2}} - \frac{\sum_{i=1}^{N_{u}}(u_{u}^{i} - x_{u}^{i} + at_{u}^{i})^{2}}{2\sigma_{u}^{2}} - \frac{a^{2}}{2}\right),
 \end{equation}
 where $\{x_{b}^{i}, t_{b}^{i}, u_{b}^{i}\}_{i=1}^{N_{b}}$ and $\{x_{u}^{i}, t_{u}^{i}, u_{u}^{i}\}_{i=1}^{N_{u}}$ denote the measurement points located on the boundary and in the interior of the domain, respectively, as depicted in Fig. \ref{transport_true}. 
 
\begin{figure}[htbp]
     \centering
     \includegraphics[width=0.45\linewidth]{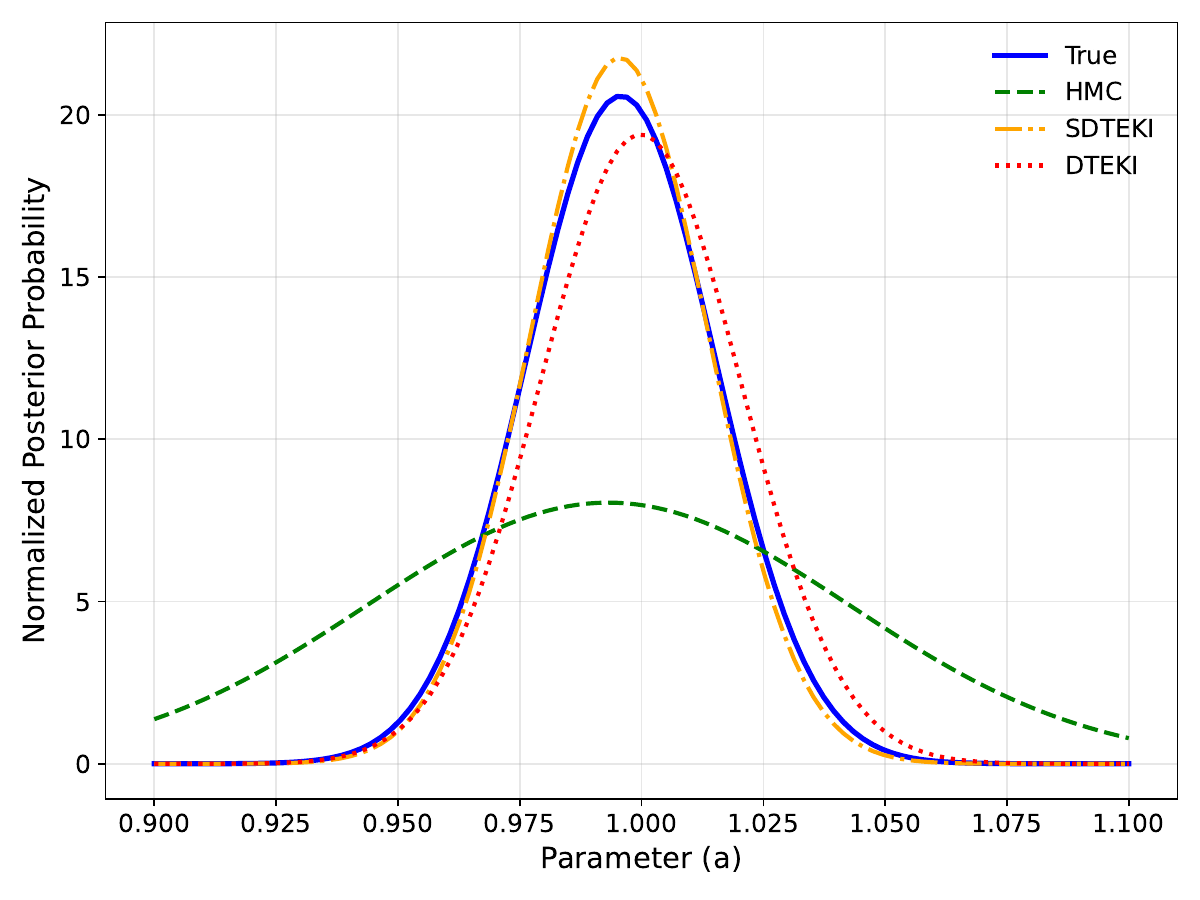}
     \caption{Transport equation: predicted posterior of $a$ obtained by different methods with cKANs.}
     \label{transport}
 \end{figure}

 Fig. \ref{transport} depicts the posterior distribution obtained by different methods compared to the true posterior.
The predicted posteriors obtained by HMC, DTEKI, and SDTEKI show little difference compared to the true posterior. However, our methods, DTEKI and SDTEKI, demonstrate better performance relative to HMC, which presents a much larger variance. This demonstrates that our method does not only achieve high accuracy but can provide more reliable uncertainty estimation. This improvement is likely attributed to the use of mini-batch training and the regularization as well as additional noises applied to the parameters, which enhance the precision of the predictions.

To further illustrate the efficiency of our methods, we present a comparison of the running time, along with the mean and standard deviation of the predictions, as summarized in Table \ref{tab:1}. It is clear that the SDTEKI outperforms all other methods not only in accuracy but also in efficiency. The active subspace reduces the number of network parameters and accelerates the training speed at the same time. Thus, overfitting is also alleviated with fewer parameters, thereby achieving better performance.

\begin{figure}
    \centering
    \begin{overpic}[width = 0.45\textwidth]{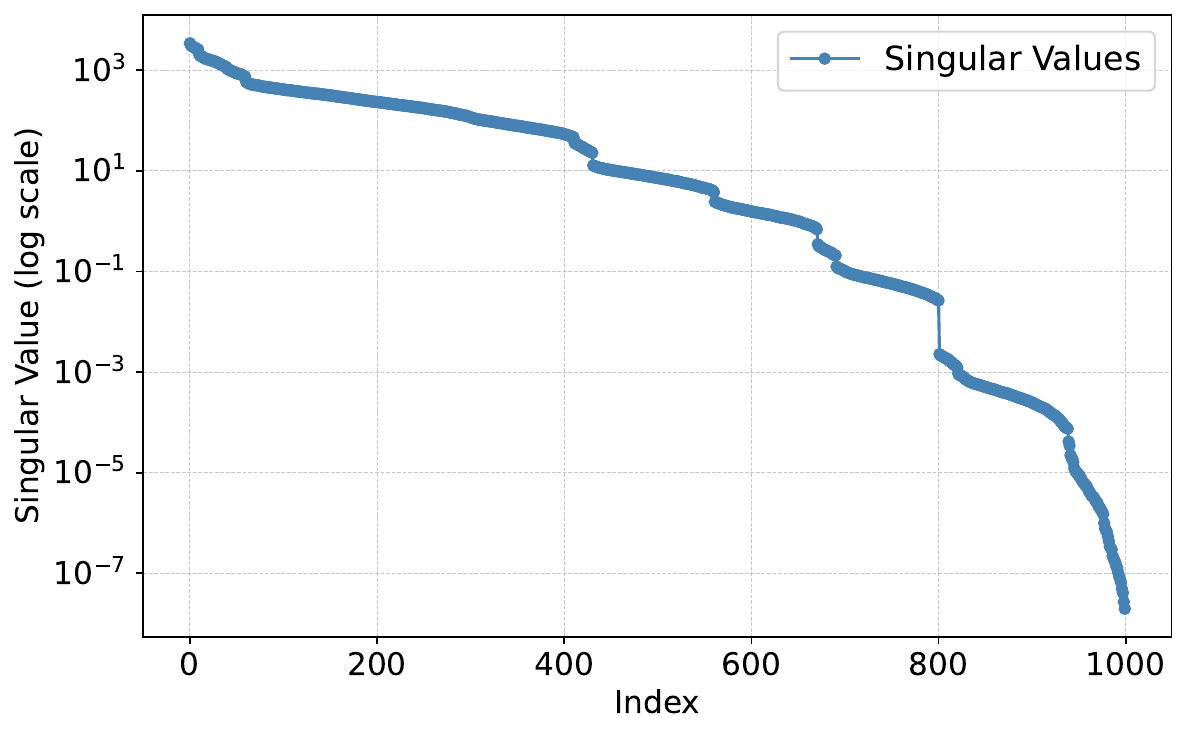}
    \end{overpic}
    \begin{overpic}[width = 0.45\textwidth]{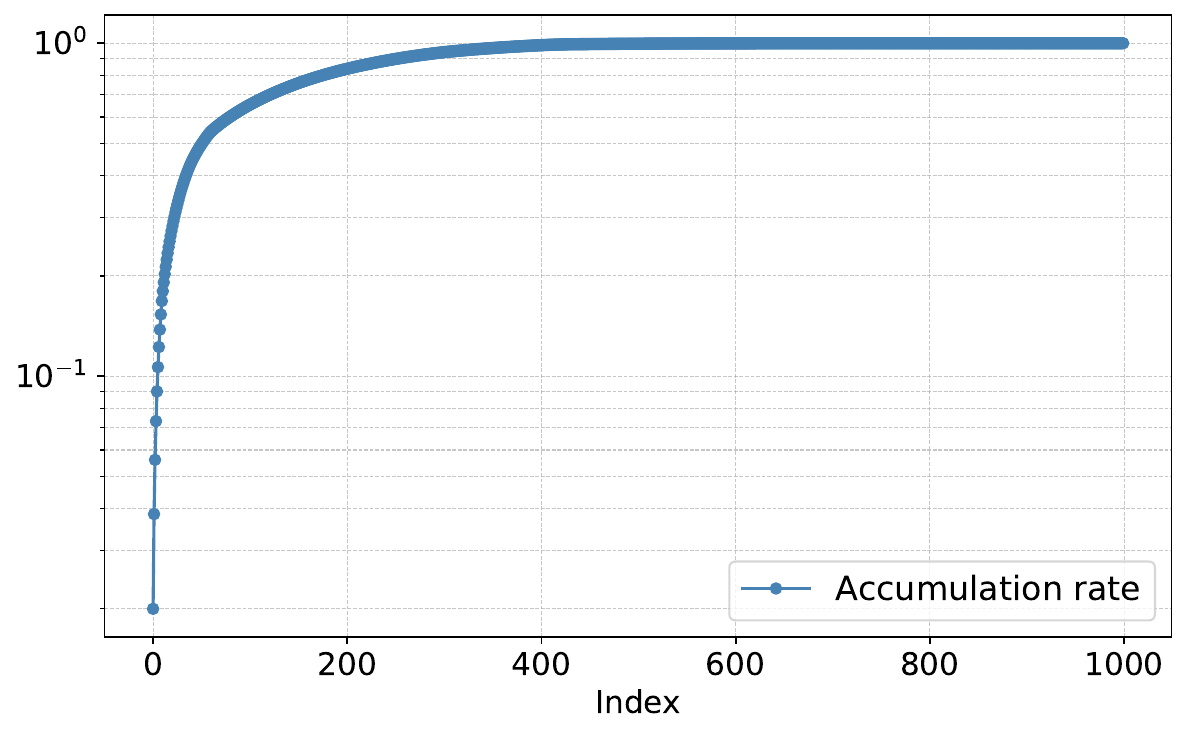}
    \end{overpic}
    \caption{Transport equation: the singular values and their accumulation rate for the decomposition.}
    \label{singular_transport}
\end{figure}

\begin{table}[t]
        \centering
        \begin{tabular}{ccccc}
\hline Method & $Mean$ & $Std$  & Walltime & Network size \\
\hline True & 0.9958& 0.0193&-&-\\ 
  HMC & $0.9932$ & $0.0432$ & 89.25s &1040  \\
DTEKI  & $0.9998$ &0.0206  & 9.38s &1040 \\
SDTEKI  & $0.9955 $ & 0.0183 & 7.53s & 347\\
\hline
\end{tabular}
        \caption{Transport equation: predicted mean and standard deviation for different methods. The running time (GPU RTX 6000) for different methods is also presented.}
        \label{tab:1}
    \end{table}
To demonstrate why active subspace works for cKANs, we plot the singular values of the decomposition and also their accumulation rate in Fig. \ref{singular_transport}. It is clear that the first one-third of the singular values account for more than $99\%$ of the total contribution, meaning that the bases of these singular values can explain the majority of the forward operator. Thus, it is reasonable that the active subspace method will work not only to accelerate training but also to reduce overfitting. 
\subsection{1D diffusion equation}
Consider the 1D diffusion equation:
\begin{equation}
    \kappa \frac{d^{2}u(x)}{dx^{2}} + D \frac{du(x)}{dx} = f(x), x\in[0,1],
\end{equation}
where the true solution is specified as $u(x) = \sin(6\pi x)\cos(4\pi x)^{2}.$ The right-hand-side $f(x)$ can be obtained by passing through the equation with $u(x)$. In this experiment, $\kappa = 0.001, D = 0.1$, which exhibits multi-scale behavior, and thus it is hard to solve. By assuming $D$ is unknown, the first goal of this problem is to infer $D$ and the solution $u(x)$ by having measurements on the boundary and interior. Specifically, we will generate $N_{u} = 6$ uniform interior points together with $N_{b} = 2$ boundary points as the observations of $u(x)$. Moreover, we will generate $N_{f} = 50$ residual points as observations of $f(x)$. The distributions of noise added to the observations are $N(0, 0.1^{2})$, that is $\sigma_{u} = \sigma_{b} = \sigma_{f} = 0.1$. 

\begin{figure}
            \centering 
            \begin{overpic}[width = 0.32\textwidth]{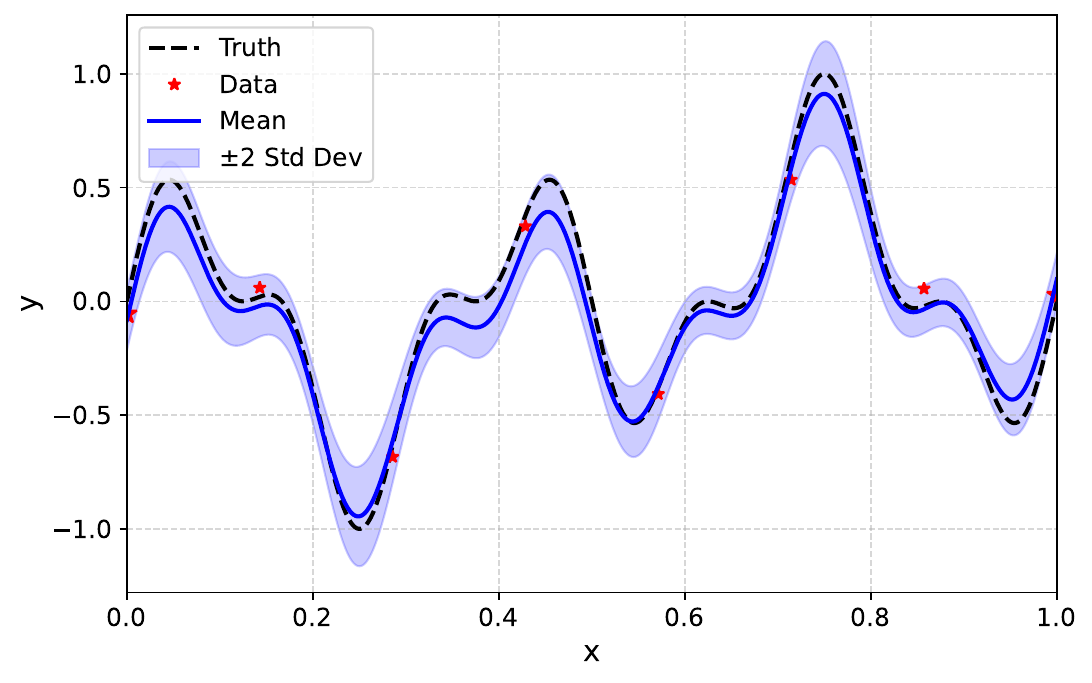}
            \put (40,63) {\textbf{\footnotesize HMC}}
            \end{overpic}
            \begin{overpic}[width = 0.32\textwidth]{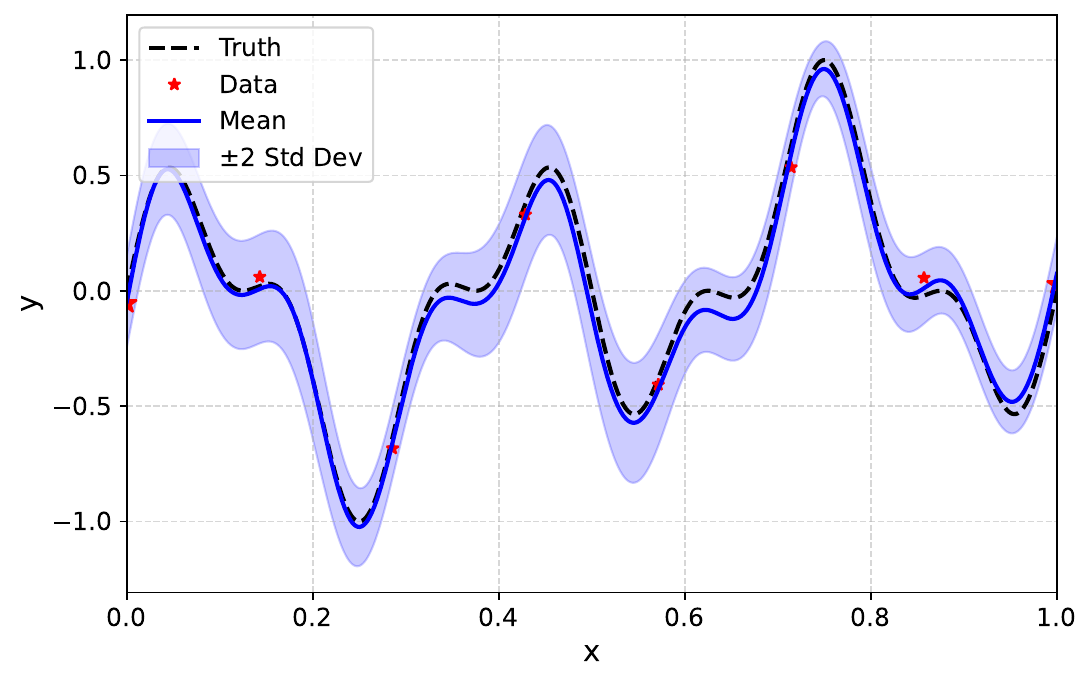}
            \put (39,63) {\textbf{\footnotesize DTEKI}}
            \end{overpic}
            \begin{overpic}[width = 0.32\textwidth]{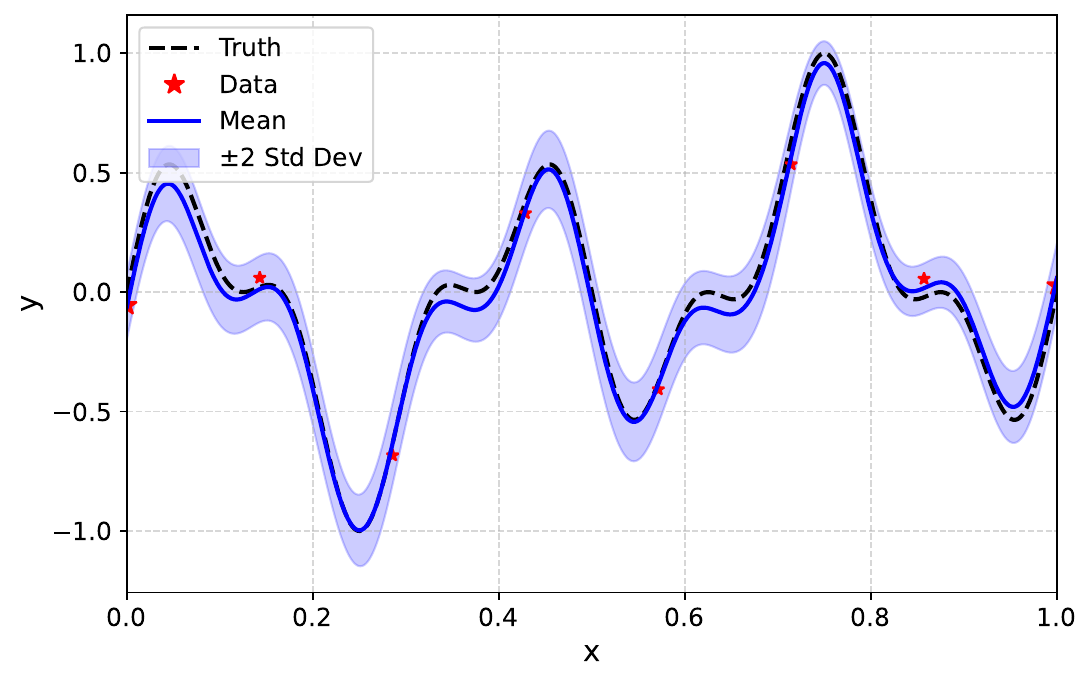}
            \put (39,63) {\textbf{\footnotesize SDTEKI}}
            \end{overpic}
            \begin{overpic}[width = 0.32\textwidth]{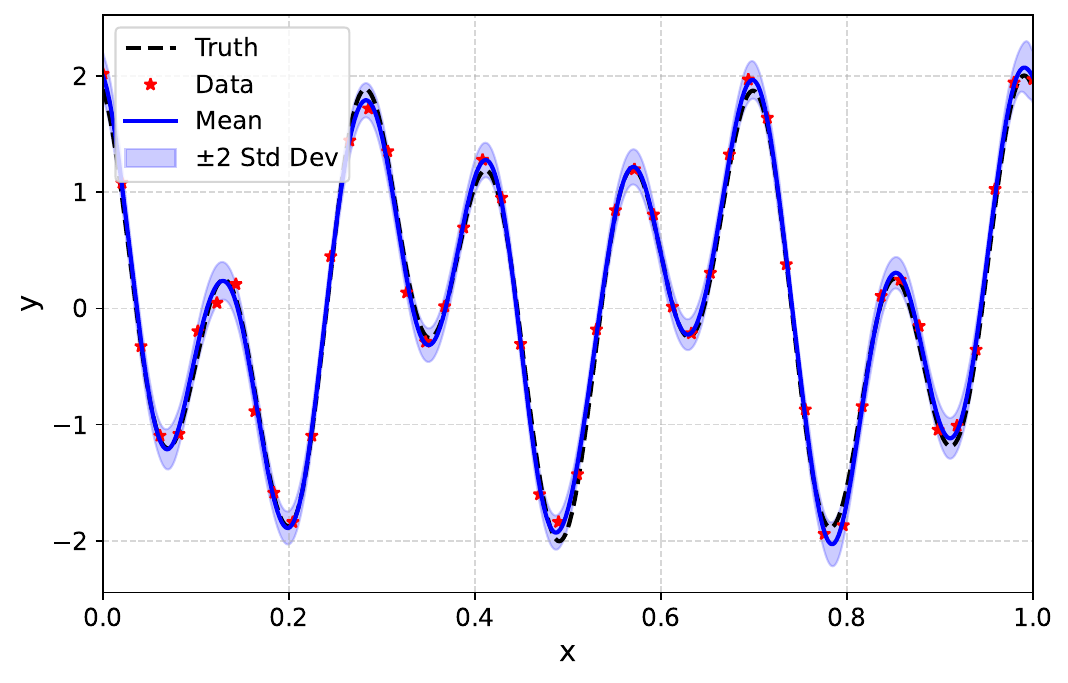}
            \end{overpic}
            \begin{overpic}[width = 0.32\textwidth]{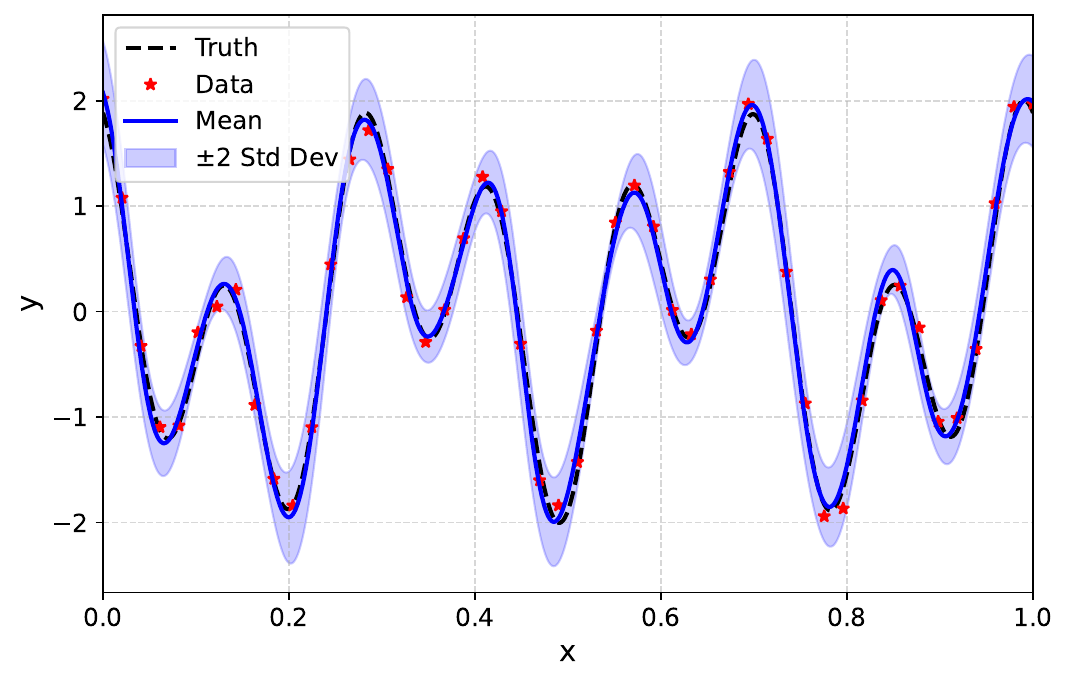}
            \end{overpic}
            \begin{overpic}[width = 0.32\textwidth]{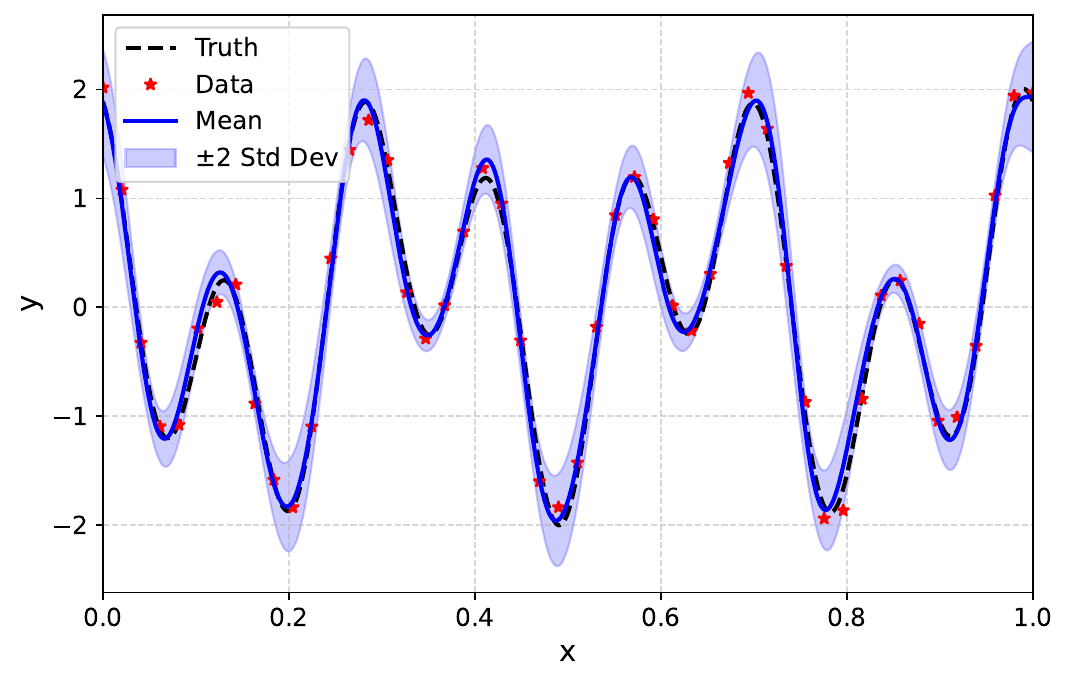}
            \end{overpic}
            \caption{Diffusion equation: noise scale is 0.1 for all data. Comparison of predicted mean $u(x)$ (first row) and predicted mean $f(x)$ (second row) for different methods.}
            \label{diffusion_prediction}
        \end{figure}

\begin{figure}
            \centering 
            \begin{overpic}[width = 0.32\textwidth]{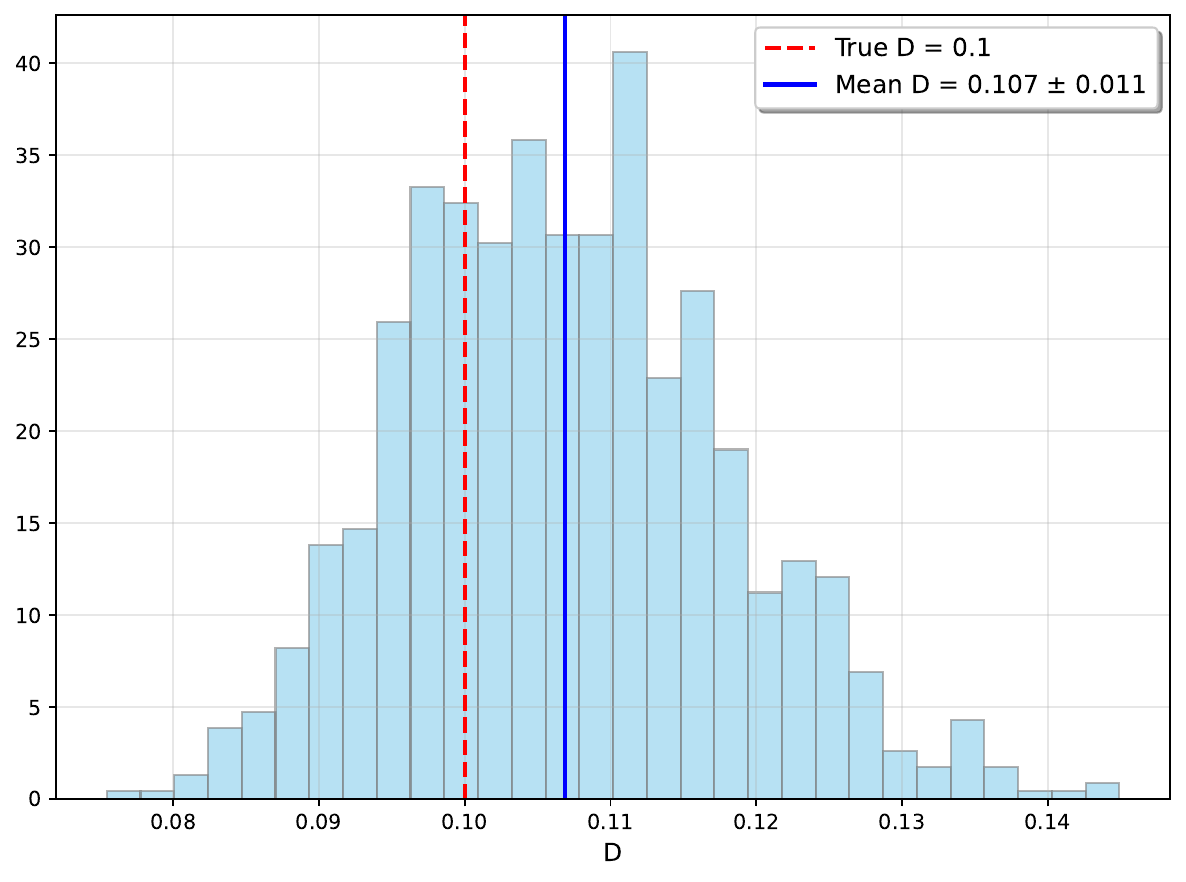}
            \put (35,75) {\textbf{\footnotesize HMC}}
            \end{overpic}
            \begin{overpic}[width = 0.32\textwidth]{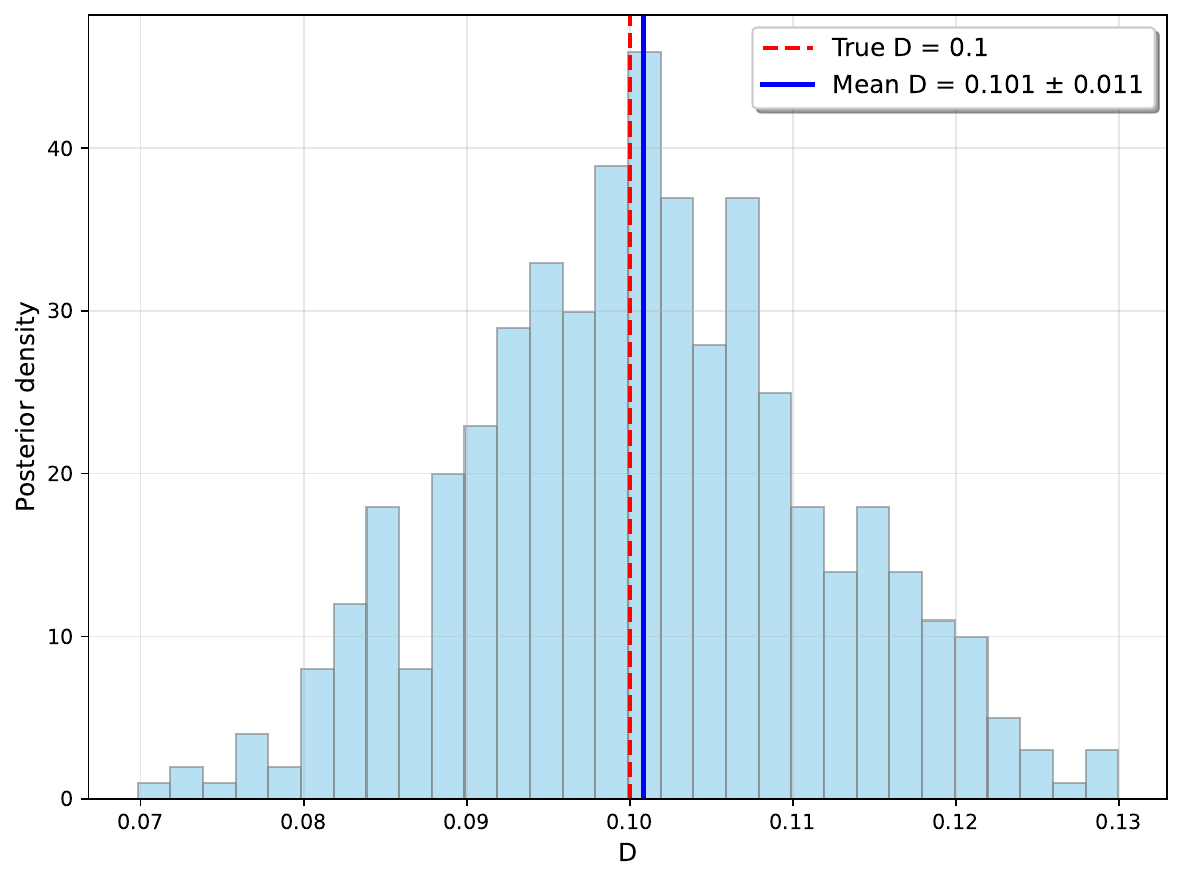}
            \put (35,75) {\textbf{\footnotesize DTEKI}}
            \end{overpic}
            \begin{overpic}[width = 0.32\textwidth]{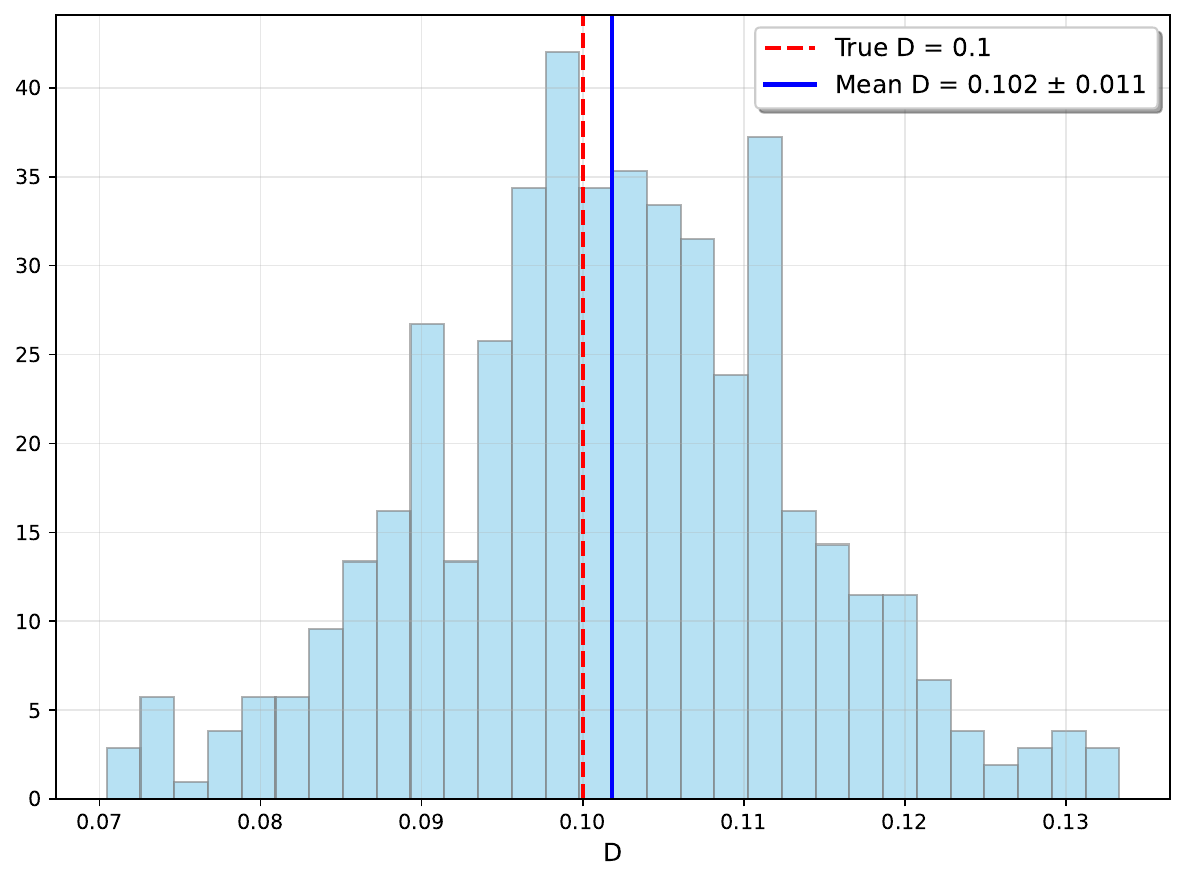}
            \put (35,75) {\textbf{\footnotesize SDTEKI}}
            \end{overpic}
            \caption{Diffusion equation: histogram of the prediction values for $D$.}
        \end{figure}
Fig. \ref{diffusion_prediction} illustrates the predicted values of $u(x)$ and 
$f(x)$ obtained using various methods. All approaches demonstrate high accuracy in predicting the solution, as reflected by the small $L_{2}$ errors reported in Table \ref{tab:2}. However, HMC exhibits an underestimated variance near 
$x=0.4$, where the confidence interval fails to encompass the true solution. This discrepancy is likely due to overfitting, as HMC optimizes a larger set of parameters without sufficient regularization. In contrast, DTEKI and SDTEKI provide more reliable uncertainty estimates across the entire domain. Additionally, SDTEKI, with its reduced number of parameters, effectively mitigates overfitting while maintaining accuracy and delivering well-calibrated confidence intervals.

For the posterior prediction of $D$, DTEKI and SDTEKI achieve superior performance, with relative errors of only $0.85\%$ and $1.70\%$, respectively. This greatly surpasses the performance of HMC, with a $L_{2}$ error of $6.88\%$, which demonstrates that our method can handle such multiscale problems with high accuracy and reliable uncertainty quantification.

The efficiency of our methods is further demonstrated in Table \ref{tab:2}. Specifically, DTEKI achieves an approximately 12-fold speed-up compared to HMC. SDTEKI, with its streamlined parameter set, is even more computationally efficient while maintaining comparable accuracy.

\begin{table}[t]
        \centering
        \begin{tabular}{ccccc}
\hline Method & $e_{u}$ & $e_{D}$  & Walltime & Network size \\
\hline  
  HMC & $12.74\%$ & $6.88\%$ & 70.85s &960  \\
DTEKI  & $11.55\%$ &$0.85\%$  & 6.03s &960 \\
SDTEKI  & $11.39\%$ & $1.70\%$ & 5.06s & 320\\
\hline
\end{tabular}
        \caption{Diffusion equation: $L_{2}$ error of predictions with respect to different methods. The running time (GPU RTX 6000) for different methods is also presented.}
        \label{tab:2}
    \end{table}  

Next, we test the ability of our proposed method to deal with large data corrupted by large noise. To this end, we generate $N_{f} = 5000$ residual measurements corrupted by white noise with noise scale $\sigma_{f} = 2$. At the same time, we still use $N_{b} = 2$ boundary measurements with noise scale $\sigma_{u} = 0.05$. Moreover, the true value of parameter $D$ is assumed to be known. The batch size is set to 100. 

\begin{figure}
            \centering 
            \begin{overpic}[width = 0.32\textwidth]{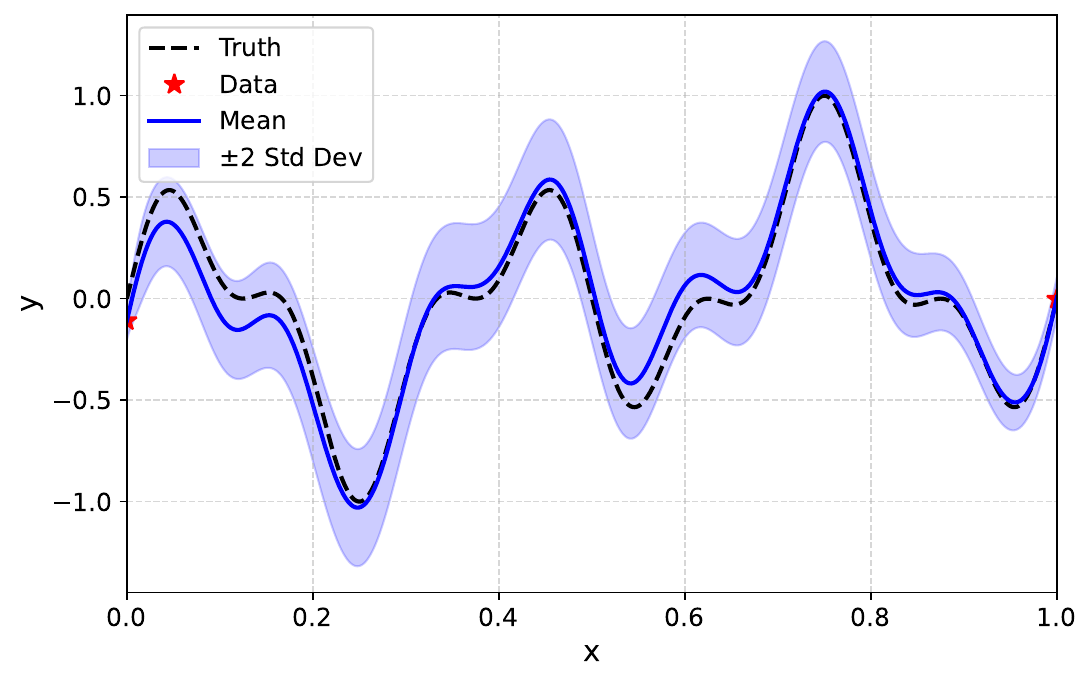}
            \put (40,63) {\textbf{\footnotesize HMC}}
            \end{overpic}
            \begin{overpic}[width = 0.32\textwidth]{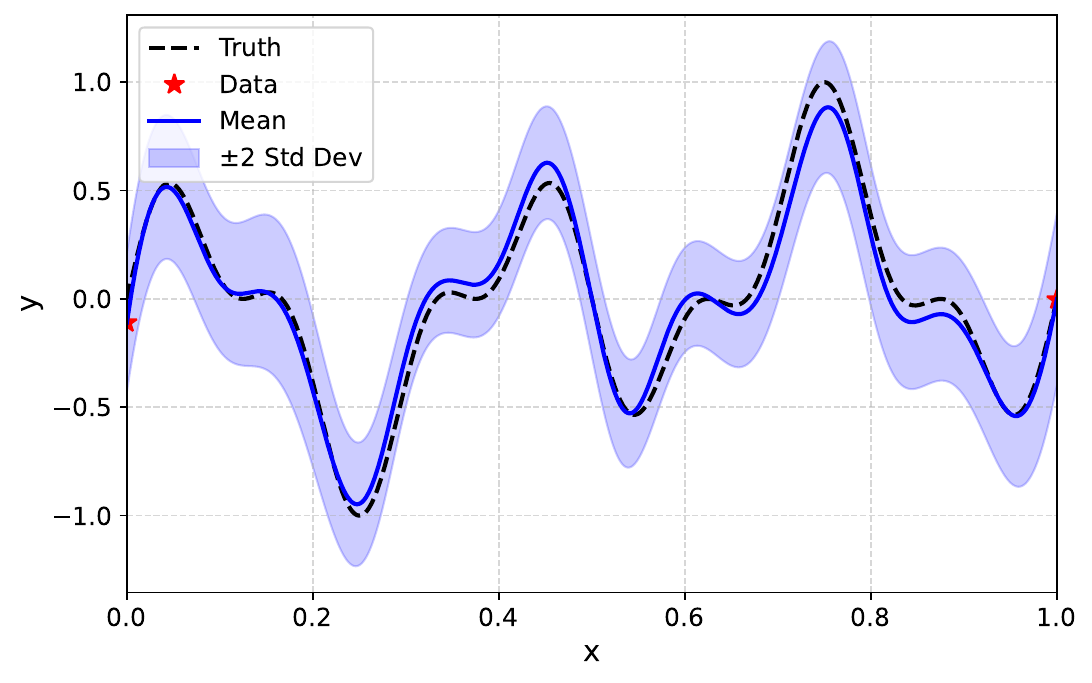}
            \put (35,63) {\textbf{\footnotesize DTEKI}}
            \end{overpic}
            \begin{overpic}[width = 0.32\textwidth]{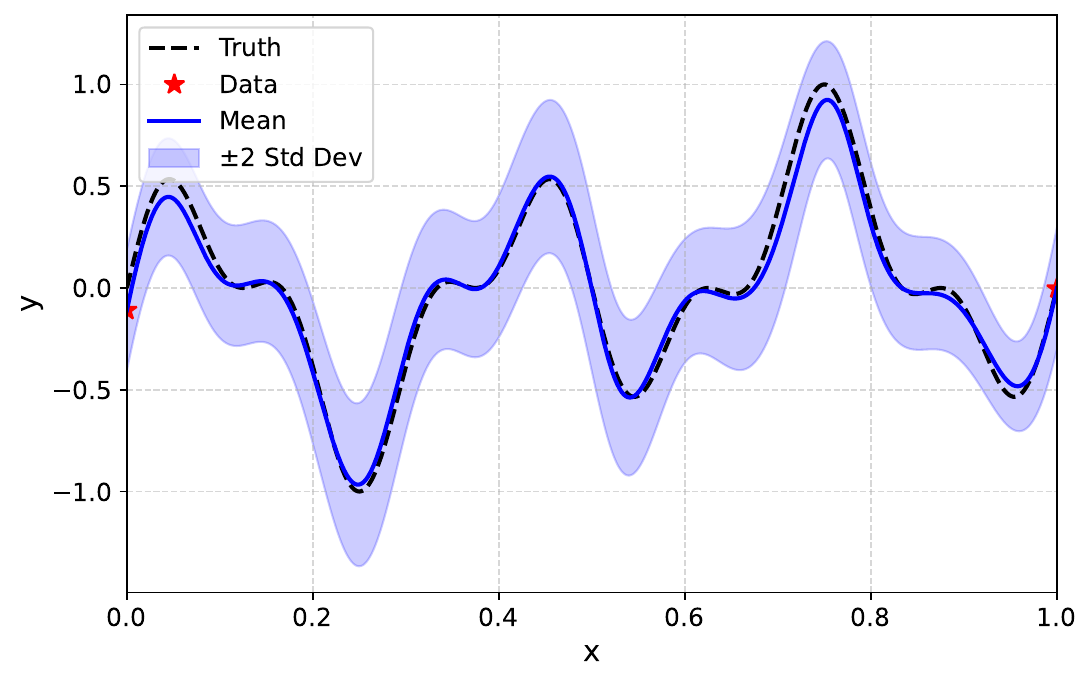}
            \put (35,63) {\textbf{\footnotesize SDTEKI}}
            \end{overpic}
            \begin{overpic}[width = 0.32\textwidth]{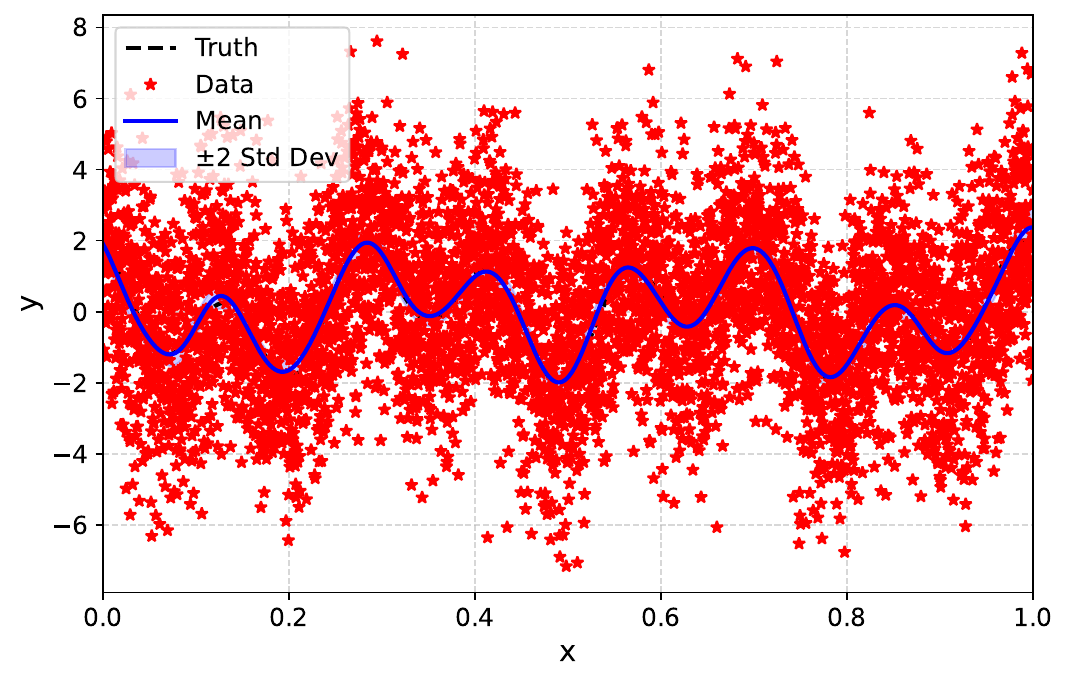}
            \end{overpic}
            \begin{overpic}[width = 0.32\textwidth]{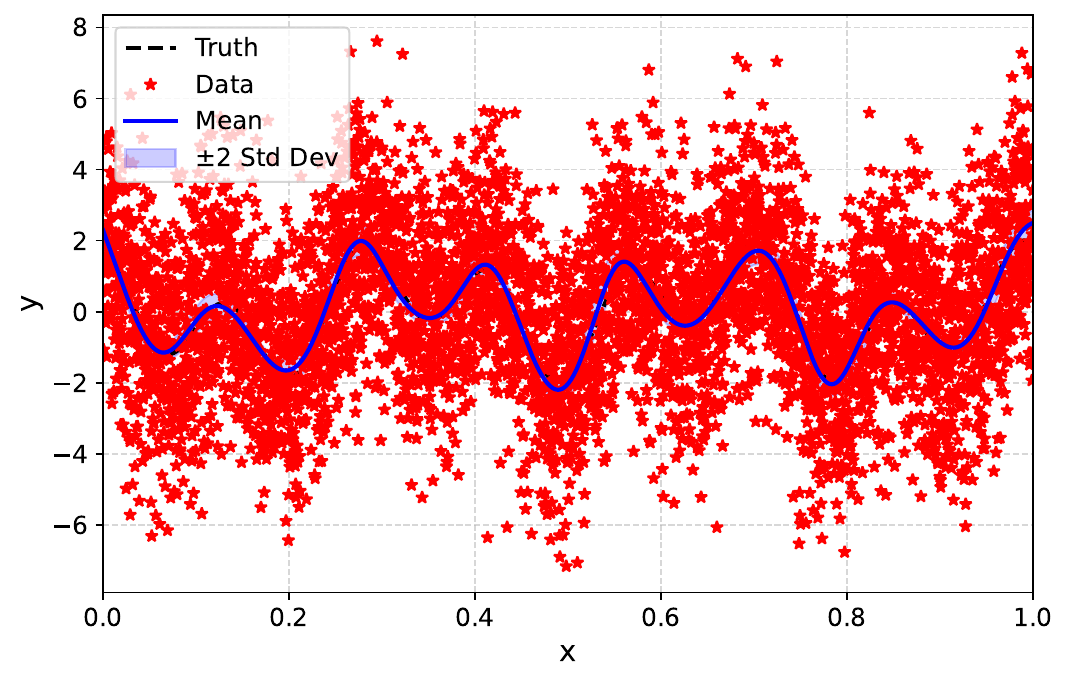}
            \end{overpic}
            \begin{overpic}[width = 0.32\textwidth]{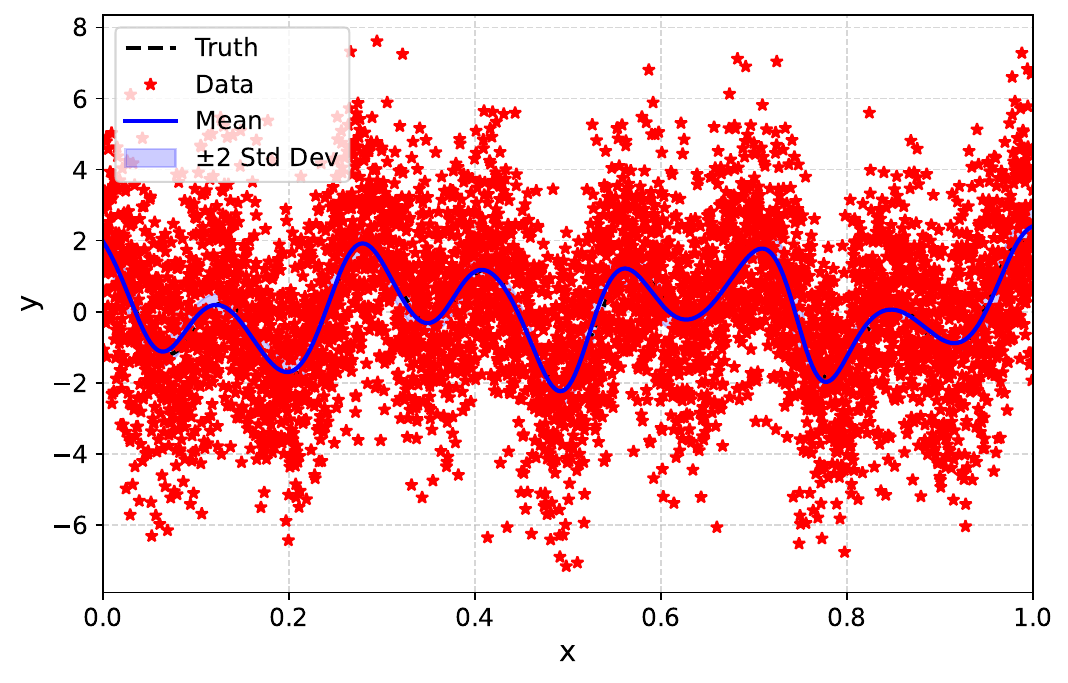}
            \end{overpic}
            \caption{Diffusion equation: noise scale $\sigma_{f} = 2$ for $f(x)$ and $\sigma_{b} = 0.05$ for boundary data. Comparison of prediction values of $u(x)$ (first row) and $f(x)$ (second row) obtained by different methods.}
            \label{diffusion_large_noise}
        \end{figure}

The results in Fig. \ref{diffusion_large_noise} illustrate that with big data corrupted by large noise levels, HMC could easily overfit with less accurate confidence intervals due to the dominance of likelihoods in such cases. On the contrary, by applying mini-batch training strategies and reducing parameters using active subspace, our methods, DTEKI and SDTEKI can both capture the information provided by the prior and likelihoods, leading to more accurate predictions and more reliable uncertainty estimation. Finally, the $L_{2}$ error of the predicted solutions are $16.76\%$ and $13.87\%$ for DTEKI and SDTEKI respectively, which outperforms the $L_{2}$ error of $20.47\%$ obtained by HMC according to Table \ref{tab:3}. Moreover, in terms of computing efficiency, our methods achieve over an order speed-up compared to HMC while preserving the prediction accuracy due to the mini-batching training and gradient-free inference. This phenomenon indicates that our method can have the potential to deal with large-scale problems with much lower computational cost.

\begin{table}[t]
        \centering
        \begin{tabular}{ccccc}
\hline Method & $e_{u}$   & Walltime & Network size \\
\hline  
  HMC & $20.47\%$  & 216.08s &960  \\
DTEKI  & $16.76\%$  & 12.56s &960 \\
SDTEKI  & $13.87\%$  & 8.87s & 320\\
\hline
\end{tabular}
        \caption{Diffusion equation: $L_{2}$ error of predictions with respect to different methods. The running time (GPU RTX 6000) and network sizes for different methods are also presented.}
        \label{tab:3}
    \end{table}  

\subsection{1D nonlinear equation}
Next, we consider the following nonlinear 1D problem:
\begin{equation}
    \lambda \frac{d^{2}u(x)}{dx^{2}} + k\tanh(u(x)) = f(x), x\in [-0.7, 0.7].
\end{equation}
The true solution is specified as $u(x) = \sin(6x)^{3}$, where the right-hand side $f(x)$ can be obtained by deriving the partial differential operator. In this experiment, $\lambda = 0.01, k = 0.7$. This nonlinear problem is used to test the ability of our method when the PDE model is more complicated. To this end, we first assume $k$ is unknown and aim to infer this parameter by having some measurements. Specifically, we will generate $N_{b} = 2$ boundary measurements and also $N_{u} = 6$ interior measurements. Also, $N_{f} = 50$ residual measurements are used to fit the model. The noises added to these measurements are Gaussian with mean zeros and deviations $\sigma_{u} = \sigma_{b} = \sigma_{f} = 0.1$. 

\begin{figure}
            \centering 
            \begin{overpic}[width = 0.32\textwidth]{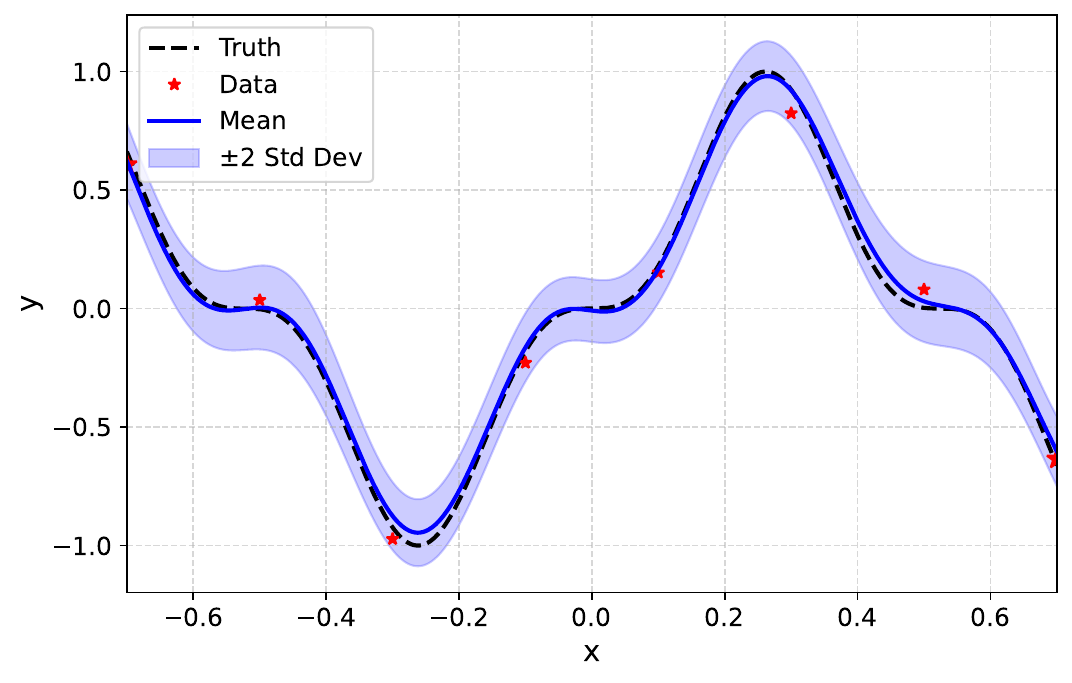}
            \put (40,65) {\textbf{\footnotesize HMC}}
            \end{overpic}
            \begin{overpic}[width = 0.32\textwidth]{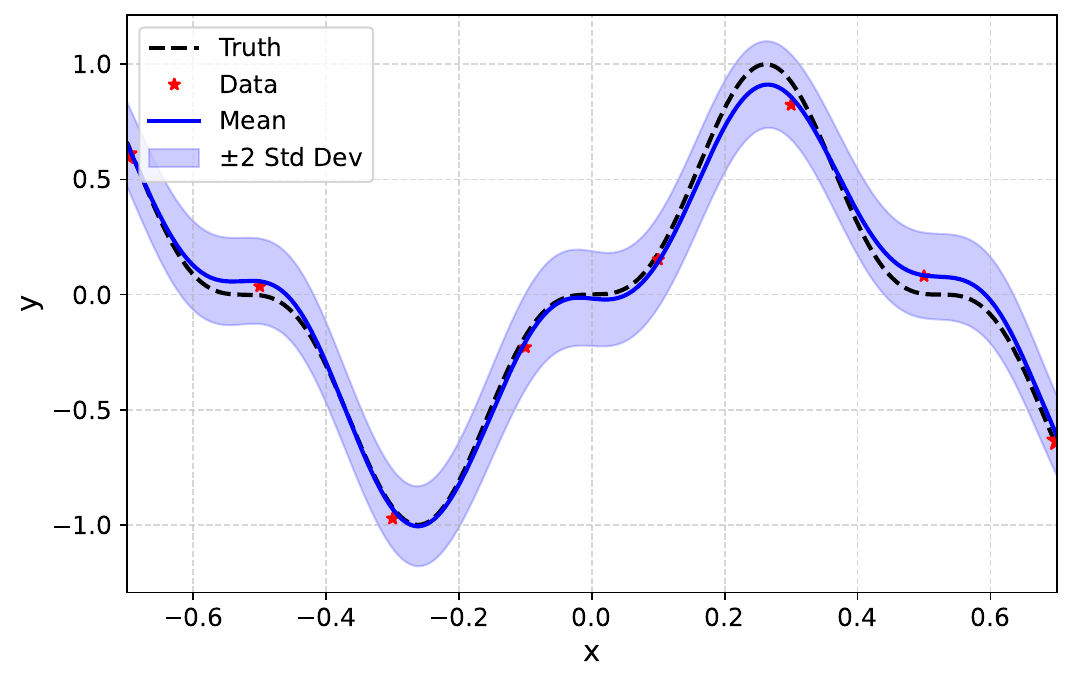}
            \put (35,65) {\textbf{\footnotesize DTEKI}}
            \end{overpic}
            \begin{overpic}[width = 0.32\textwidth]{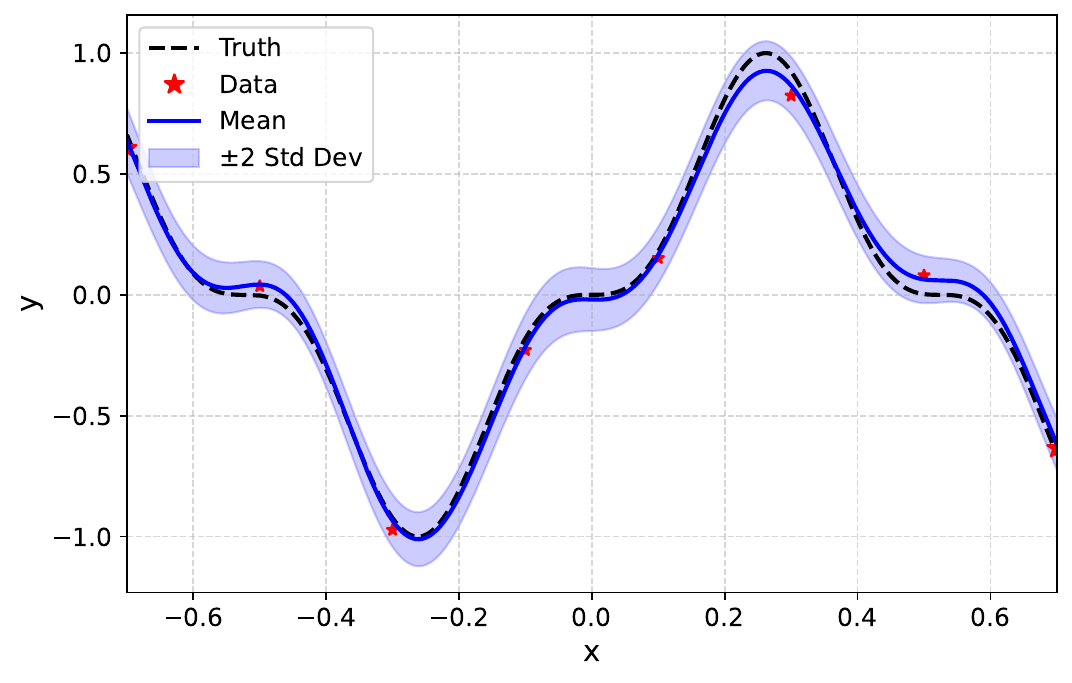}
            \put (35,65) {\textbf{\footnotesize SDTEKI}}
            \end{overpic}
            \begin{overpic}[width = 0.32\textwidth]{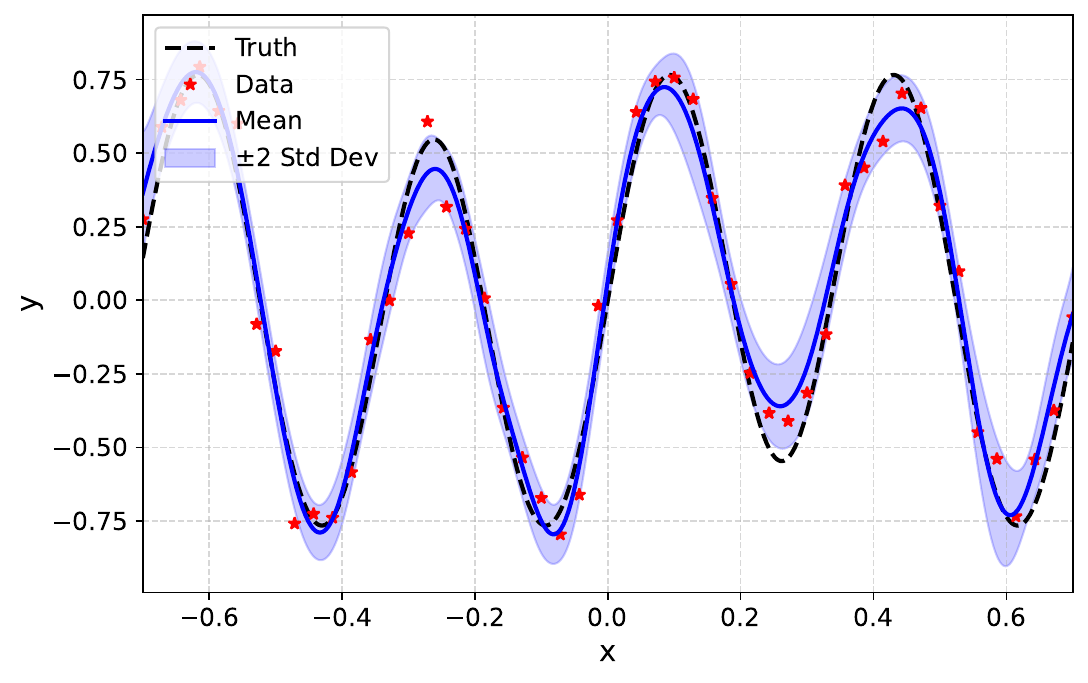}
            \end{overpic}
            \begin{overpic}[width = 0.32\textwidth]{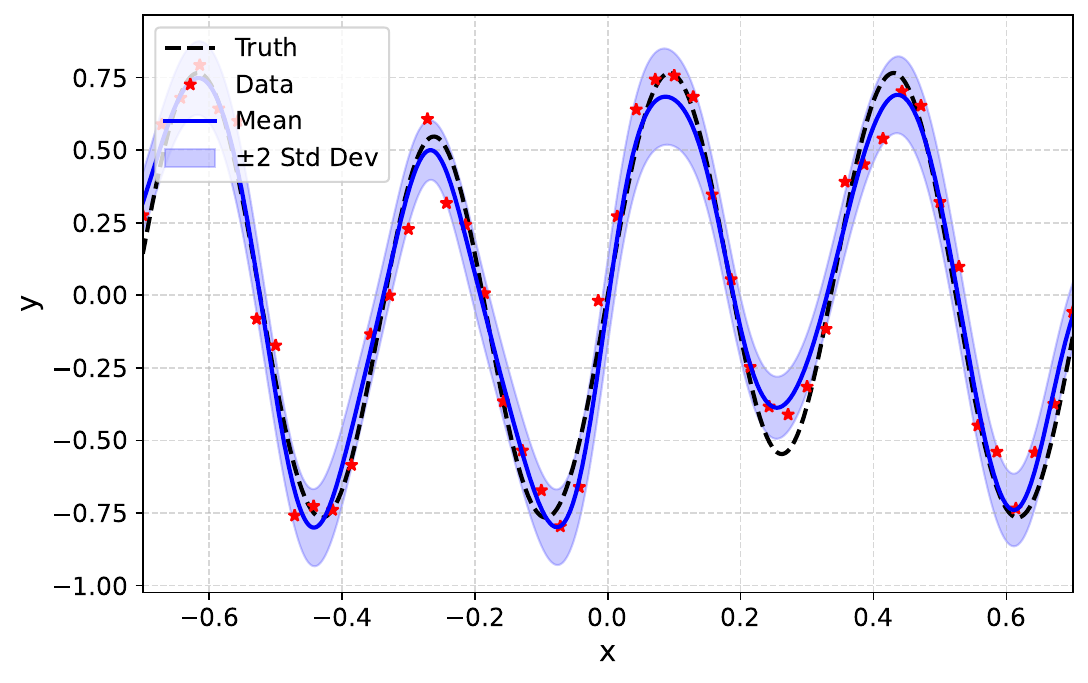}
            \end{overpic}
            \begin{overpic}[width = 0.32\textwidth]{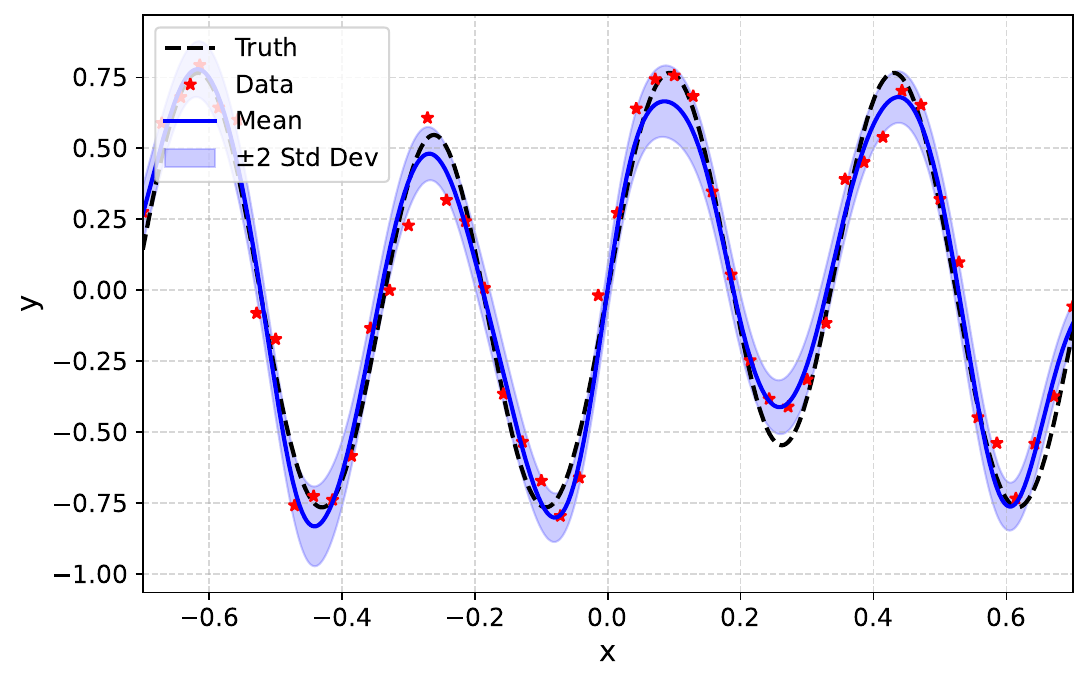}
            \end{overpic}
            \caption{Nonlinear equation: noise scale is 0.1 for all data. Comparison of $u(x)$ (first row) and $f(x)$ (second row) for different methods.}
            \label{nonlinear_prediction}
        \end{figure}

From Fig. \ref{nonlinear_prediction}, it is clear that all methods can achieve good performance with high accuracy. SDTEKI uses fewer parameters by looking for the active subspace in the parameter space. Thus, the overfitting phenomenon is reduced by showing high accuracy and a more accurate confidence interval compared to HMC. Finally, the $L_{2}$ prediction errors of $u(x)$ for these three methods are $5.95\%, 9.54\%, 7.63\%$ respectively as depicted in Table \ref{tab:4}. Furthermore, we can see from Table \ref{tab:4} that $SDTEKI$ is much more efficient than HMC with less running time, which demonstrates that our method can still obtain an accurate prediction with high efficiency for nonlinear problems. 
\begin{table}[t]
        \centering
        \begin{tabular}{ccccc}
\hline Method & $e_{u}$ &$e_{k}$   & Walltime & Network size \\
\hline  
  HMC & $5.95\%$ &$6.11\%$  & 70.62s &960  \\
DTEKI  & $9.54\%$ &$4.83\%$  & 6.36s &960 \\
SDTEKI  & $7.63\%$ &$2.84\%$  & 5.31s & 320\\
\hline
\end{tabular}
        \caption{Nonlinear equation: $L_{2}$ error of predictions with respect to different methods for the inverse problem. The running time (GPU RTX 6000) and network sizes for different methods are also presented.}
        \label{tab:4}
    \end{table}  

\begin{figure}[htbp]
            \centering 
            \begin{overpic}[width = 0.32\textwidth]{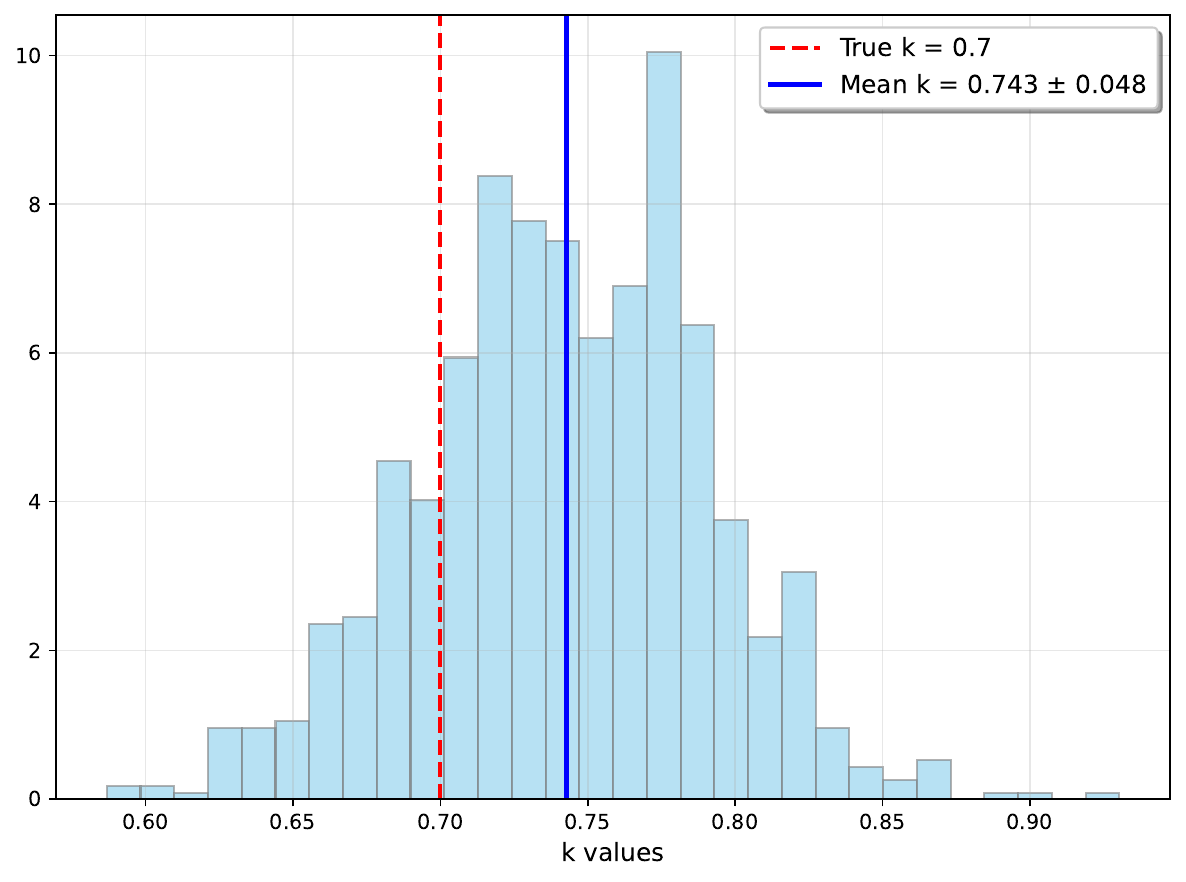}
            \put (35,77) {\textbf{\footnotesize HMC}}
            \end{overpic}
            \begin{overpic}[width = 0.32\textwidth]{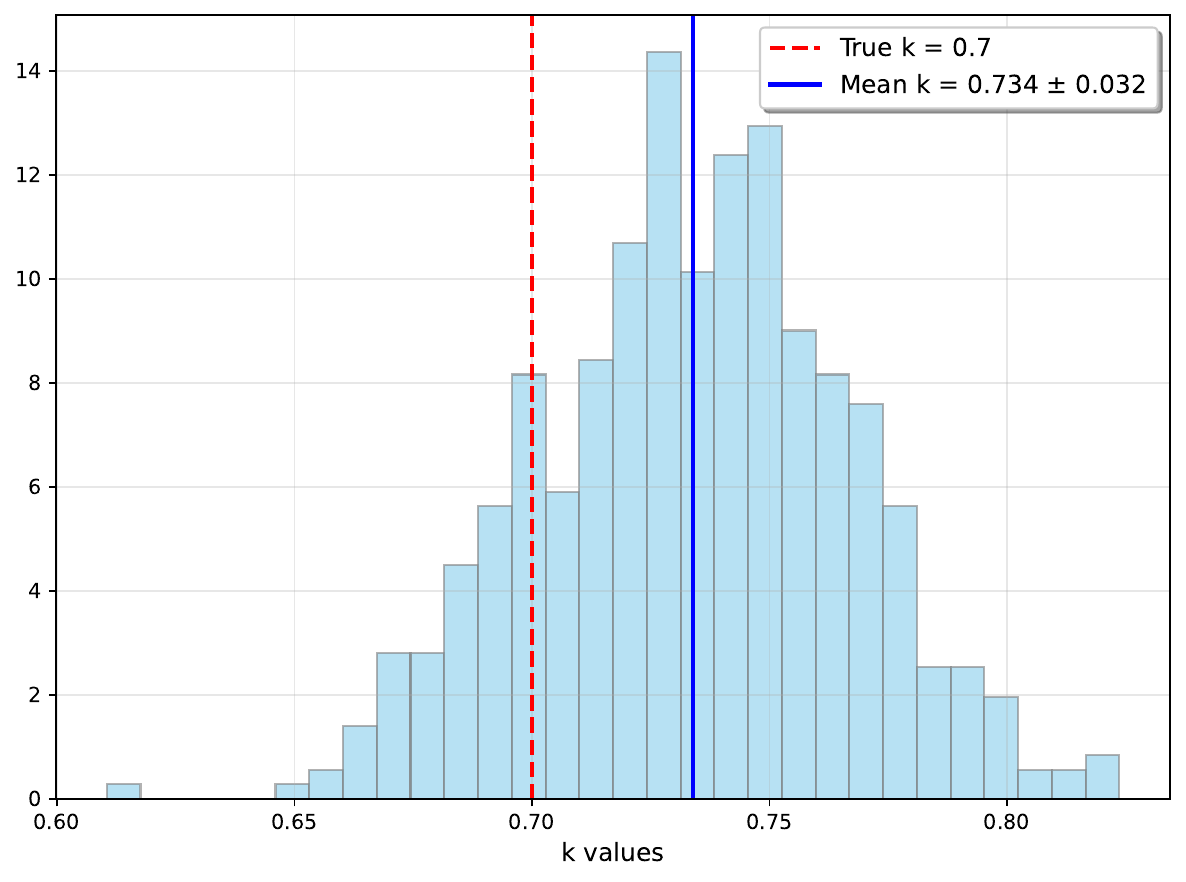}
            \put (35,77) {\textbf{\footnotesize DTEKI}}
            \end{overpic}
            \begin{overpic}[width = 0.32\textwidth]{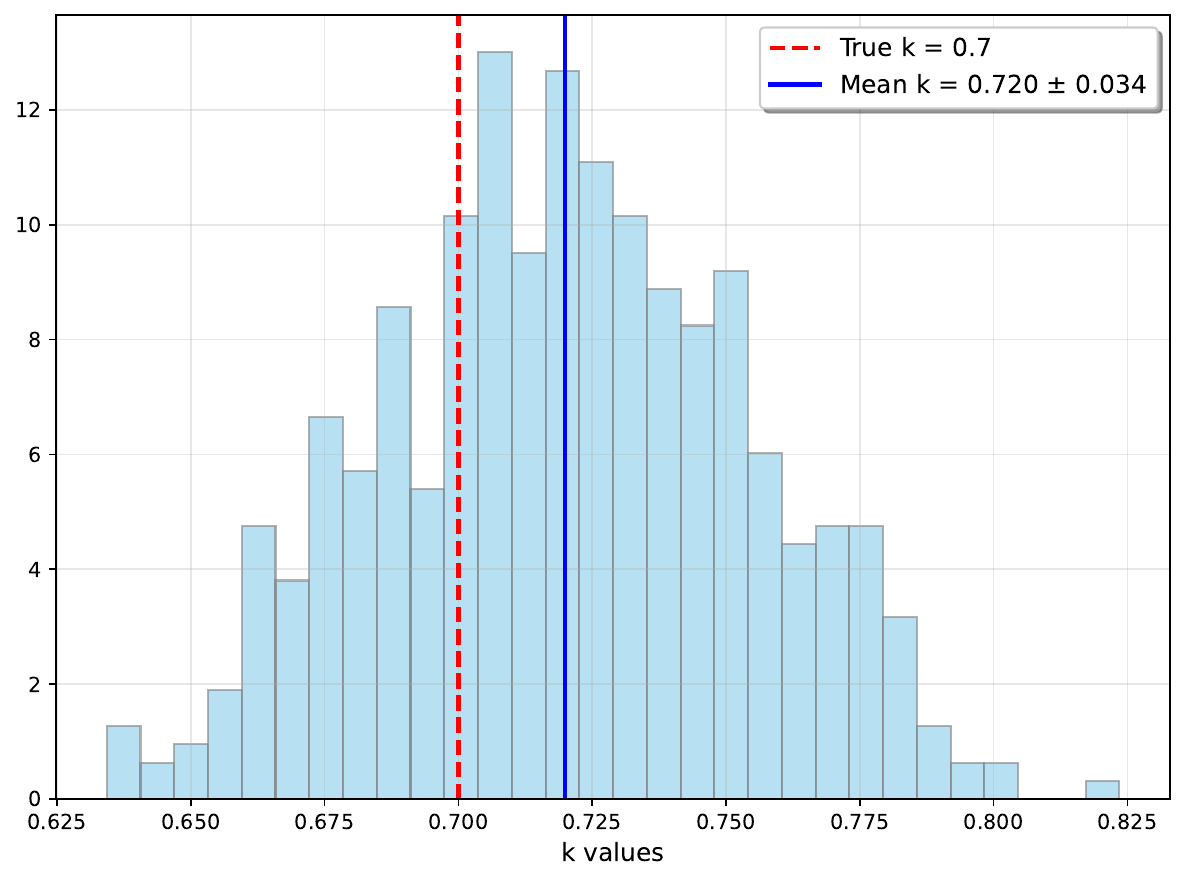}
            \put (35,77) {\textbf{\footnotesize SDTEKI}}
            \end{overpic}
            \caption{Nonlinear equation: histogram plot of prediction values for $k$ with different methods.}
            \label{estimate_k}
        \end{figure}

Regarding the posterior prediction of $k$ as depicted in Fig. \ref{estimate_k}, SDTEKI still has the most accurate prediction, with a relative $L_{2}$ error $2.84\%$ compared to $4.83\%$ and $6.1\%$ for DTEKI and HMC, respectively. This demonstrates that our method can even outperform the gold-standard HMC method with higher computing efficiency. Moreover, the prediction is reliable due to the small deviation as reflected in the figure. 

\begin{figure}[htbp]
            \centering 
            \begin{overpic}[width = 0.32\textwidth]{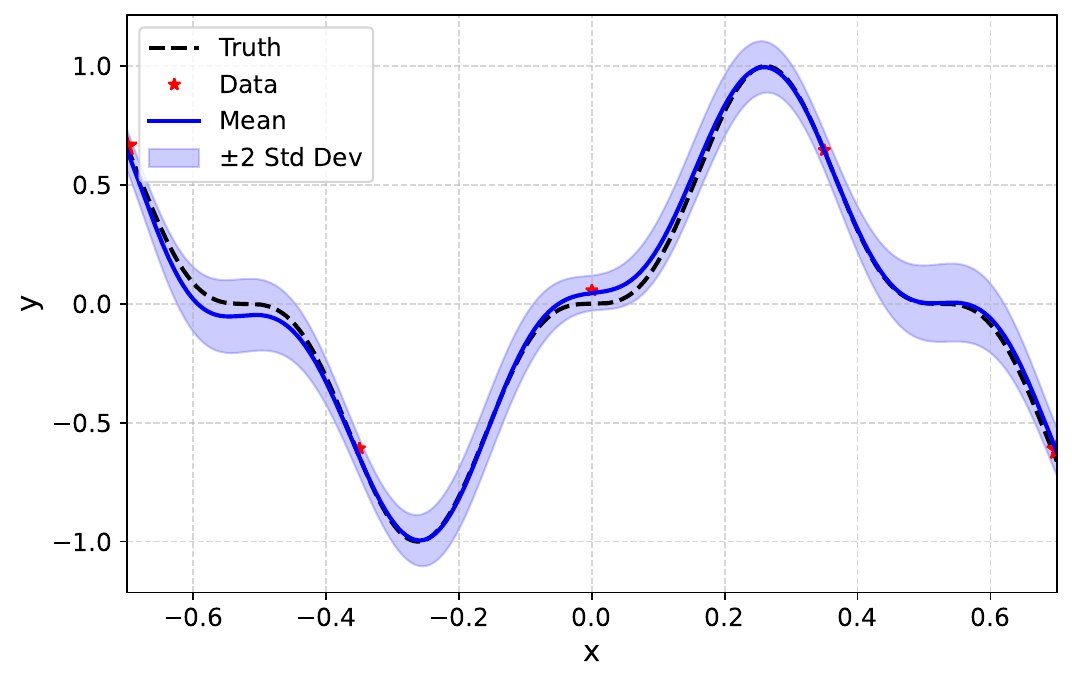}
            \put (40,65) {\textbf{\footnotesize HMC}}
            \end{overpic}
            \begin{overpic}[width = 0.32\textwidth]{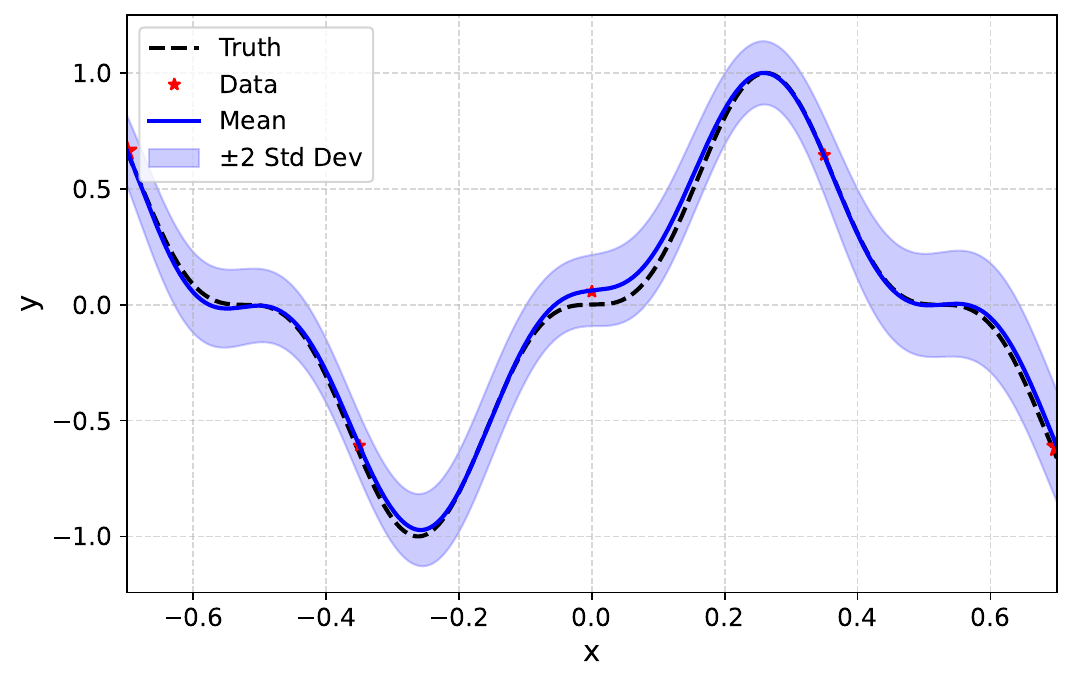}
            \put (35,65) {\textbf{\footnotesize DTEKI}}
            \end{overpic}
            \begin{overpic}[width = 0.32\textwidth]{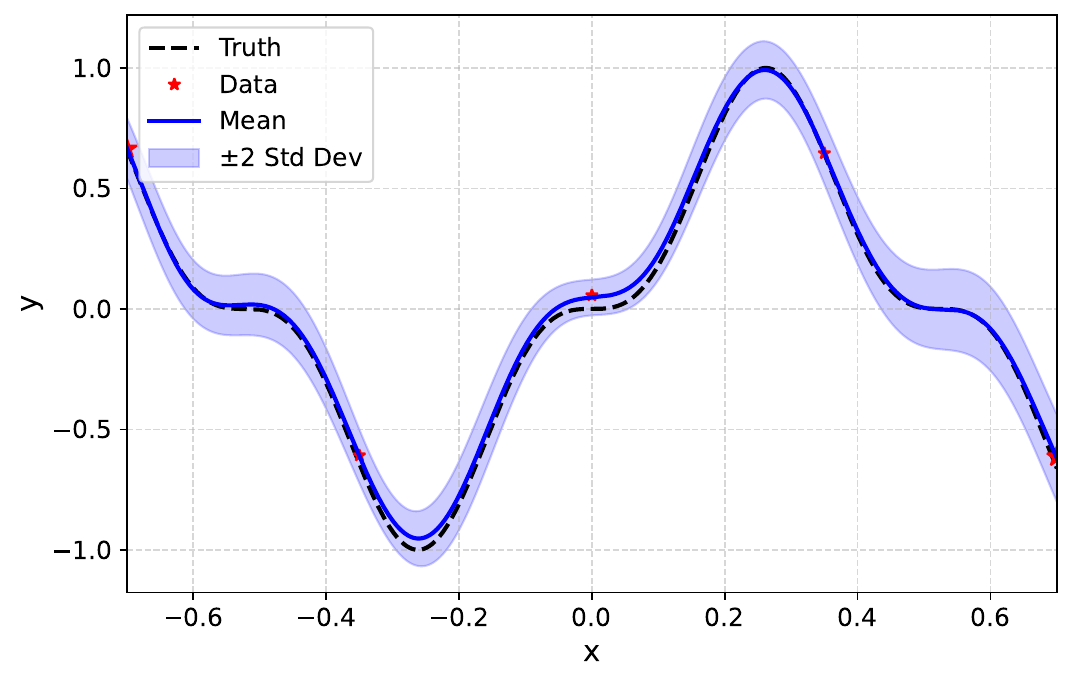}
            \put (35,65) {\textbf{\footnotesize SDTEKI}}
            \end{overpic}
            \begin{overpic}[width = 0.32\textwidth]{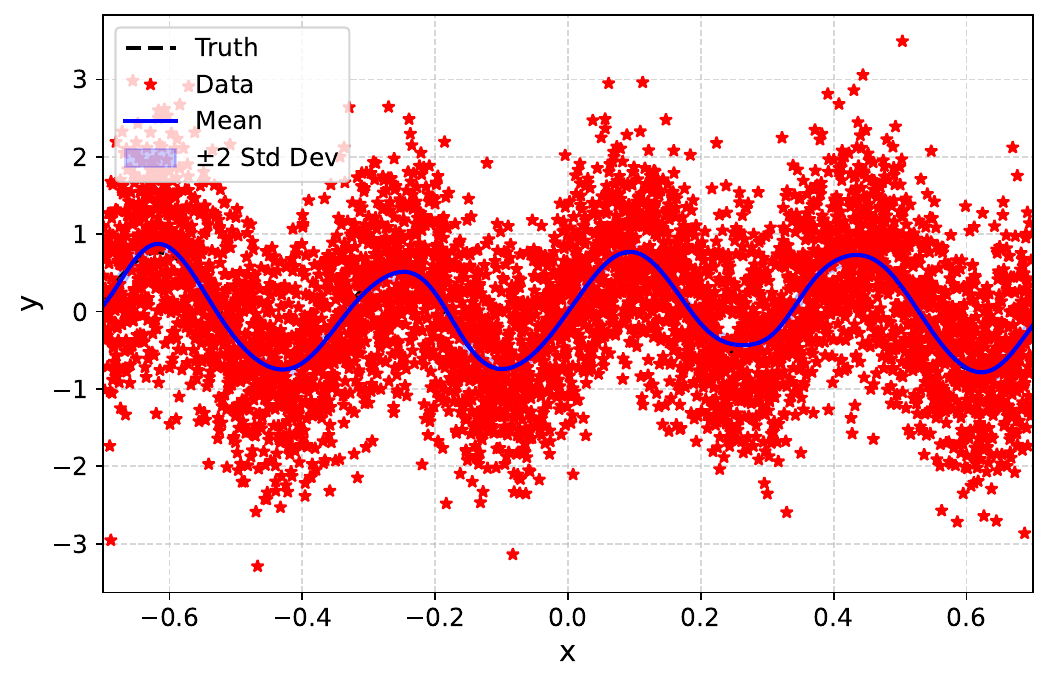}
            \end{overpic}
            \begin{overpic}[width = 0.32\textwidth]{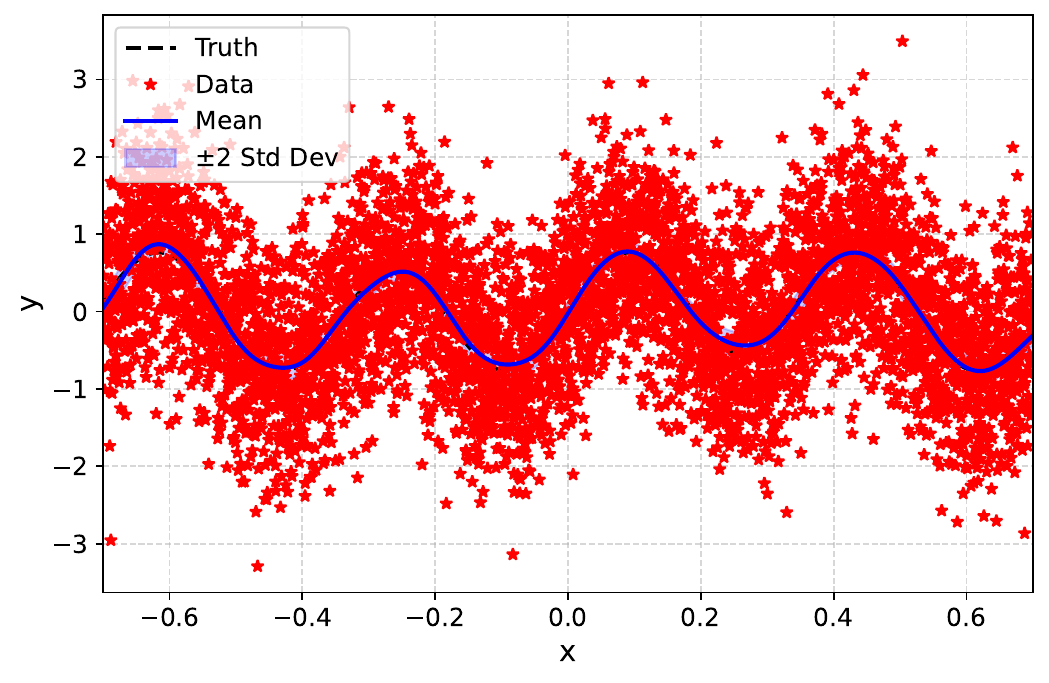}
            \end{overpic}
            \begin{overpic}[width = 0.32\textwidth]{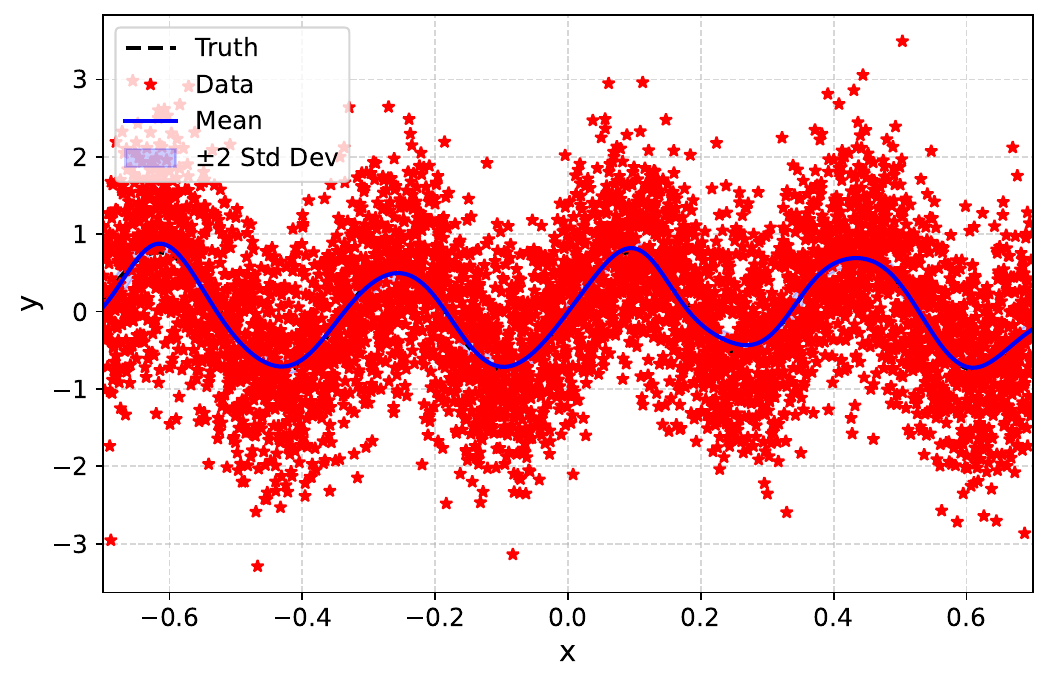}
            \end{overpic}
            \caption{Nonlinear equation: noise scale $\sigma_{f} = 0.8$ for $f(x)$ and $\sigma_{b} = 0.05$ for boundary data. Comparison of prediction values of $u(x)$ (first row) and $f(x)$ (second row) obtained by different methods.}
            \label{nonlinear_large_noise}
        \end{figure}

To demonstrate the ability of our proposed method to deal with big data, we assume that the parameter $k$ is known and add large noise to the residual measurements. Specifically, we generate $N_{u} = 3$ and $N_{b} = 2$ interior and boundary measurements corrupted by white noise with zero means and noise scale $\sigma_{u} = \sigma_{b} = 0.05$. Moreover, to counteract the large noise added to the residual measurements, we generate $N_{f} = 5000$ residual points corrupted with noise scale $\sigma_{f} = 0.8$, which is almost 
three times of the largest value of $f(x)$. In Fig. \ref{nonlinear_large_noise} we plot the predicted values of $u(x)$ and $f(x)$ in the first row and second row, respectively. Both DTEKI and SDTEKI can obtain accurate predictions with reliable uncertainty quantification compared to HMC. Moreover, choosing parameters in the major varying directions in SDTEKI can further reduce the uncertainty estimate and provide a more accurate estimate. In terms of computational time, SDTEKI is nearly 50 times faster than HMC, greatly improving the efficiency with comparable accuracy as demonstrated in Table \ref{tab:4}
\begin{table}[htbp]
        \centering
        \begin{tabular}{ccccc}
\hline Method & $e_{u}$   & Walltime & Network size \\
\hline  
  HMC & $6.43\%$  & 204.68s &960  \\
DTEKI  & $7.01\%$  & 12.71s &960 \\
SDTEKI  & $6.56\%$  & 4.53s & 320\\
\hline
\end{tabular}
        \caption{Nonlinear equation: $L_{2}$ error of predictions with respect to different methods. The running time for different methods is also presented.}
        \label{tab:5}
    \end{table}  
\subsection{2D Darcy flow}
\begin{figure}[htbp]
    \centering
    \begin{overpic}[width = 0.35\textwidth]{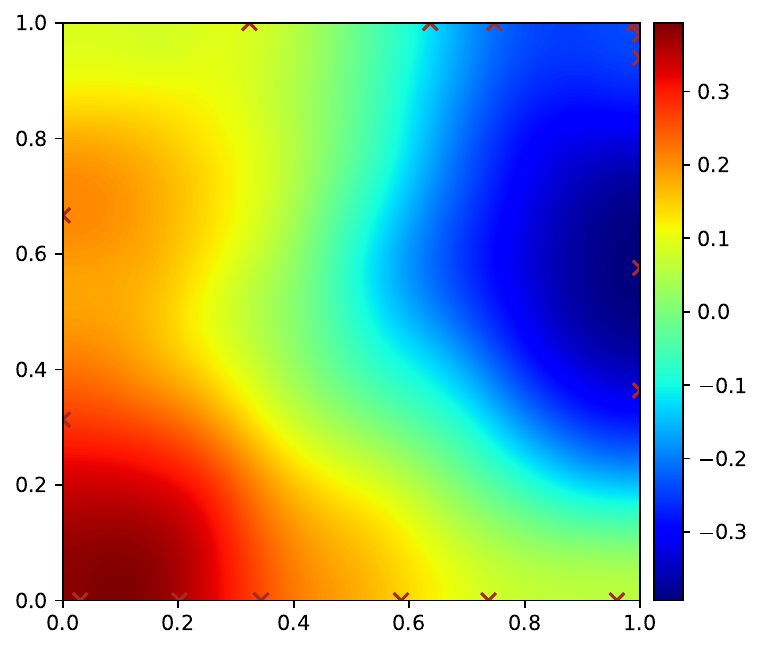}
    \end{overpic}
    \begin{overpic}[width = 0.35\textwidth]{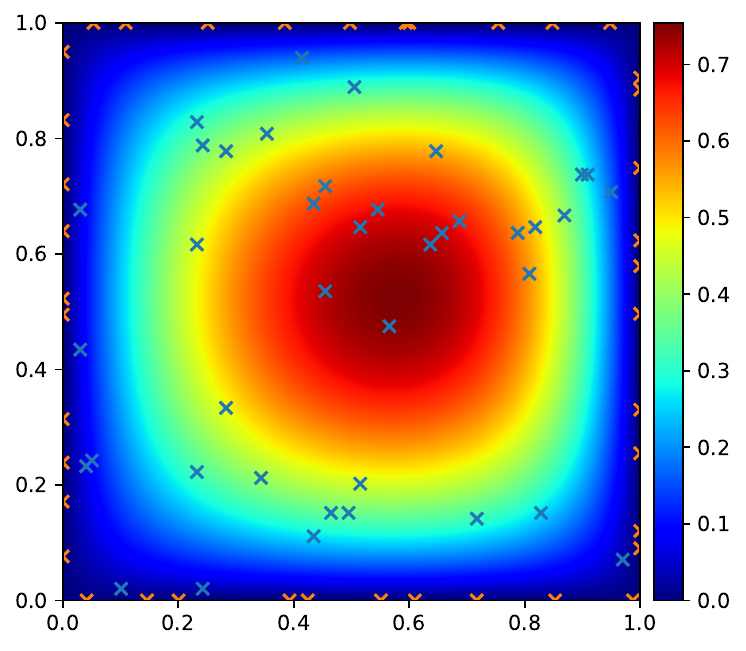}
    \end{overpic}
    \caption{Darcy equation: the true permeability field with boundary measurements (left), and the true solution field with interior and boundary measurements (right).}
    \label{true_darcy}
\end{figure}

Consider a 2D high-dimensional Darcy flow problem \cite{huang2022iterated,gao2024adaptive}, defined as 
\begin{equation}
\begin{split}
    -\nabla\cdot (\exp(\kappa(\mathbf{x},\mathbf{\lambda}))\nabla u(\mathbf{x})) &= f(\mathbf{x}),\quad  \mathbf{x}\in \Omega,\\ 
    u(\mathbf{x}) &= 0, \quad\quad  \mathbf{x}\in \partial \Omega.
\end{split}
\end{equation}
The right-hand side is $f(x)=10$ for simplicity. Here $\Omega = [0,1]^{2}$. The goal of this problem is to infer the true permeability field $\kappa(\mathbf{x}, \mathbf{\lambda})$ and also the solution field $u(\mathbf{x})$ from noisy measurements. In this experiment, the prior of permeability field $\kappa(\mathbf{x}, \lambda)$ is assumed to follow a centered Gaussian with covariance
\begin{equation*}
    C = (-\Delta + \tau^{2})^{-d},
\end{equation*}
where $-\Delta$ is the Laplacian operator on the space of spatial zero mean functions with homogeneous Neumann boundary conditions on $\Omega$;  $\tau>0$ is the inverse length scale of the random field and $d$ controls the regularity of the random field. In this study, we set $\tau = 5, d = 4$ to generate smooth priors. To sample from the random field, we consider the following Karhunen-Lo\`eve (KL) expansion:
\begin{equation}
    \log \kappa(\mathbf{x},\mathbf{\lambda})=\sum_{l\in Z}\lambda_{(l)}\sqrt{\mu_l}\psi_l(\mathbf{x}),
\end{equation}
where $Z = \mathbb{Z}^{0+}\times \mathbb{Z}^{0+}\backslash \{0,0\}$, $\lambda_{l}\sim \mathcal{N}(0,1)$ i.i.d. and $\mu_{l}$ and $\psi_{l}$ are the eigenpairs of the expansion. In practice, we truncate it with the first $256$ terms to generate the ground truth and the true solution field as demonstrated in Fig. \ref{true_darcy}.  

In this experiment, the permeability field is an infinite-dimensional Gaussian stochastic process. Thus, we will use two separate networks with the same architecture to approximate $\kappa(\mathbf{x}, \lambda)$ and $u(\mathbf{x})$, respectively. The networks have two hidden layers with 10 latent neurons, where the polynomial degree is at most 7. The regularization parameter $\alpha$ is set to $5$ to regularize the network with a much larger network size.

To recover the true permeability field, we will generate $N_{u} = 40$ interior measurements together with $N_{b} = 40$ boundary measurements corrupted by Gaussian noise white zero means and noise scales $\sigma_{u} = \sigma_{b} = 0.01$. We use $N_{f} = 5000$ residual measurements with noise scale $\sigma_{f} = 0.1$ to constrain the equation to provide physical information. Moreover, to ensure well-posedness, we will generate 16 boundary measurements of $\kappa$ with noise scale 0.01.

\begin{figure}[htbp]
            \centering 
            \begin{overpic}[width = 0.31\textwidth]{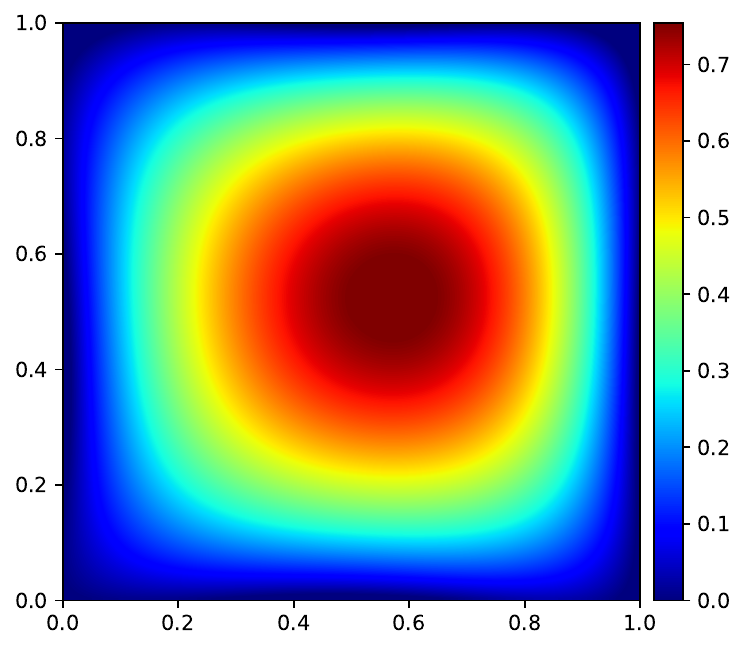}
            \put (33,87) {\textbf{\footnotesize HMC}}
            \end{overpic}
            \begin{overpic}[width = 0.31\textwidth]{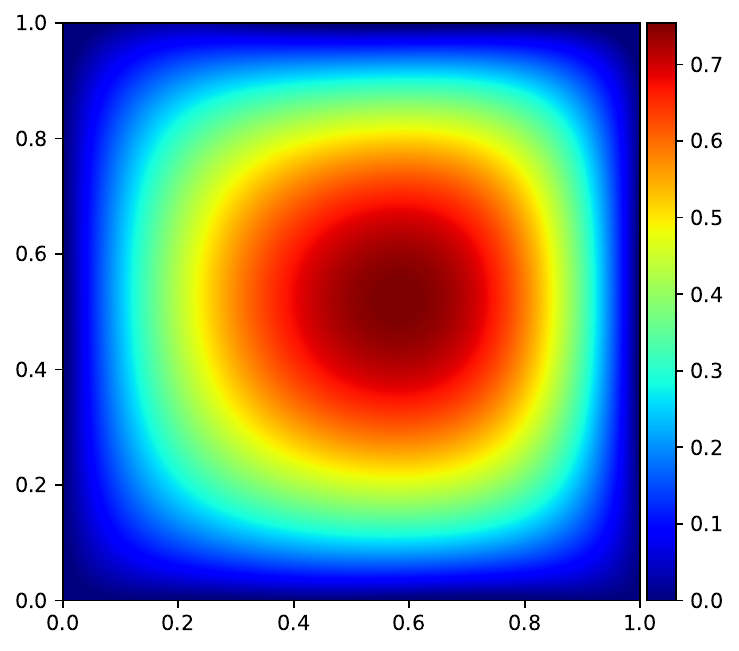}
            \put (33,87) {\textbf{\footnotesize DTEKI}}
            \end{overpic}
            \begin{overpic}[width = 0.31\textwidth]{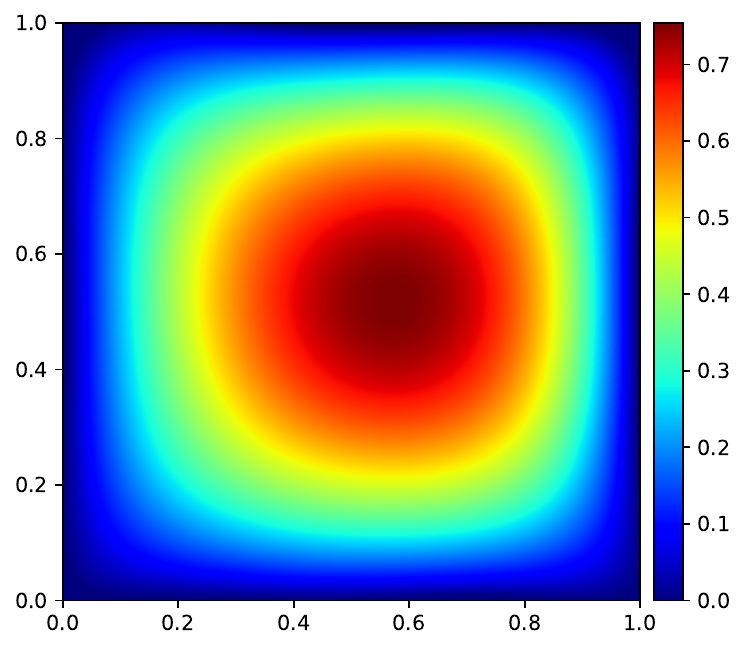}
            \put (33,87) {\textbf{\footnotesize SDTEKI}}
            \end{overpic}
            \begin{overpic}[width = 0.31\textwidth]{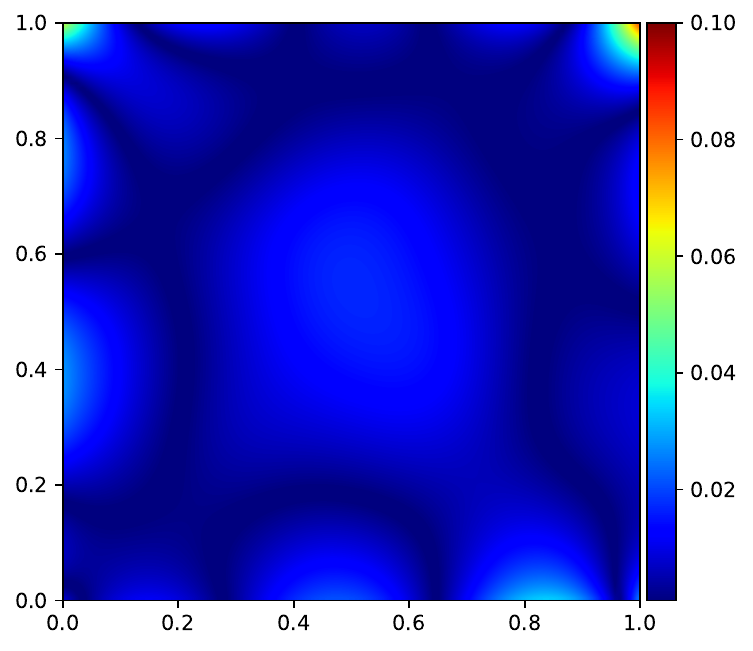}
            \end{overpic}
            \begin{overpic}[width = 0.31\textwidth]{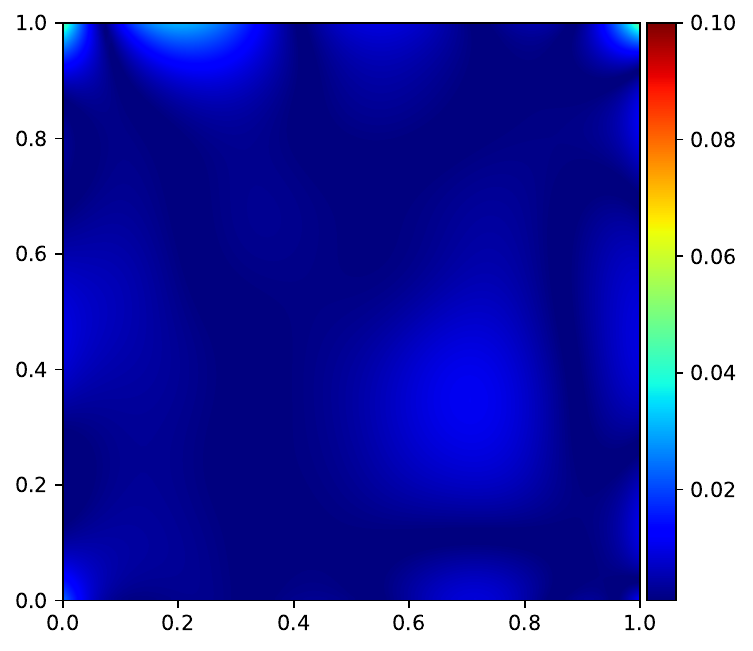}
            \end{overpic}
            \begin{overpic}[width = 0.31\textwidth]{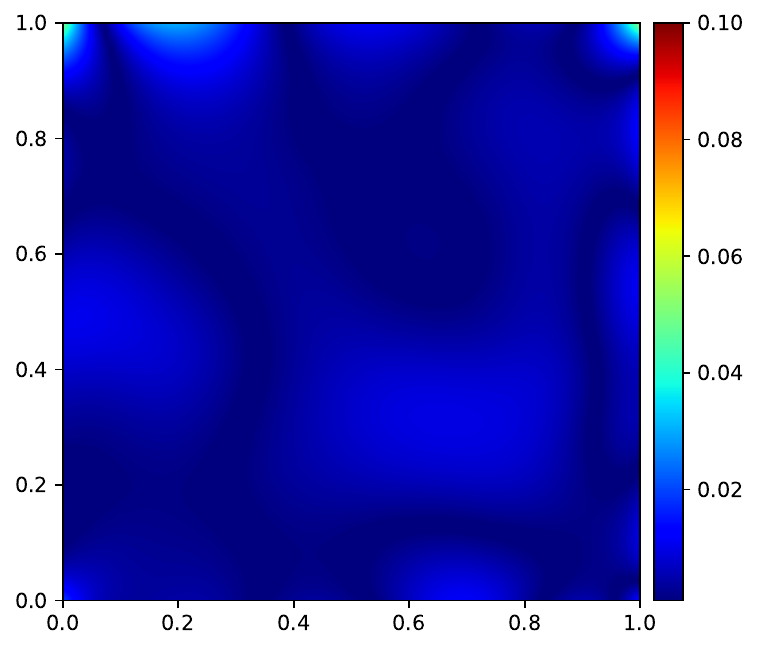}
            \end{overpic}
            \begin{overpic}[width = 0.31\textwidth]{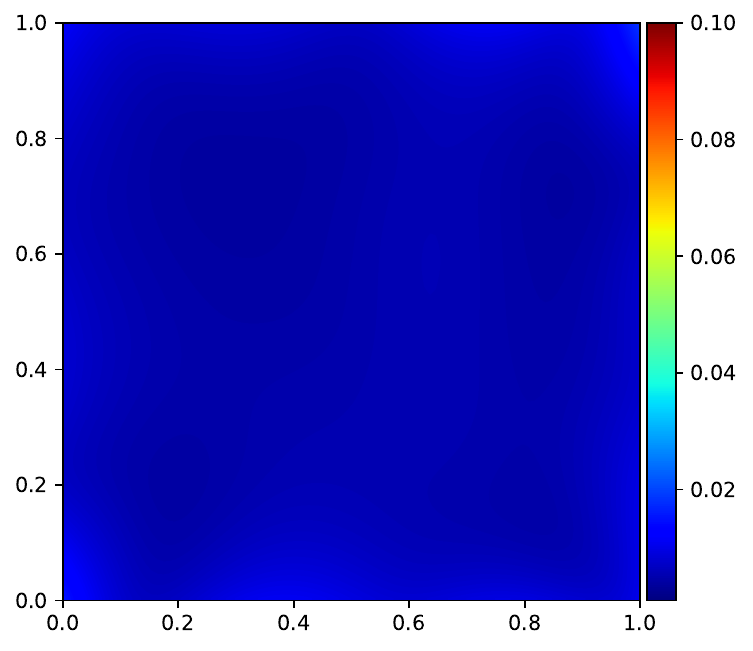}
            \end{overpic}
            \begin{overpic}[width = 0.31\textwidth]{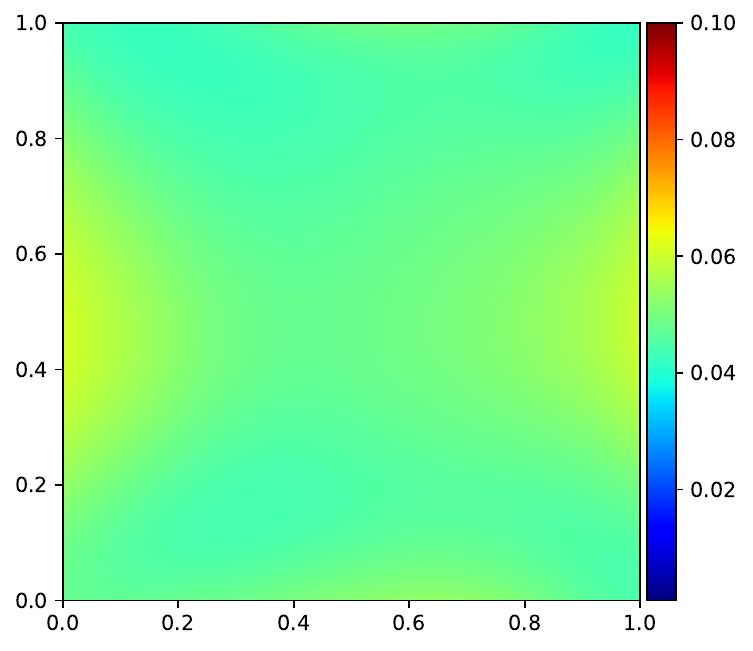}
            \end{overpic}
            \begin{overpic}[width = 0.31\textwidth]{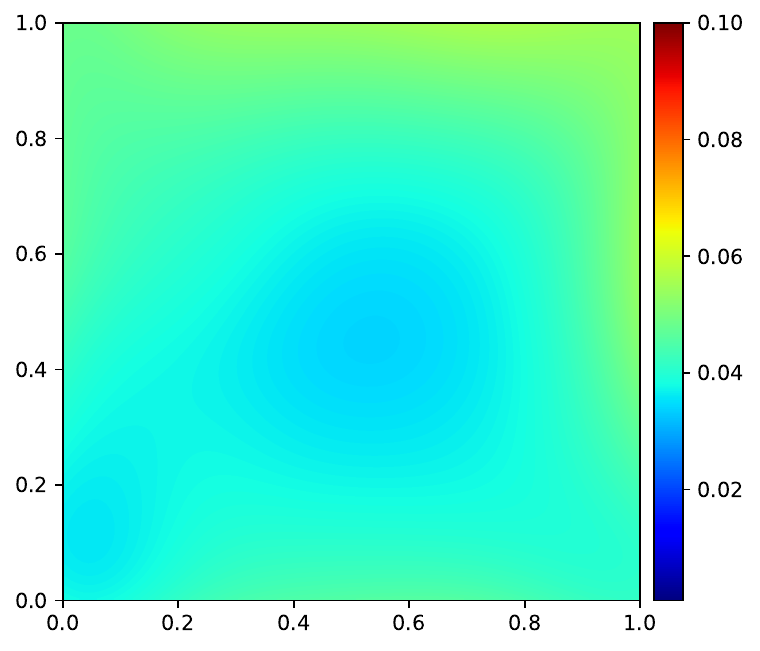}
            \end{overpic}
            \caption{Darcy equation: noise scales $\sigma_{f} = 0.1$ for $f(x)$, $\sigma_{b} = 0.01$ for boundary measurements and $\sigma_{u} = 0.01$ for the interior measurements. Comparison of mean prediction values of $u(x)$ (first row), absolute errors (second row), and 2 times standard deviations (third row) obtained by different methods.}
            \label{darcy_u}
        \end{figure}

In Fig. \ref{darcy_u}, we show the recovered mean solution field (first row) and their corresponding standard deviations (second row) obtained using HMC, DTEKI, and SDTEKI, respectively. Both DTEKI and SDTEKI demonstrate superior performance compared to HMC, achieving higher prediction accuracy with smaller $L_2$ errors. For instance, as shown in Table \ref{tab:6}, DTEKI and SDTEKI achieve solution errors ($e_u$) of $1.28\%$ and $1.33\%$, respectively, outperforming HMC's $2.24\%$. Regarding the posterior estimation of $\kappa$, DTEKI and SDTEKI both reduce the relative errors to approximately $13.05\%$ and $15.04\%$, compared to HMC's $22.84\%$. Moreover, it is clear that the standard deviations of DTEKI and SDTEKI can cover the absolute error both for the solution field and permeability field, while HMC is not capable of this, possibly due to the dominance of the likelihood. This indicates that our methods not only improve the accuracy of the solution field but also provide more reliable uncertainty quantification.

One of the key reasons for the improved performance of DTEKI and SDTEKI lies in the mini-batch training, which effectively mitigates overfitting, a known issue for HMC, especially in high-dimensional parameter spaces. Moreover, the ensemble-based estimation in the EKI framework further enhances the robustness of our methods.

The efficiency advantage of our methods is particularly noteworthy. As detailed in Table~\ref{tab:6}, DTEKI and SDTEKI achieve significant computational speed-ups, with SDTEKI completing the inference process in just $18.76$ seconds, compared to HMC's $576.83$ seconds—a nearly 30-fold improvement. This efficiency stems from two factors: the use of mini-batch updates in DTEKI and the dimensionality reduction achieved through the subspace method in SDTEKI. The latter also allows SDTEKI to use a smaller network size (692 parameters) while maintaining a comparable or even superior accuracy to HMC, which requires 2080 parameters.

In summary, DTEKI and SDTEKI not only deliver higher prediction accuracy than HMC but also exhibit significantly better efficiency. The combination of mini-batch training and dimensionality reduction enables our methods to handle complex, large-scale inverse problems effectively, saving computational costs and achieving reliable predictions simultaneously.

\begin{figure}[htbp]
            \centering 
            \begin{overpic}[width = 0.31\textwidth]{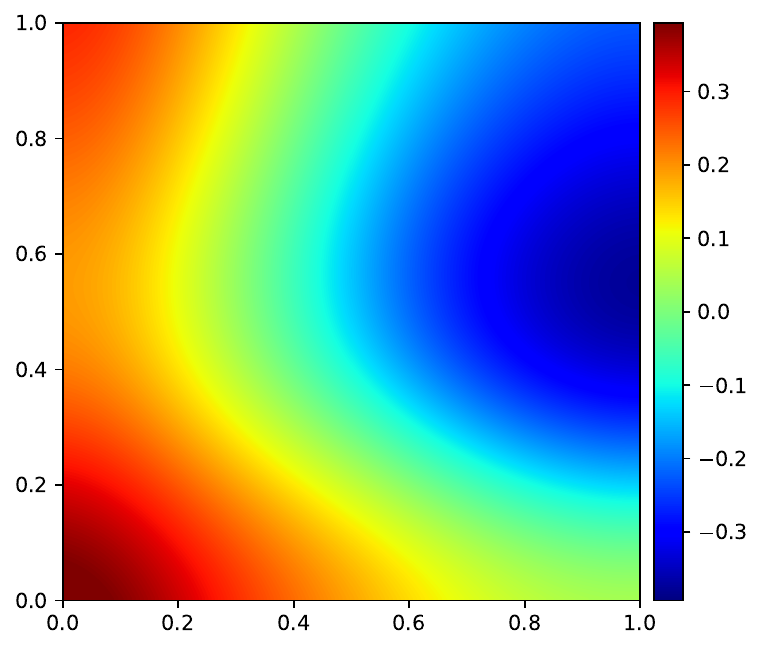}
            \put (33,85) {\textbf{\footnotesize HMC}}
            \end{overpic}
            \begin{overpic}[width = 0.31\textwidth]{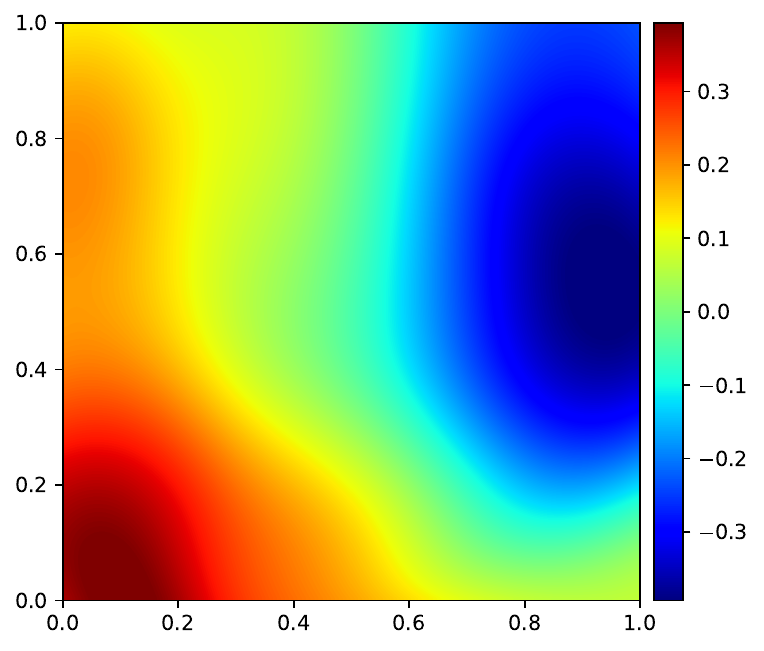}
            \put (33,85) {\textbf{\footnotesize DTEKI}}
            \end{overpic}
            \begin{overpic}[width = 0.31\textwidth]{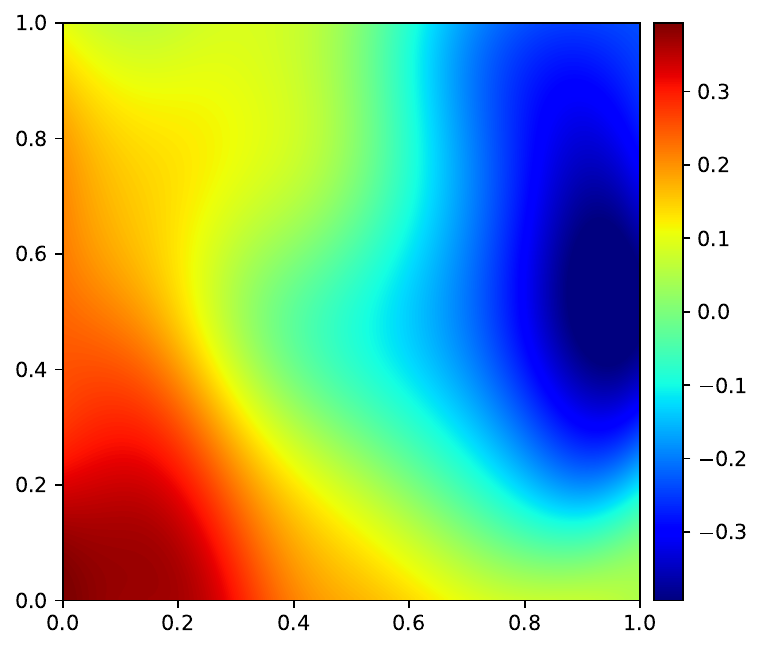}
            \put (33,85) {\textbf{\footnotesize SDTEKI}}
            \end{overpic}
            \begin{overpic}[width = 0.31\textwidth]{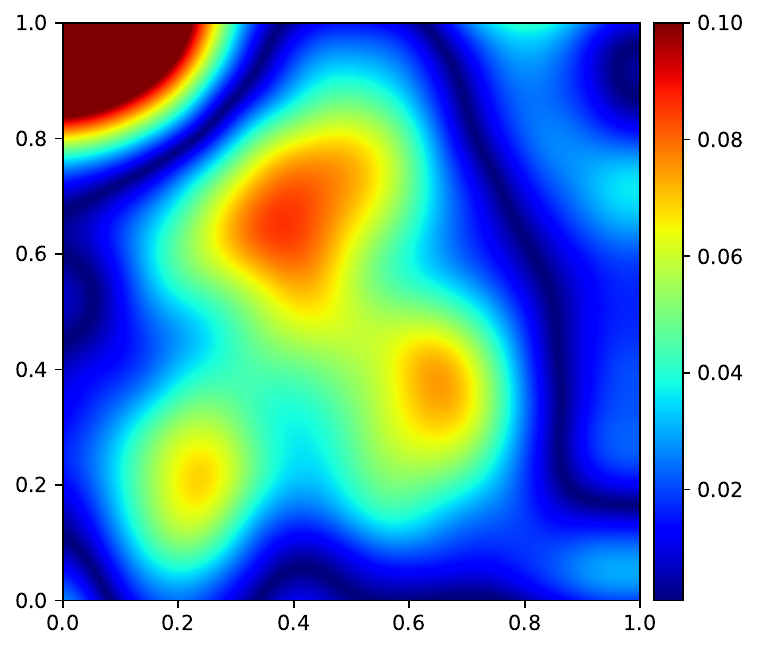}
            \end{overpic}
            \begin{overpic}[width = 0.31\textwidth]{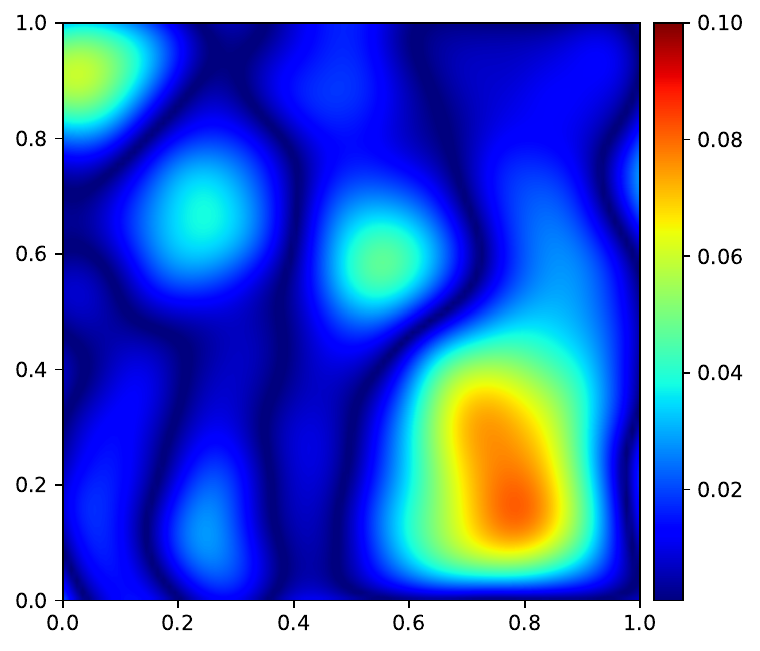}
            \end{overpic}
            \begin{overpic}[width = 0.31\textwidth]{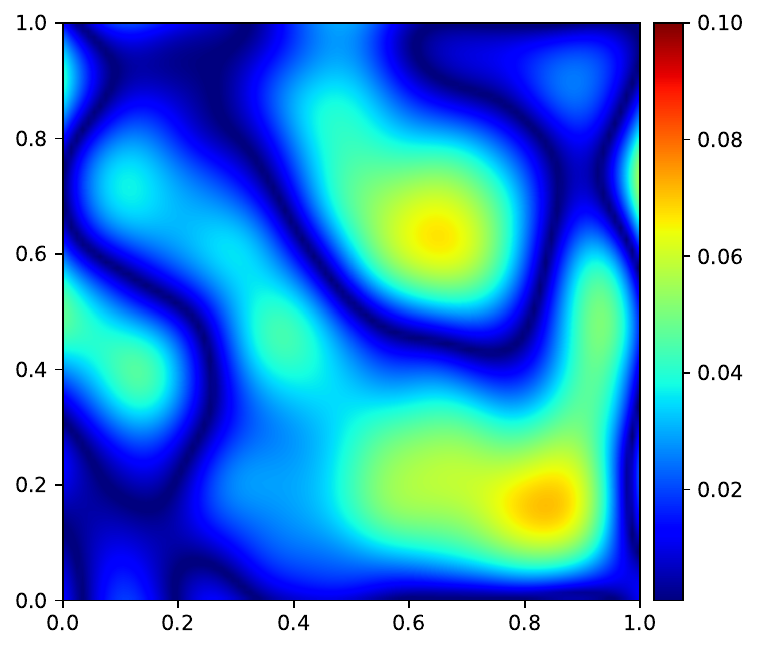}
            \end{overpic}
            \begin{overpic}[width = 0.31\textwidth]{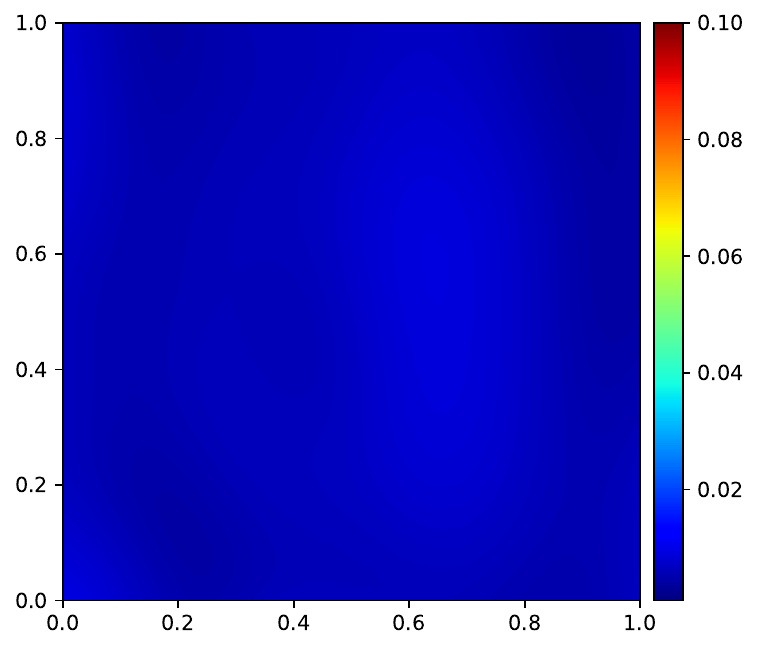}
            \end{overpic}
            \begin{overpic}[width = 0.31\textwidth]{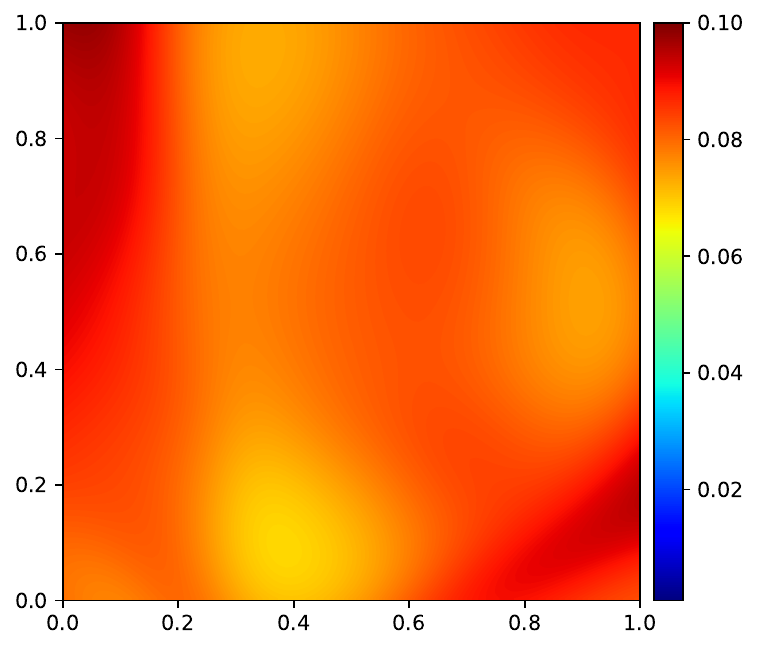}
            \end{overpic}
            \begin{overpic}[width = 0.31\textwidth]{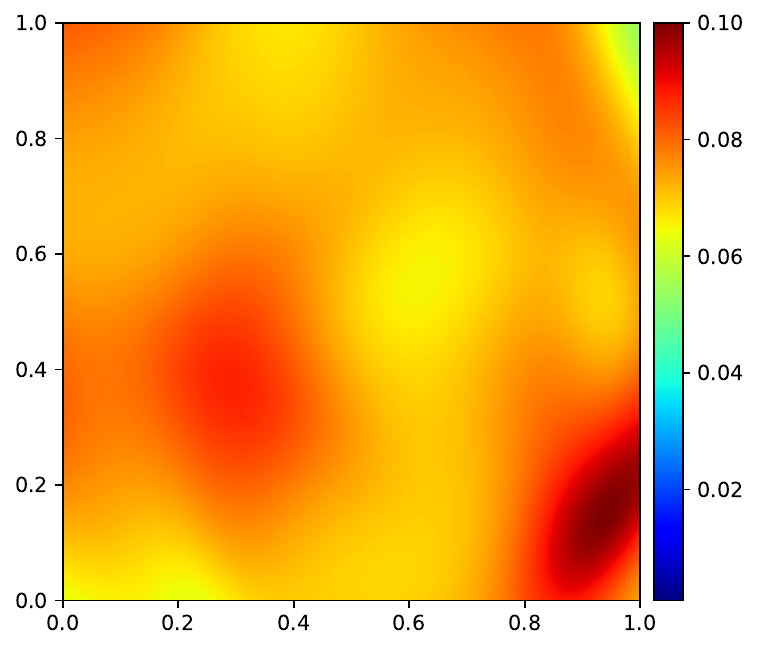}
            \end{overpic}
            \caption{Darcy equation: noise scales $\sigma_{f} = 0.1$ for $f(x)$, $\sigma_{b} = 0.01$ for boundary measurements and $\sigma_{u} = 0.01$ for the interior measurements. Comparison of mean prediction values of $u(x)$ (first row), absolute errors (second row), and 2 times standard deviations (third row) obtained by different methods.}
            \label{darcy_kappa}
        \end{figure}

\begin{table}[t]
        \centering
        \begin{tabular}{ccccc}
\hline Method & $e_{u}$ & $e_{\kappa}$  & Walltime & Network size \\
\hline  
  HMC & $2.24\%$ & $22.84\%$ & 576.83s &2080  \\
DTEKI  & $1.28\%$ &$13.05\%$  & 28.80s &2080 \\
SDTEKI  & $1.33\%$ & $15.04\%$ & 18.76s & 692\\
\hline
\end{tabular}
        \caption{Darcy equation: $L_{2}$ error of predictions with respect to different methods. The running time (GPU RTX 6000) and network sizes for different methods are also presented.}
        \label{tab:6}
    \end{table} 
    
\section{Summary}
In this paper, we proposed utilizing cKANs to replace MLPs in the process of training for B-PINNs, resulting in Bayesian PIKANs, which could reduce the number of parameters and thus accelerate training. To avoid the computation of back-propagation, we adopted the ensemble Kalman inversion method together with its Tikhonov version to achieve gradient-free inference. Moreover, to overcome the rank deficiency of the covariance matrix, we applied the dropout constraint to the ensemble. Furthermore, to reduce the parameter dimension and reduce overfitting, we employed the active subspace method to construct a low-dimensional parameter space during the iteration. This could greatly reduce the number of parameters to accelerate training while preserving the prediction accuracy. Our experiments demonstrated that both DTEKI and SDTEKI can achieve comparable accuracy compared to HMC, while SDTEKI is the fastest and can also scale well to large-scale problems. 
For future work, we consider adding momentum to the dynamics of DTEKI, which could further accelerate convergence rate. Moreover, functional priors can be trained to provide more useful information during the inference process.

\section*{Acknowledgments}
The first author would like to thank Professor Xuhui Meng for helpful discussions and comments regarding the applications of the proposed method in different problems.
This work was supported by the MURI grant (FA9550-20-1-0358), the ONR Vannevar Bush Faculty Fellowship (N0001422-1-2795), and the U.S. Department of Energy, Advanced Scientific Computing Research program, under the Scalable, Efficient and Accelerated Causal Reasoning Operators, Graphs and Spikes for Earth and Embedded Systems (SEA-CROGS) project, (DE-SC0023191).
Additional funding was provided by GPU Cluster for Neural PDEs and Neural Operators to support MURI Research and Beyond, under Award \#FA9550-23-1-0671.

\appendix
\counterwithin{figure}{section}
\counterwithin{table}{section}

\section{Ablation study}
To compare the performance when using MLP and cKANs, we run DTEKI separately with nearly the same number of parameters. Other settings are kept the same to make sure equality. The network sizes used for MLPs are listed in Table \ref{tab:network_size} for different problems. For simplicity, we only test the inverse problems when the physical parameter is unknown. The results show that for all problems the running time is similar. 
\begin{table}[htbp]
\centering
\begin{tabular}{@{}ccc@{}}
\toprule
\textbf{Problem}     & \textbf{Network Structure}&
\textbf{Network Size}\\ \midrule
Transport            & \(2 \times 30 \times 30 \times 1\) &1051\\
Diffusion  &  \(1 \times 29 \times 29 \times 1\)& 958 \\
Nonlinear  &  \(1 \times 29 \times 29 \times 1\)& 958 \\
Darcy      & \(2 \times 30 \times 30 \times 1\) & 2102                 \\ \bottomrule
\end{tabular}
\caption{Network structure for different problems with MLP and the corresponding total number of parameters.}
\label{tab:network_size}
\end{table}

The performance of DTEKI with cKANs is much better than MLPs. For the transport equation, Fig. \ref{fig:transport_mlp} demonstrates the predicted posterior distribution of $k$ compared to the true posterior by different methods. The predicted mean of $k$ using MLP-DTEKI is 1.032 and the standard deviation is $0.0190$. Note that using MLP, the prediction posterior can still be close to the true posterior, while the accuracy of the prediction is less than cKANs, indicating that cKANs can achieve higher accuracy than MLPs.
\setcounter{figure}{0}  
\begin{figure}[htbp]
    \centering
\includegraphics[width=0.5\linewidth]{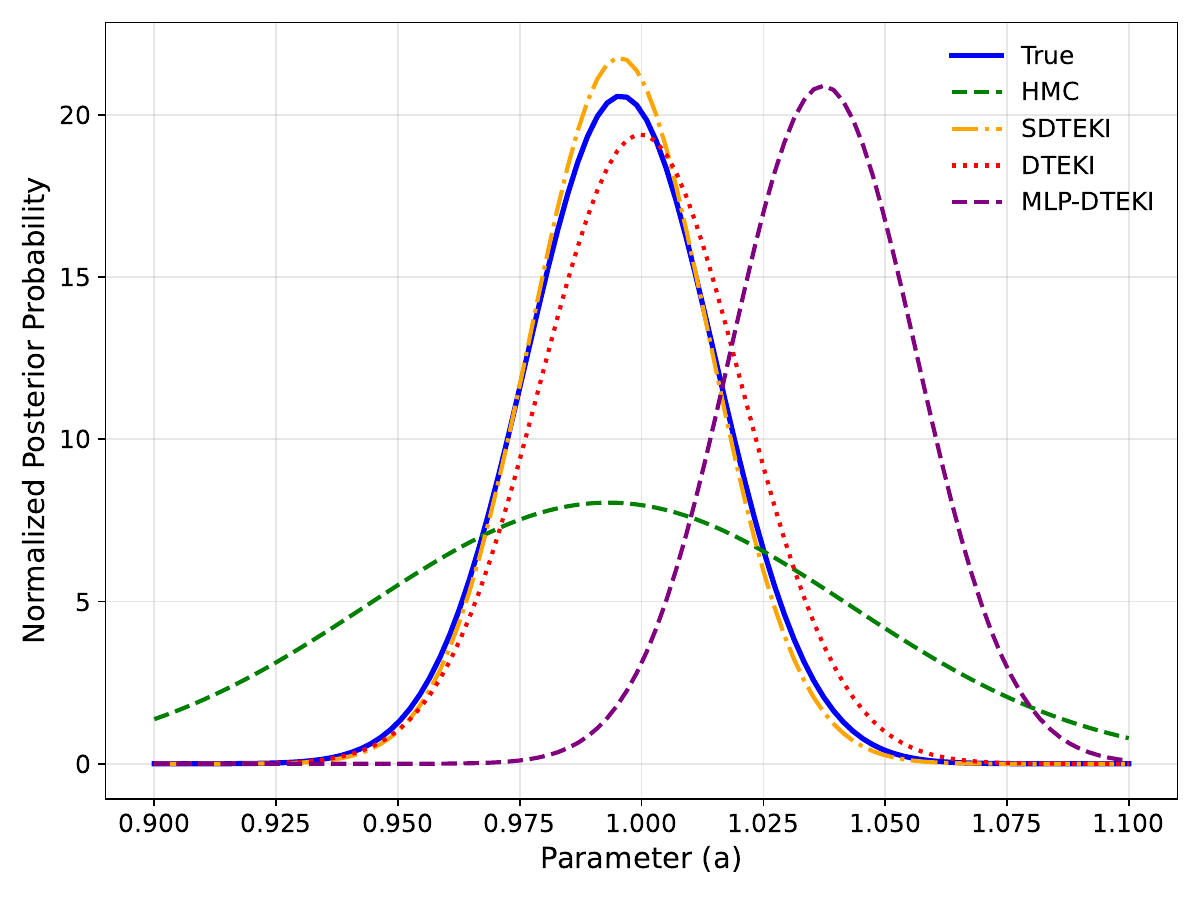}
    \caption{Transport equation: comparison of predicted posteriors and the true posterior distribution obtained by different methods.}
    \label{fig:transport_mlp}
\end{figure}

\begin{figure}[htbp]
\centering 
\begin{overpic}[width = 0.32\textwidth]{figures/diffusion/HMC_kan_u_inverse_ensemble.pdf}
\put (-5, 25) {\rotatebox{90}{\footnotesize HMC}}
\end{overpic}
\begin{overpic}[width = 0.32\textwidth]{figures/diffusion/HMC_kan_f_inverse_ensemble.pdf}
\end{overpic}
\begin{overpic}[width = 0.27\textwidth]{figures/diffusion/HMC_kan_k_inverse_ensemble.pdf}
\end{overpic}
\begin{overpic}[width = 0.32\textwidth]{figures/diffusion/EKI_kan_u_inverse_ensemble.pdf}
\put (-5, 25) {\rotatebox{90}{\footnotesize DTEKI}}
\end{overpic}
\begin{overpic}[width = 0.32\textwidth]{figures/diffusion/EKI_kan_f_inverse_ensemble.pdf}
\end{overpic}
\begin{overpic}[width = 0.27\textwidth]{figures/diffusion/EKI_kan_k_inverse_ensemble.pdf}
\end{overpic}
\begin{overpic}[width = 0.32\textwidth]{figures/diffusion/EKI_kan_u_inverse_subspace_ensemble.pdf}
\put (-5, 25) {\rotatebox{90}{\footnotesize SDTEKI}}
\end{overpic}
\begin{overpic}[width = 0.32\textwidth]{figures/diffusion/EKI_kan_f_inverse_subspace_ensemble.pdf}
\end{overpic}
\begin{overpic}[width = 0.27\textwidth]{figures/diffusion/EKI_kan_k_inverse_subspace_ensemble.pdf}
\end{overpic}
\begin{overpic}[width = 0.32\textwidth]{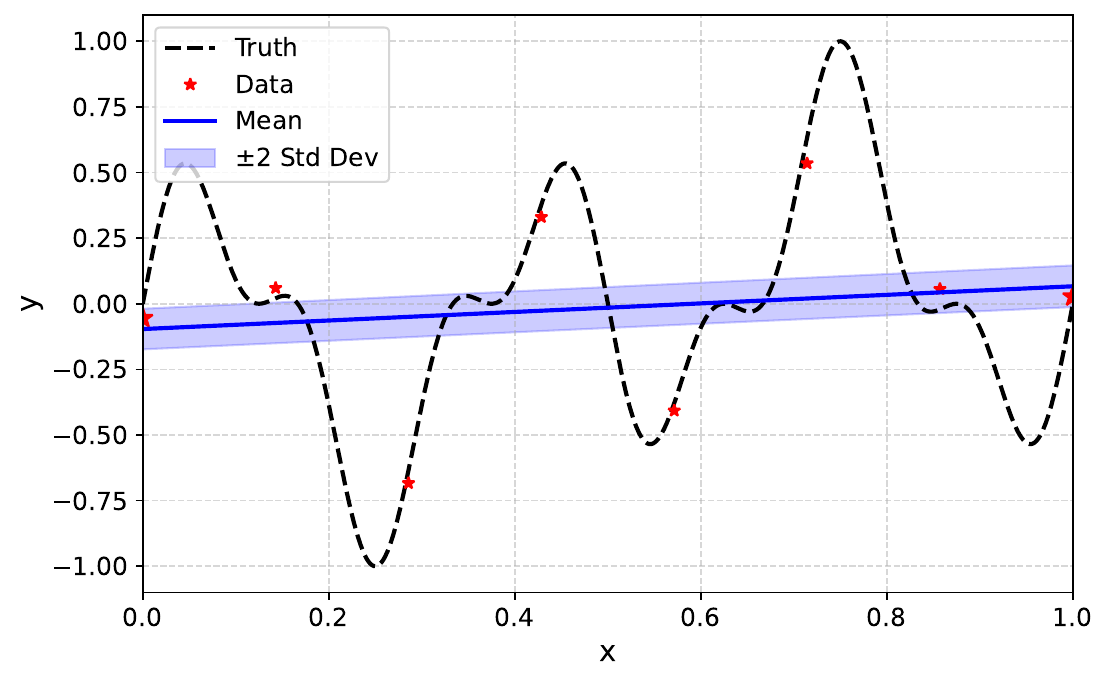}
\put (-5, 15) {\rotatebox{90}{\footnotesize MLP-DTEKI}}
\end{overpic}
\begin{overpic}[width = 0.32\textwidth]{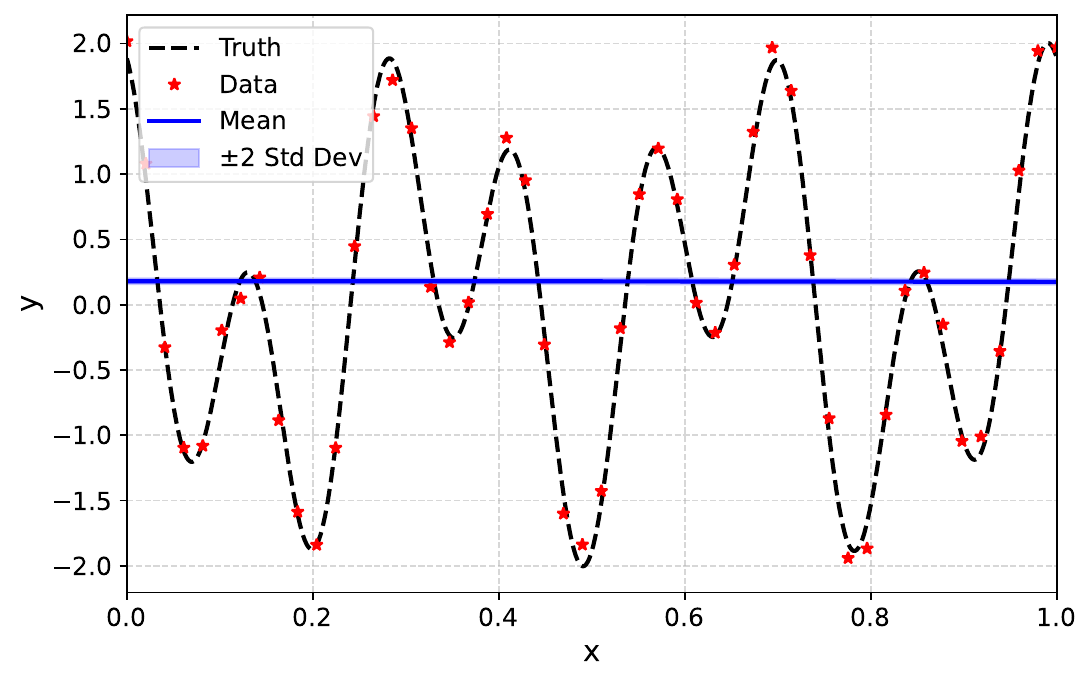}
\end{overpic}
\begin{overpic}[width = 0.27\textwidth]{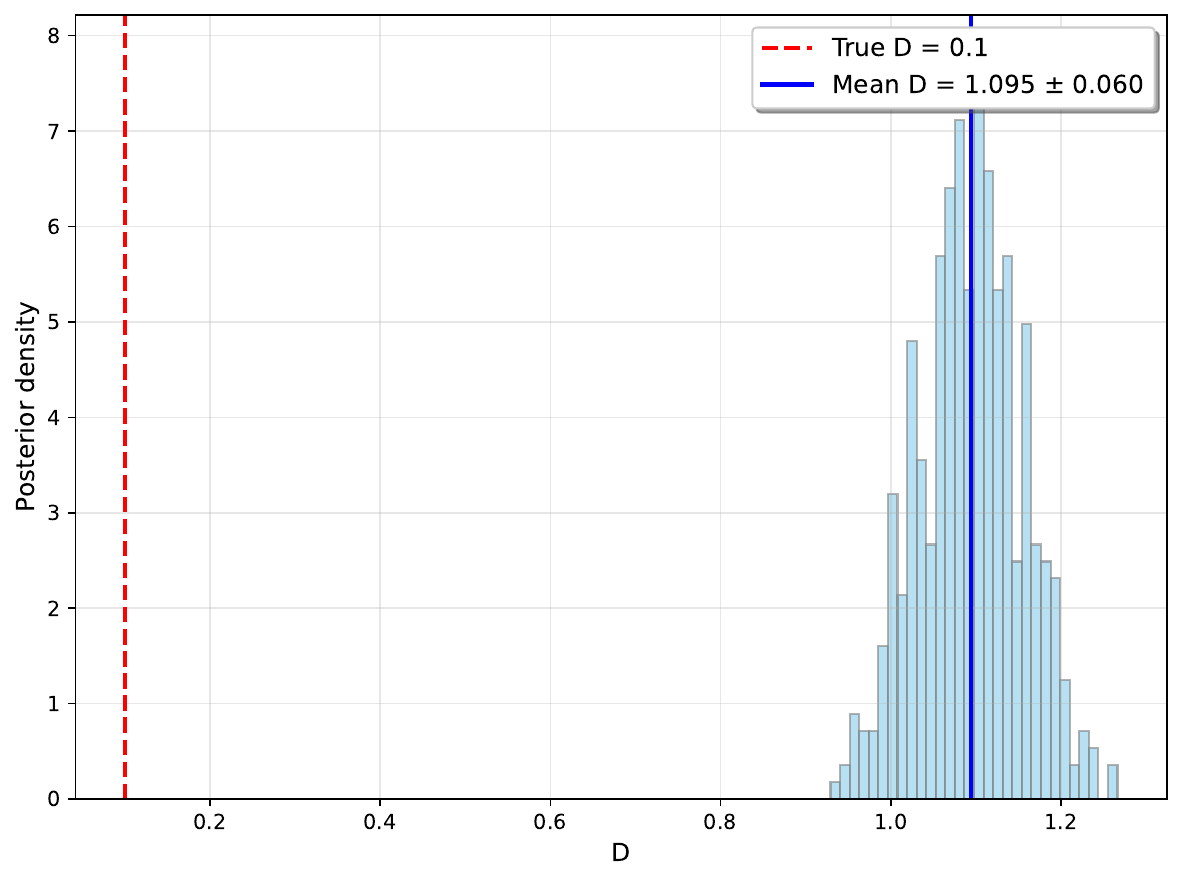}
\end{overpic}
\caption{Diffusion equation: predicted mean values of $u(x)$ and $f(x)$ and predicted posterior distribution from left to right using MLP-DTEKI.}
\label{diffusion_mlp}
\end{figure}

\begin{figure}[htbp]
\centering 
\begin{overpic}[width = 0.32\textwidth]{figures/nonlinear/HMC_kan_u_inverse_ensemble.pdf}
\put (-5, 25) {\rotatebox{90}{\footnotesize HMC}}
\end{overpic}
\begin{overpic}[width = 0.32\textwidth]{figures/nonlinear/HMC_kan_f_inverse_ensemble.pdf}
\end{overpic}
\begin{overpic}[width = 0.27\textwidth]{figures/nonlinear/HMC_kan_k_inverse_ensemble.pdf}
\end{overpic}
\begin{overpic}[width = 0.32\textwidth]{figures/nonlinear/EKI_kan_u_inverse_ensemble.pdf}
\put (-5, 25) {\rotatebox{90}{\footnotesize DTEKI}}
\end{overpic}
\begin{overpic}[width = 0.32\textwidth]{figures/nonlinear/EKI_kan_f_inverse_ensemble.pdf}
\end{overpic}
\begin{overpic}[width = 0.27\textwidth]{figures/nonlinear/EKI_kan_k_inverse_ensemble.pdf}
\end{overpic}
\begin{overpic}[width = 0.32\textwidth]{figures/nonlinear/EKI_kan_u_inverse_subspace_ensemble.pdf}
\put (-5, 25) {\rotatebox{90}{\footnotesize SDTEKI}}
\end{overpic}
\begin{overpic}[width = 0.32\textwidth]{figures/nonlinear/EKI_kan_f_inverse_subspace_ensemble.pdf}
\end{overpic}
\begin{overpic}[width = 0.27\textwidth]{figures/nonlinear/EKI_kan_k_inverse_subspace_ensemble.pdf}
\end{overpic}
\begin{overpic}[width = 0.32\textwidth]{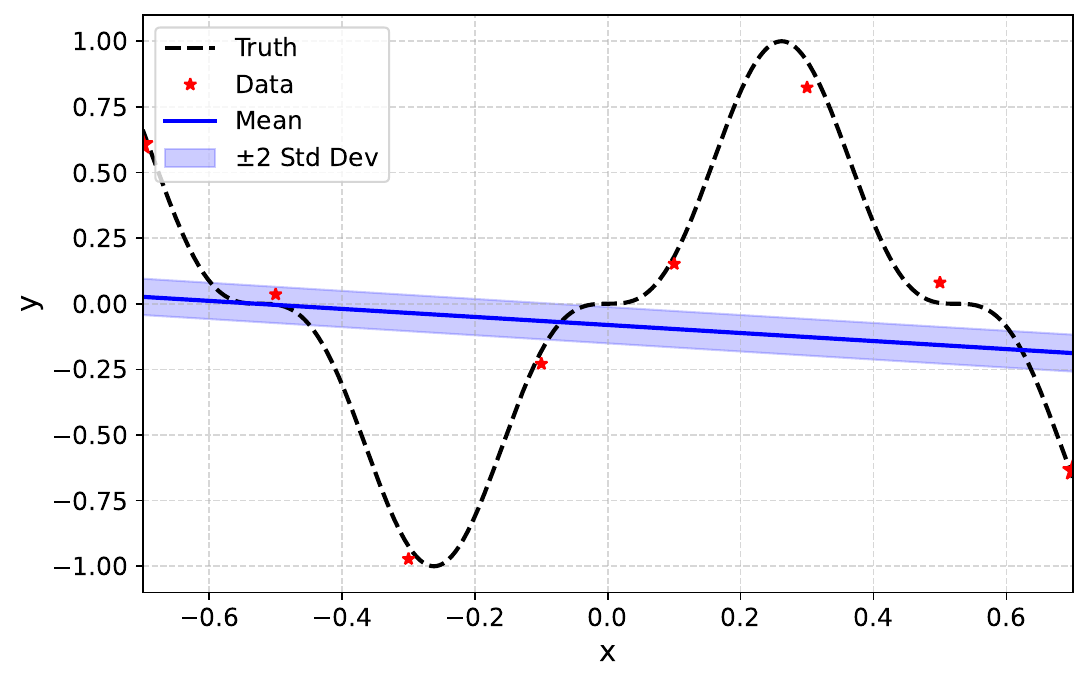}
\put (-5, 15) {\rotatebox{90}{\footnotesize MLP-DTEKI}}
\end{overpic}
\begin{overpic}[width = 0.32\textwidth]{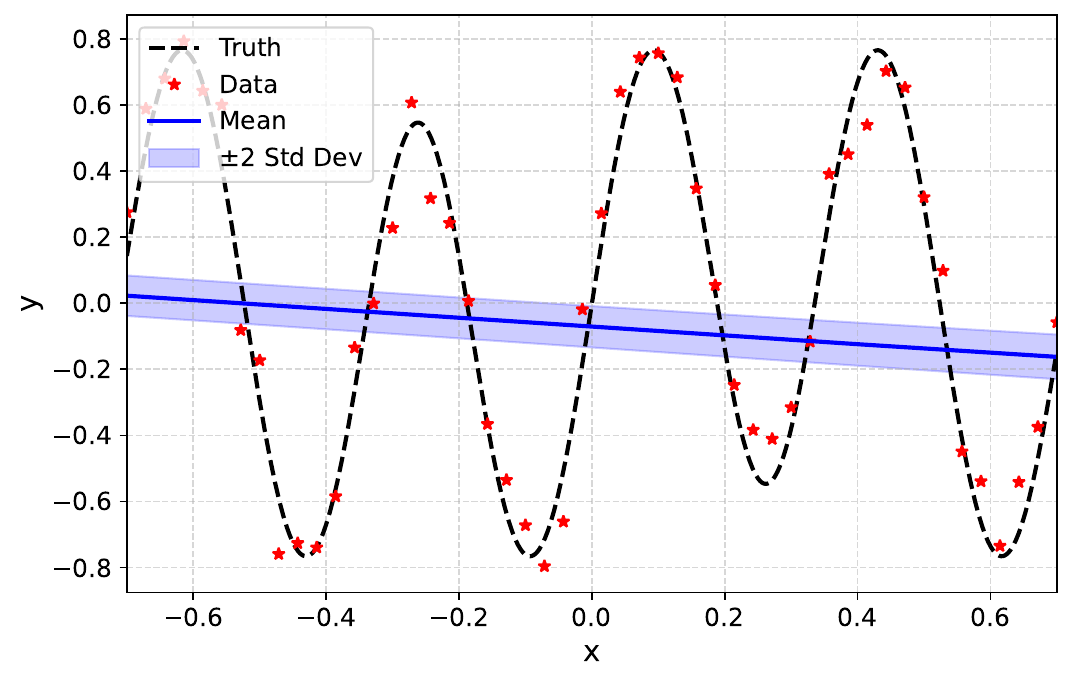}
\end{overpic}
\begin{overpic}[width = 0.27\textwidth]{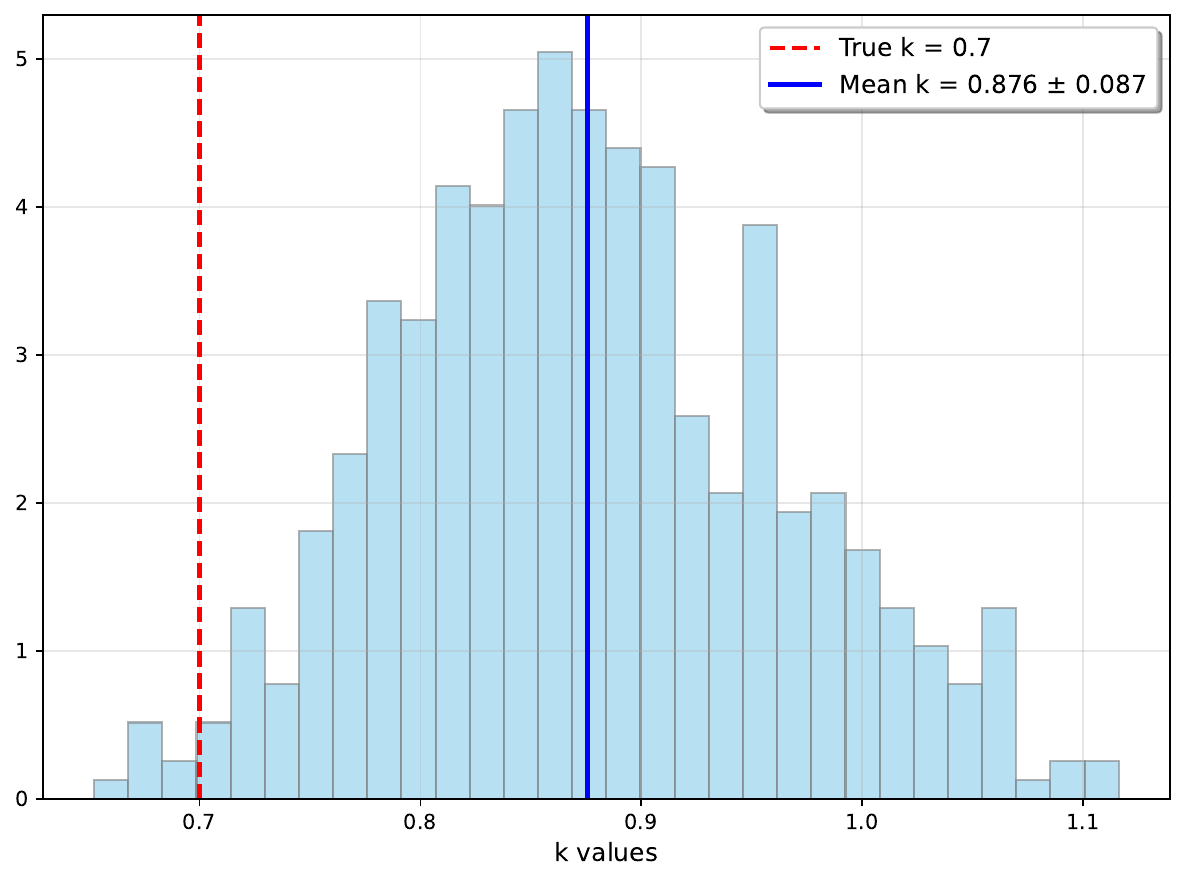}
\end{overpic}
\caption{Nonlinear equation: predicted mean values of $u(x)$ and $f(x)$ and predicted posterior distribution from left to right using MLP-DTEKI.}
\label{nonlinear_mlp}
\end{figure}

\begin{figure}[htbp]
\centering 
\begin{overpic}[width = 0.32\textwidth]{figures/darcy/u_HMC.pdf}
\put (-5, 40) {\rotatebox{90}{\footnotesize HMC}}
\end{overpic}
\begin{overpic}[width = 0.32\textwidth]{figures/darcy/u_error_HMC.pdf}
\end{overpic}
\begin{overpic}[width = 0.32\textwidth]{figures/darcy/u_std_HMC.pdf}
\end{overpic}
\begin{overpic}[width = 0.32\textwidth]{figures/darcy/u_EKI.pdf}
\put (-5, 40) {\rotatebox{90}{\footnotesize DTEKI}}
\end{overpic}
\begin{overpic}[width = 0.32\textwidth]{figures/darcy/u_error_EKI.pdf}
\end{overpic}
\begin{overpic}[width = 0.32\textwidth]{figures/darcy/u_std_EKI.pdf}
\end{overpic}
\begin{overpic}[width = 0.32\textwidth]{figures/darcy/u_subspace_EKI.pdf}
\put (-5, 40) {\rotatebox{90}{\footnotesize SDTEKI}}
\end{overpic}
\begin{overpic}[width = 0.32\textwidth]{figures/darcy/u_error_subspace_EKI.pdf}
\end{overpic}
\begin{overpic}[width = 0.32\textwidth]{figures/darcy/u_std_subspace_EKI.pdf}
\end{overpic}
\begin{overpic}[width = 0.32\textwidth]{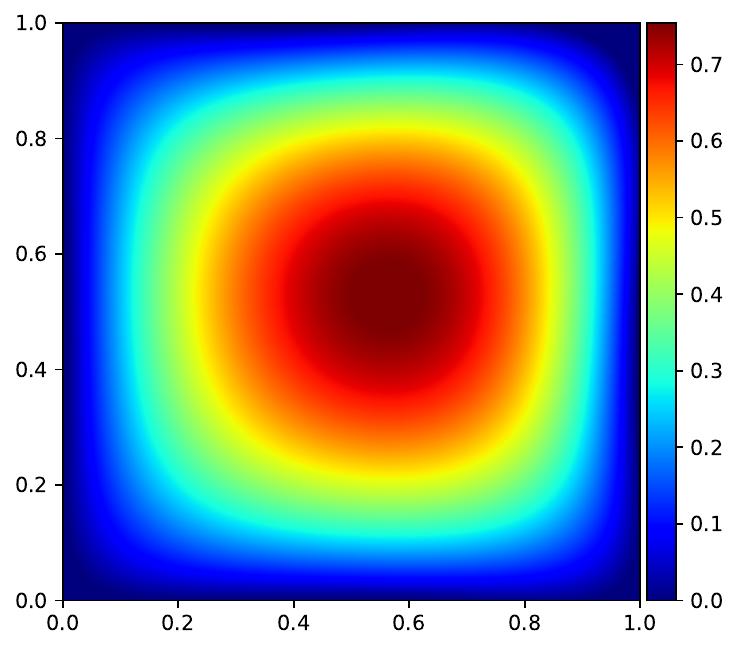}
\put (-5, 32) {\rotatebox{90}{\footnotesize MLP-DTEKI}}
\end{overpic}
\begin{overpic}[width = 0.32\textwidth]{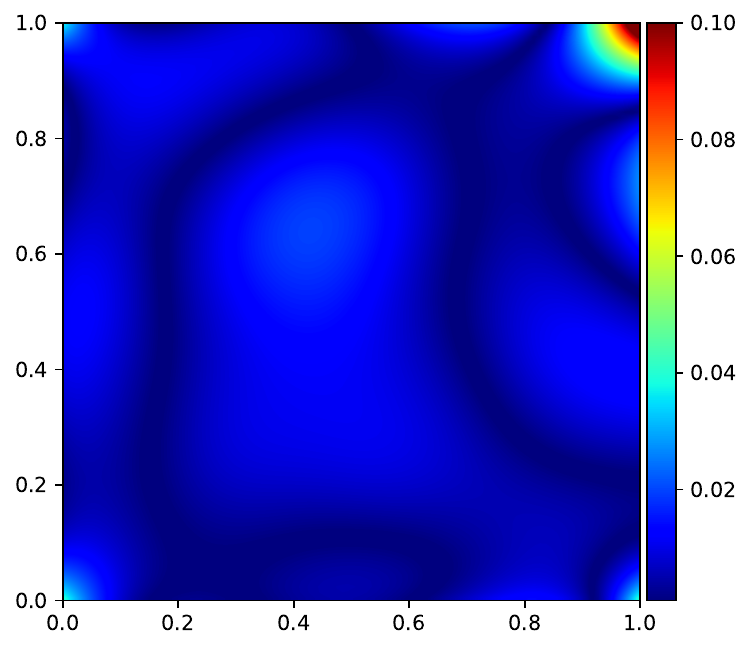}
\end{overpic}
\begin{overpic}[width = 0.32\textwidth]{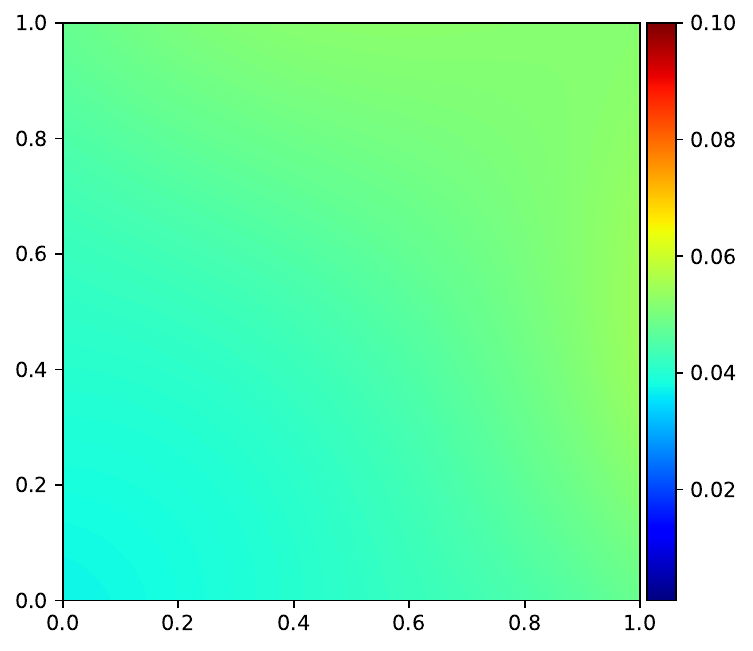}
\end{overpic}
\caption{Darcy equation: predicted mean values (first column) of $u(x)$, the absolute errors (second column), and the 2 times standard deviations (third column) using MLP-DTEKI. }
\label{darcy_u_mlp}
\end{figure}

\begin{figure}[htbp]
\centering 
\begin{overpic}[width = 0.32\textwidth]{figures/darcy/kappa_HMC.pdf}
\put (-5, 40) {\rotatebox{90}{\footnotesize HMC}}
\end{overpic}
\begin{overpic}[width = 0.32\textwidth]{figures/darcy/kappa_error_HMC.pdf}
\end{overpic}
\begin{overpic}[width = 0.32\textwidth]{figures/darcy/kappa_std_HMC.pdf}
\end{overpic}
\begin{overpic}[width = 0.32\textwidth]{figures/darcy/kappa_EKI.pdf}
\put (-5, 40) {\rotatebox{90}{\footnotesize DTEKI}}
\end{overpic}
\begin{overpic}[width = 0.32\textwidth]{figures/darcy/kappa_error_EKI.pdf}
\end{overpic}
\begin{overpic}[width = 0.32\textwidth]{figures/darcy/kappa_std_EKI.pdf}
\end{overpic}
\begin{overpic}[width = 0.32\textwidth]{figures/darcy/kappa_subspace_EKI.pdf}
\put (-5, 40) {\rotatebox{90}{\footnotesize SDTEKI}}
\end{overpic}
\begin{overpic}[width = 0.32\textwidth]{figures/darcy/kappa_error_subspace_EKI.pdf}
\end{overpic}
\begin{overpic}[width = 0.32\textwidth]{figures/darcy/kappa_std_subspace_EKI.pdf}
\end{overpic}
\begin{overpic}[width = 0.32\textwidth]{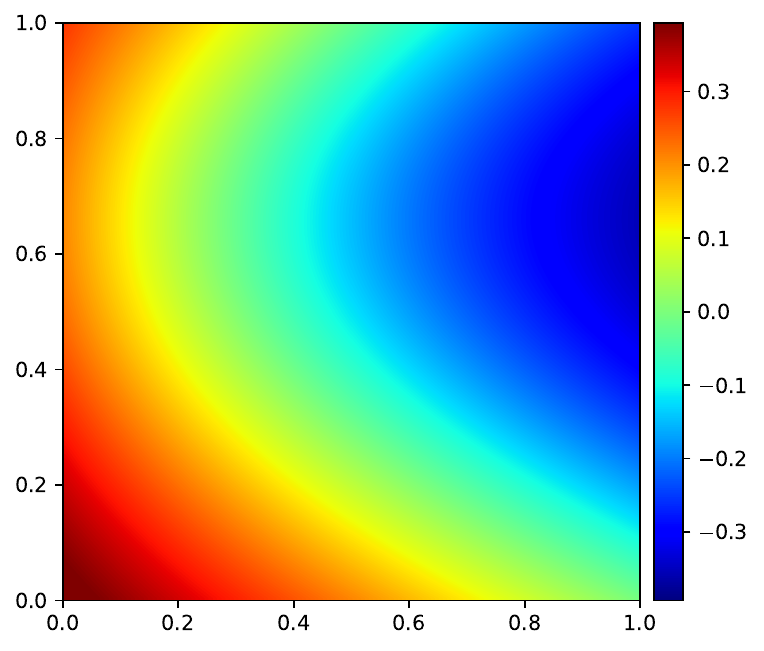}
\put (-5, 25) {\rotatebox{90}{\footnotesize MLP-DTEKI}}
\end{overpic}
\begin{overpic}[width = 0.32\textwidth]{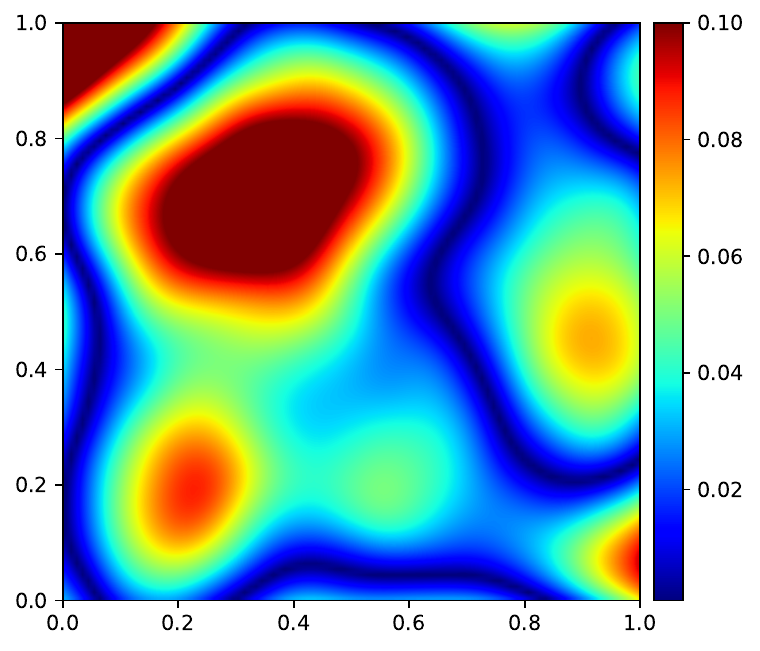}
\end{overpic}
\begin{overpic}[width = 0.32\textwidth]{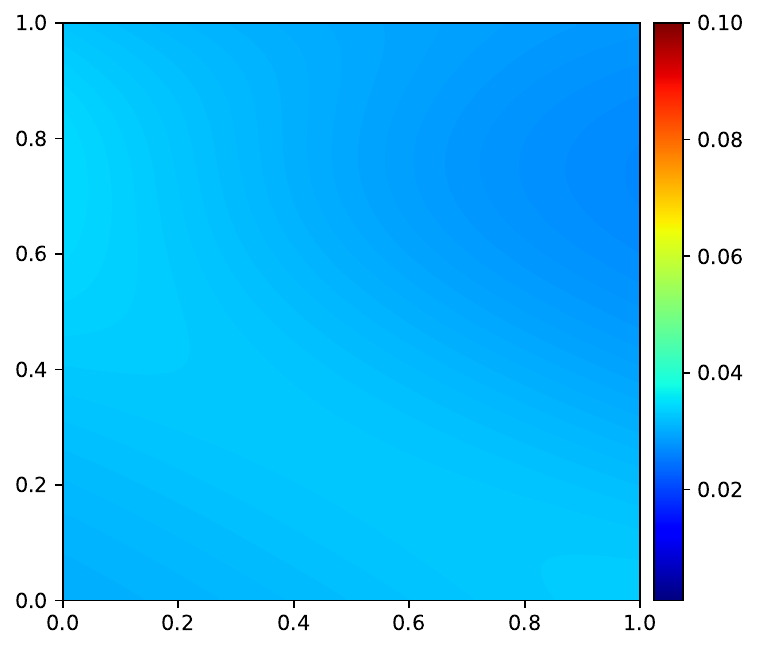}
\end{overpic}
\caption{Darcy equation: predicted mean values (first column) of $\kappa(x,\lambda)$, the absolute errors (second column), and the 2 times standard deviations (third column) using MLP-DTEKI.}
\label{darcy_kappa_mlp}
\end{figure}

For other problems, we summarize the performance of MLP-DTEKI in Table \ref{tab:performance_comparison}. It is clear that when using MLPs as surrogates in DTEKI, the $L_{2}$ errors obtained by DTEKI are much larger than those obtained by using cKANs as surrogates. To be clearer, we plot the predicted values of $u(x)$ and $f(x)$ as well as the predicted posterior distribution for the physical parameter in Figs.\ref{diffusion_mlp}-\ref{darcy_kappa_mlp}.

\begin{table}[htbp]
\centering
\begin{tabular}{@{}lcccc@{}}
\toprule
                 & \multicolumn{2}{c}{\textbf{MLP-DTEKI}} & \multicolumn{2}{c}{\textbf{DTEKI}} \\ \cmidrule(lr){2-3} \cmidrule(lr){4-5}
\textbf{Problem} & $e_{u}$        & $e_{\lambda}$        & $e_{u}$        & $e_{\lambda}$        \\ \midrule
Diffusion        & $99.98\%$     & $994.54\%$           & $11.55\%$      & $0.85\%$           \\
Nonlinear        & $104.82\%$     & $25.08\%$            & $9.54\%$     & $4.83\%$            \\
Darcy            & $2.49\%$       & $27.14\%$            & $1.28\%$       & $13.05\%$            \\ \bottomrule
\end{tabular}
\caption{The $L_{2}$ errors of the prediction for different problems obtained by two methods.}
\label{tab:performance_comparison}
\end{table}

For 1D problems, they involve multi-scale behavior and thus are hard to solve. Therefore, MLPs give nearly wrong predictions, without any tendencies when the data is sparse. While for the 2D Darcy problem, though the forward prediction of the solution is very accurate, the inverse prediction of the permeability field is not as good as cKANs, indicating that MLPs are hard to express the field. Moreover, for the 2D Darcy problem, the standard deviation of $u(x)$ is too large, providing nearly useless information to the error bound. In summary, the new method we propose using cKANs could not only improve the prediction accuracy but can scale well to provide reliable uncertainty estimates for large-scale inverse problems. 

\begin{figure}[t]
\centering 
\begin{overpic}[width = 0.32\textwidth]{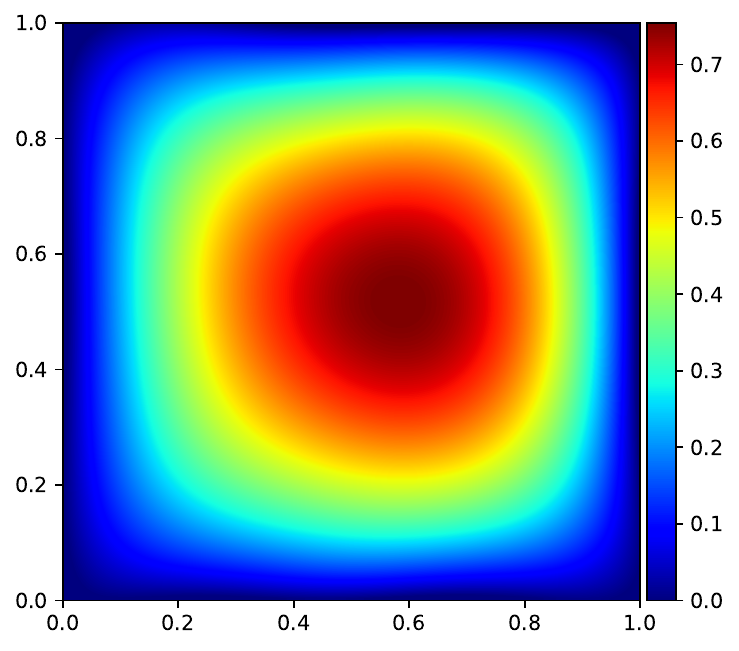}
\end{overpic}
\begin{overpic}[width = 0.32\textwidth]{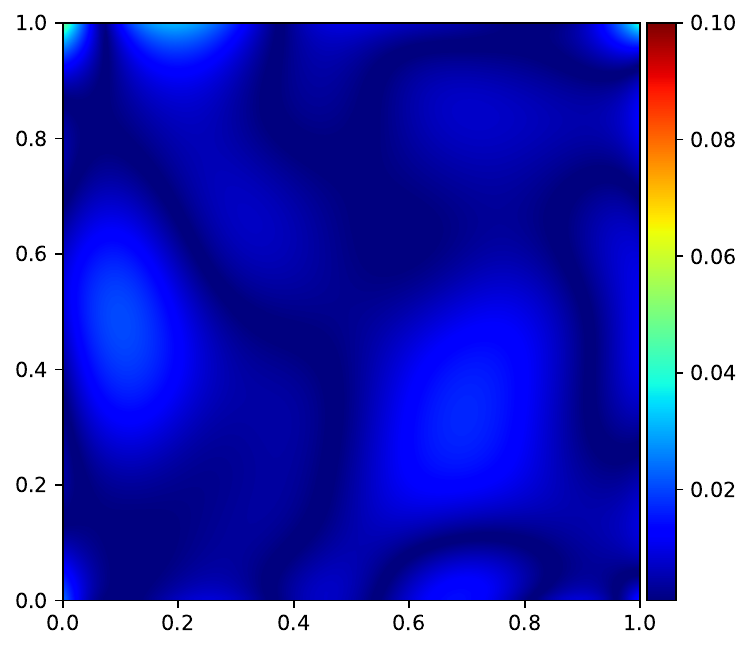}
\end{overpic}
\begin{overpic}[width = 0.32\textwidth]{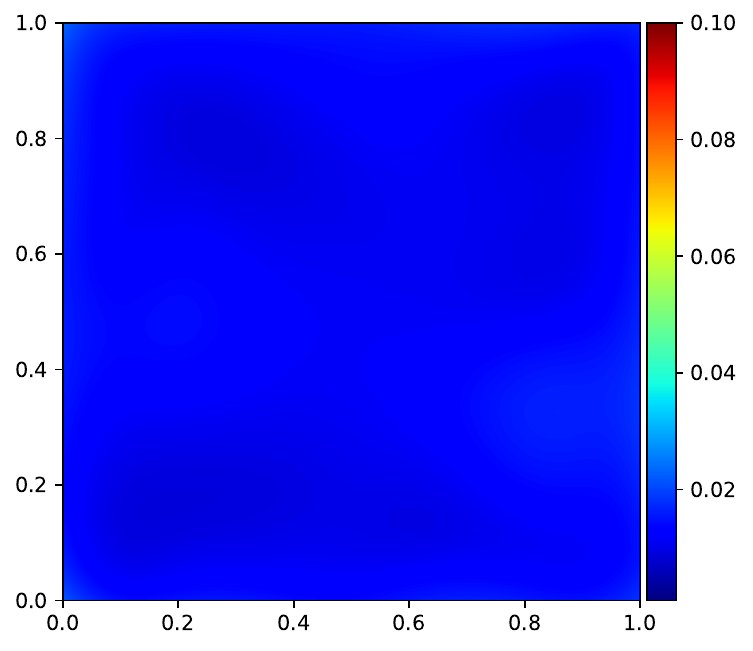}
\end{overpic}
\begin{overpic}[width = 0.32\textwidth]{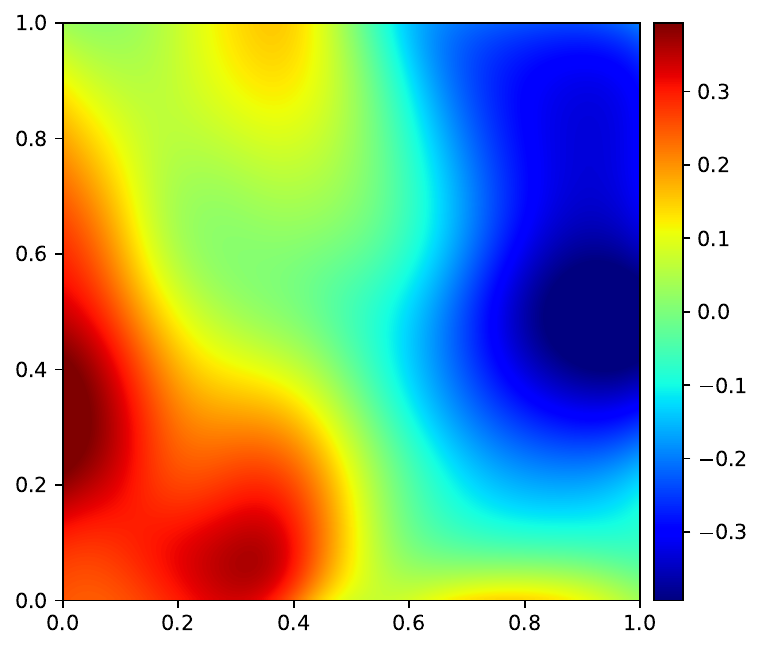}
\end{overpic}
\begin{overpic}[width = 0.32\textwidth]{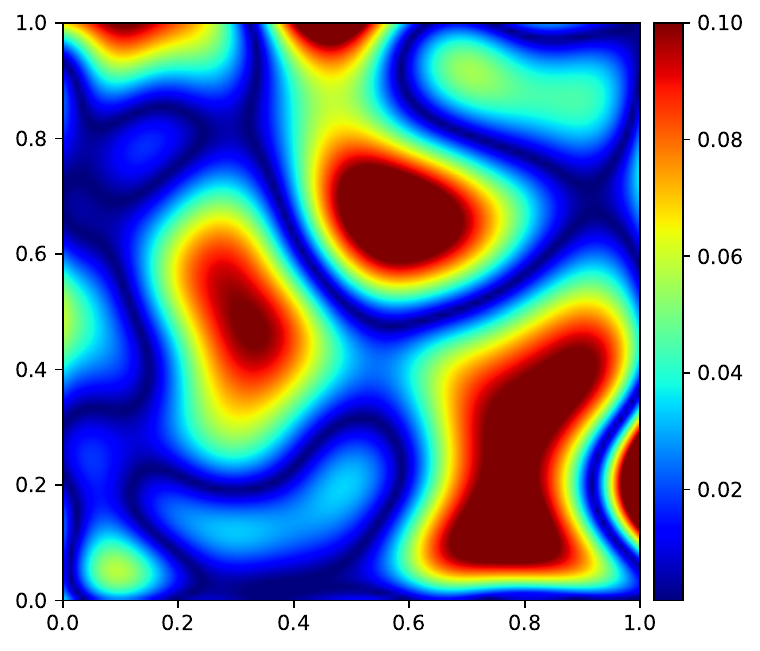}
\end{overpic}
\begin{overpic}[width = 0.32\textwidth]{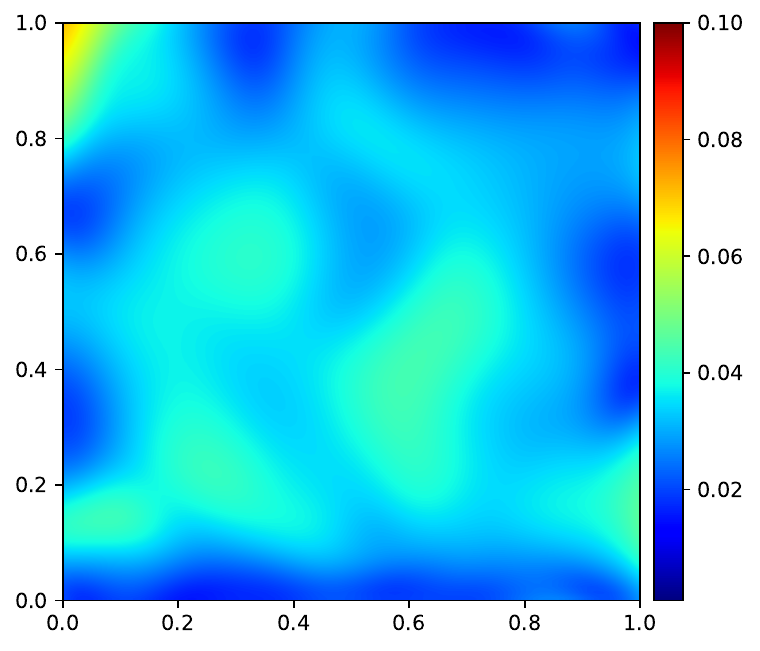}
\end{overpic}
\caption{Darcy equation: predicted solution field and permeability field (first column), the absolute errors (second column), and 2 times standard deviations (third column) using vanilla EKI without dropout and Tikhonov regularization.}
\label{vanilla_EKI}
\end{figure}

To test the importance of dropout and Tikhonov regularization compared to vanilla EKI, we test the 2D Darcy problem with simple EKI, the results are plotted in Fig. \ref{vanilla_EKI}. As depicted, without dropout and Tikhonov regularization, the performance will degrade, and also the standard deviation can not cover the absolute error. In detail, the $L_{2}$ errors for the solution field $u(x)$ and permeability field $\kappa(x, \lambda)$ are $1.58\%$ and $26.74\%$ respectively, which are much worse than results obtained by combining these two techniques. The reason behind this is that the network could possibly overfit with a large set of parameters and thus the estimation of the covariance matrices is not as accurate as adding these two techniques. In contrast, it is verified by combining dropout and Tikhonov regularization, it is reasonable for the performance to be improved. 

\bibliographystyle{siamplain}
\bibliography{references}
\end{document}